\documentclass[11pt, a4paper]{article}

\usepackage[cp1251]{inputenc}
\usepackage[russian]{babel}
\usepackage{amsmath, stmaryrd, latexsym, longtable} \usepackage{euscript}
\usepackage{mymatrx5}

\Linestrue\Autonumtrue
\usepackage{mathrsfs}
\usepackage{amssymb, amscd}
\usepackage[matrix, arrow, curve]{xy}
\begin{document}
\newcounter{bnomer} \newcounter{snomer}
\newcounter{bsnomer}
\setcounter{bnomer}{0}
\renewcommand{\thesnomer}{\thebnomer.\arabic{snomer}}
\renewcommand{\thebsnomer}{\thebnomer.\arabic{bsnomer}}
\renewcommand{\refname}{\begin{center}\large{\textbf{References}}\end{center}}
\renewcommand{\contentsname}{\begin{center}\large{\textbf{Table of contents}}\end{center}}

\setcounter{MaxMatrixCols}{14}

\newcommand{\sect}[1]{%
\setcounter{snomer}{0}\setcounter{bsnomer}{0}
\refstepcounter{bnomer}
\par\bigskip\begin{center}\large{\textbf{\arabic{bnomer}. {#1}}}\end{center}}
\newcommand{\fakesect}{%
\setcounter{snomer}{0}\setcounter{bsnomer}{0}
\refstepcounter{bnomer}}
\newcommand{\sst}[1]{%
\refstepcounter{bsnomer}
\par\bigskip\textbf{\arabic{bnomer}.\arabic{bsnomer}. {#1}.}}
\newcommand{\fakesst}{%
\refstepcounter{bsnomer}}
\newcommand{\defi}[1]{%
\refstepcounter{snomer}
\par\medskip\textbf{Definition \arabic{bnomer}.\arabic{snomer}. }{#1}\par\medskip}
\newcommand{\theo}[2]{%
\refstepcounter{snomer}
\par\textbf{Теорема \arabic{bnomer}.\arabic{snomer}. }{#2} {\emph{#1}}\hspace{\fill}$\square$\par}
\newcommand{\mtheop}[2]{%
\refstepcounter{snomer}
\par\textbf{Theorem \arabic{bnomer}.\arabic{snomer}. }{\emph{#1}}
\par\textsc{Proof}. {#2}\hspace{\fill}$\square$\par}
\newcommand{\mcorop}[2]{%
\refstepcounter{snomer}
\par\textbf{Corollary \arabic{bnomer}.\arabic{snomer}. }{\emph{#1}}
\par\textsc{Proof}. {#2}\hspace{\fill}$\square$\par}
\newcommand{\mtheo}[1]{%
\refstepcounter{snomer}
\par\medskip\textbf{Theorem \arabic{bnomer}.\arabic{snomer}. }{\emph{#1}}\par\medskip}
\newcommand{\mlemm}[1]{%
\refstepcounter{snomer}
\par\medskip\textbf{Lemma \arabic{bnomer}.\arabic{snomer}. }{\emph{#1}}\par\medskip}
\newcommand{\mprop}[1]{%
\refstepcounter{snomer}
\par\medskip\textbf{Proposition \arabic{bnomer}.\arabic{snomer}. }{\emph{#1}}\par\medskip}
\newcommand{\theobp}[2]{%
\refstepcounter{snomer}
\par\textbf{Теорема \arabic{bnomer}.\arabic{snomer}. }{#2} {\emph{#1}}\par}
\newcommand{\theop}[2]{%
\refstepcounter{snomer}
\par\textbf{Theorem \arabic{bnomer}.\arabic{snomer}. }{\emph{#1}}
\par\textsc{Proof}. {#2}\hspace{\fill}$\square$\par}
\newcommand{\theosp}[2]{%
\refstepcounter{snomer}
\par\textbf{Теорема \arabic{bnomer}.\arabic{snomer}. }{\emph{#1}}
\par\textbf{Схема доказательства}. {#2}\hspace{\fill}$\square$\par}
\newcommand{\exam}[1]{%
\refstepcounter{snomer}
\par\medskip\textbf{Example \arabic{bnomer}.\arabic{snomer}. }{#1}\par\medskip}
\newcommand{\deno}[1]{%
\refstepcounter{snomer}
\par\textbf{Definition \arabic{bnomer}.\arabic{snomer}. }{#1}\par}
\newcommand{\post}[1]{%
\refstepcounter{snomer}
\par\textbf{Предложение \arabic{bnomer}.\arabic{snomer}. }{\emph{#1}}\hspace{\fill}$\square$\par}
\newcommand{\postp}[2]{%
\refstepcounter{snomer}
\par\medskip\textbf{Proposition \arabic{bnomer}.\arabic{snomer}. }{\emph{#1}}%
\ifhmode\par\fi\textsc{Proof}. {#2}\hspace{\fill}$\square$\par\medskip}
\newcommand{\lemm}[1]{%
\refstepcounter{snomer}
\par\textbf{Lemma \arabic{bnomer}.\arabic{snomer}. }{\emph{#1}}\hspace{\fill}$\square$\par}
\newcommand{\lemmp}[2]{%
\refstepcounter{snomer}
\par\medskip\textbf{Lemma \arabic{bnomer}.\arabic{snomer}. }{\emph{#1}}
\par\textsc{Proof}. {#2}\hspace{\fill}$\square$\par\medskip}
\newcommand{\coro}[1]{%
\refstepcounter{snomer}
\par\textbf{Corollary \arabic{bnomer}.\arabic{snomer}. }{\emph{#1}}\hspace{\fill}$\square$\par}
\newcommand{\mcoro}[1]{%
\refstepcounter{snomer}
\par\textbf{Corollary \arabic{bnomer}.\arabic{snomer}. }{\emph{#1}}\par\medskip}
\newcommand{\corop}[2]{%
\refstepcounter{snomer}
\par\textbf{Corollary \arabic{bnomer}.\arabic{snomer}. }{\emph{#1}}
\par\textsc{Proof}. {#2}\hspace{\fill}$\square$\par}
\newcommand{\nota}[1]{%
\refstepcounter{snomer}
\par\medskip\textbf{Remark \arabic{bnomer}.\arabic{snomer}. }{#1}\par\medskip}
\newcommand{\propp}[2]{%
\refstepcounter{snomer}
\par\medskip\textbf{Proposition \arabic{bnomer}.\arabic{snomer}. }{\emph{#1}}
\par\textsc{Proof}. {#2}\hspace{\fill}$\square$\par\medskip}
\newcommand{\hypo}[1]{%
\refstepcounter{snomer}
\par\medskip\textbf{Conjecture \arabic{bnomer}.\arabic{snomer}. }{\emph{#1}}\par\medskip}
\newcommand{\prop}[1]{%
\refstepcounter{snomer}
\par\textbf{Proposition \arabic{bnomer}.\arabic{snomer}. }{\emph{#1}}\hspace{\fill}$\square$\par}

\def\sideremark#1{\ifvmode\leavevmode\fi\vadjust{\vbox to0pt{\vss
 \hbox to 0pt{\hskip\hsize\hskip1em
\vbox{\hsize1.6cm\tiny\raggedright\pretolerance10000 
 \noindent #1\hfill}\hss}\vbox to8pt{\vfil}\vss}}}

\newcommand\restr[2]{{
  \left.\kern-\nulldelimiterspace 
  #1 
  \right|_{#2} 
}}

\newcommand{\Ind}[3]{%
\mathrm{Ind}_{#1}^{#2}{#3}}
\newcommand{\Res}[3]{%
\mathrm{Res}_{#1}^{#2}{#3}}
\newcommand{\epsi}{\varepsilon}
\newcommand{\tri}{\triangleleft}
\newcommand{\Supp}[1]{%
\mathrm{Supp}(#1)}
\newcommand{\Sing}[1]{%
\mathrm{Sing}(#1)}

\newcommand{\Pic}[1]{%
\mathrm{Pic}\,#1}

\newcommand{\Hom}[1]{%
\mathrm{Hom}(#1)}

\newcommand{\reg}{\mathrm{reg}}
\newcommand{\sreg}{\mathrm{sreg}}
\newcommand{\codim}{\mathrm{codim}\,}
\newcommand{\chara}{\mathrm{char}\,}
\newcommand{\rk}{\mathrm{rk}\,}
\newcommand{\chr}{\mathrm{ch}\,}
\newcommand{\Ann}{\mathrm{Ann}\,}
\newcommand{\Aut}{\mathrm{Aut}\,}
\newcommand{\Lie}{\mathrm{Lie}\,}
\newcommand{\id}{\mathrm{id}}
\newcommand{\Ad}{\mathrm{Ad}}
\newcommand{\St}{\mathrm{St}}
\newcommand{\col}{\mathrm{col}}
\newcommand{\row}{\mathrm{row}}
\newcommand{\low}{\mathrm{low}}
\newcommand{\pho}{\hphantom{\quad}\vphantom{\mid}}
\newcommand{\fho}[1]{\vphantom{\mid}\setbox0\hbox{00}\hbox to \wd0{\hss\ensuremath{#1}\hss}}
\newcommand{\wt}{\widetilde}
\newcommand{\wh}{\widehat}
\newcommand{\ad}[1]{\mathrm{ad}_{#1}}
\newcommand{\Grr}[2]{\mathrm{Gr}(#1,#2)}
\newcommand{\Grb}[2]{\mathrm{Gr}^{\beta}(#1,#2)}
\newcommand{\Grbn}[3]{\mathrm{Gr}^{#3}(#1,#2)}
\newcommand{\Grbi}[3]{\mathrm{Gr}_{#3}^{\beta}(#1,#2)}
\newcommand{\Gra}[1]{\mathrm{Gr}(#1)}
\newcommand{\Grad}[1]{\mathrm{Gr}'(#1)}
\newcommand{\Grab}[1]{\mathrm{Gr}^{\beta}(#1)}
\newcommand{\tr}{\mathrm{tr}\,}
\newcommand{\Cent}{\mathrm{Cent}\,}
\newcommand{\GL}{\mathrm{GL}}
\newcommand{\Or}{\mathrm{O}}
\newcommand{\SO}{\mathrm{SO}}
\newcommand{\Spin}{\mathrm{Spin}}
\newcommand{\SU}{\mathrm{SU}}
\newcommand{\SL}{\mathrm{SL}}
\newcommand{\Sp}{\mathrm{Sp}}
\newcommand{\Mat}{\mathrm{Mat}}
\newcommand{\Pf}{\mathrm{Pf}}
\newcommand{\Prim}{\mathrm{Prim}\,}
\newcommand{\Ker}{\mathrm{Ker}\,}
\newcommand{\Stab}[2]{\mathrm{Stab}_{#1}{#2}}
\newcommand{\ilm}{\varinjlim}
\newcommand{\plm}{\varprojlim}
\newcommand{\soc}[1]{\mathrm{soc}^{#1}\,}
\newcommand{\gee}{\geqslant}
\newcommand{\lee}{\leqslant}

\newcommand{\vfi}{\varphi}
\newcommand{\teta}{\vartheta}
\newcommand{\Bfi}{\Phi}
\newcommand{\uto}{\uparrow}
\newcommand{\Fp}{\mathbb{F}}
\newcommand{\Pp}{\mathbb{P}}
\newcommand{\Np}{\mathbb{N}}
\newcommand{\Rp}{\mathbb{R}}
\newcommand{\Zp}{\mathbb{Z}}
\newcommand{\Cp}{\mathbb{C}}
\newcommand{\Qp}{\mathbb{Q}}
\newcommand{\Hp}{\mathbb{H}}
\newcommand{\Ap}{\mathbb{A}}
\newcommand{\ut}{\mathfrak{u}}
\newcommand{\at}{\mathfrak{a}}
\newcommand{\nt}{\mathfrak{n}}
\newcommand{\mt}{\mathfrak{m}}
\newcommand{\htt}{\mathfrak{h}}
\newcommand{\ttt}{\mathfrak{t}}
\newcommand{\spt}{\mathfrak{sp}}
\newcommand{\slt}{\mathfrak{sl}}
\newcommand{\sot}{\mathfrak{so}}
\newcommand{\ot}{\mathfrak{o}}
\newcommand{\glt}{\mathfrak{gl}}
\newcommand{\rt}{\mathfrak{r}}
\newcommand{\zt}{\mathfrak{z}}
\newcommand{\rad}{\mathfrak{rad}}
\newcommand{\bt}{\mathfrak{b}}
\newcommand{\gt}{\mathfrak{g}}
\newcommand{\vt}{\mathfrak{v}}
\newcommand{\pt}{\mathfrak{p}}
\newcommand{\Xt}{\mathfrak{X}}
\newcommand{\Po}{\mathcal{P}}
\newcommand{\Uo}{\EuScript{U}}
\newcommand{\Fo}{\EuScript{F}}
\newcommand{\Do}{\mathcal{D}}
\newcommand{\Eo}{\EuScript{E}}
\newcommand{\Vo}{\EuScript{V}}
\newcommand{\Iu}{\mathcal{I}}
\newcommand{\Mo}{\mathcal{M}}
\newcommand{\Nu}{\mathcal{N}}
\newcommand{\Ro}{\mathcal{R}}
\newcommand{\Co}{\mathcal{C}}
\newcommand{\Go}{\mathcal{G}}
\newcommand{\Lo}{\mathcal{L}}
\newcommand{\So}{\mathcal{S}}
\newcommand{\Ou}{\mathcal{O}}
\newcommand{\Uu}{\mathcal{U}}
\newcommand{\Au}{\mathcal{A}}
\newcommand{\Du}{\mathcal{D}}
\newcommand{\Vu}{\mathcal{V}}
\newcommand{\Bu}{\mathcal{B}}
\newcommand{\Sy}{\mathcal{Z}}
\newcommand{\Sb}{\mathcal{F}}
\newcommand{\Gr}{\mathcal{G}}
\newcommand{\Fl}{\EuScript{F}\ell}
\newcommand{\rtc}[1]{C_{#1}^{\mathrm{red}}}

\title{Ind-varieties of generalized flags: a survey of results\footnote{The first author was partially supported by RFBR grants no. 14--01--97017 and 16--01--00154, and by the Ministry of Science and Education of the Russian Federation, project no. 204. A part of this work was done during the stay of the first author at Jacobs University Bremen; the first author thanks this institution for its hospitality. The first and the second authors were partially supported by DFG, grant no. PE 980/6--1.}}
\date{}
\author{Mikhail V. Ignatyev\and Ivan Penkov}
\maketitle

\begin{center}
\begin{tabular}{p{15cm}}
\small{$\hphantom{1111}$\textsc{Abstract}. This is a review of results on the structure of the homogeneous ind-varieties $G/P$ of the ind-groups $G=\GL_{\infty}(\Cp)$, $\SL_{\infty}(\Cp)$, $\SO_{\infty}(\Cp)$, $\Sp_{\infty}(\Cp)$, subject to the condition that $G/P$ is a inductive limit of compact homogeneous spaces $G_n/P_n$. In this case the subgroup $P\subset G$ is a splitting parabolic subgroup of $G$, and the ind-variety $G/P$ admits a ``flag realization''. Instead of ordinary flags, one considers generalized flags which are, generally infinite, chains $\Co$ of subspaces in the natural representation $V$ of $G$ which satisfy a certain condition: roughly speaking, for each nonzero vector~$v$ of~$V$ there must be a largest space in~$\Co$ which does not contain $v$, and a smallest space in~$\Co$ which contains $v$.

$\hphantom{1111}$We start with a review of the construction of the ind-varieties of generalized flags, and then show that these ind-varieties are homogeneous ind-spaces of the form $G/P$ for splitting parabolic ind-subgroups $P\subset G$. We also briefly review the characterization of more general, i.e. non-splitting, parabolic ind-subgroups in terms of generalized flags. In the special case of an ind-grassmannian $X$, we give a purely algebraic-geometric construction of $X$. Further topics discussed are the Bott--Borel--Weil Theorem for ind-varieties of generalized flags, finite-rank vector bundles on ind-varieties of generalized flags, the theory of Schubert decomposition of $G/P$ for arbitrary splitting parabolic ind-subgroups $P\subset G$, as well as the orbits of real forms on $G/P$ for $G=\SL_{\infty}(\Cp)$.}
\end{tabular}
\end{center}


\maketitle

\tableofcontents

\section{Introduction}\label{sect:introduction}\fakesect

Flag varieties play a fundamental role in geometry. Projective spaces and grassmannians go back to the roots of modern geometry and appear in almost every aspect of global geometric studies. Full flag varieties (varieties of maximal flags) also play a significant role in geometry in general, but are most often associated with geometric representation theory. This is because they are universal compact homogeneous spaces of the group $\GL_n(\Cp)$, and almost any branch of the representation theory of an algebraic group $G$ has deep results connected with the geometry of the flag variety $G/B$. The Bott--Borel--Weil Theorem is a basic example in this direction. Another example is Beilinson--Bernstein's Localization Theorem, and there is a long list of further striking results.

When one thinks what an ``infinite-dimensional flag variety'' should be, one realizes that there are many possibilities. One option is to study the homogeneous ind-spaces $G/B$ for Kac--Moody groups $G$. These ind-varieties play a prominent role in representation theory and geometry since 1980's, see for instance \cite{Mathieu1}, \cite{PressleySegal1} and \cite{Kumar1}.

There is another option which we consider in the present paper. Our main topic are homogeneous ind-spaces of the locally linear ind-groups $\GL_{\infty}(\Cp)$, $\SL_{\infty}(\Cp)$, $\SO_{\infty}(\Cp)$ and $\Sp_{\infty}(\Cp)$. These ``infinite-dimensional flag varieties'' are smooth locally projective ind-varieties (in particular, ind-schemes), and their points are certain chains of subspaces in a given countable-dimensional complex vector space $V\cong\Cp^{\infty}$. Moreover, our ``infinite-dimensional flag ind-varieties'' are exhausted by of usual finite-dimensional flag varieties. What is not obvious is which type of chains of subspaces in $V$ one has to consider so that there is a bijection between such chains and Borel (or, more generally, parabolic) subgroups of the group $\GL_{\infty}(\Cp)$. The answer to this question, obtained in \cite{DimitrovPenkov1}, leads to the definition of a generalized flag.

A generalized flag is a chain $\Fo$ of subspaces which satisfies the condition that every vector $v\in V$ determines a unique pair $F_v'\subset F_v''$ of subspaces from $\Fo$ so that $F_v'$ is the immediate predecessor of $F_v''$ and $v\in F_v''\setminus F_v'$ (the precise definition see in Subsection~\ref{sst:gen_flags}). What is important is that it is not required that every space $F\in\Fo$ have an immediate predecessor \emph{and} immediate successor; instead it is required that only an immediate predecessor \emph{or} an immediate successor exist. The appearance of these somewhat complicated linear orders on the chains $\Fo$ is of course related with the fact that the Borel subgroups of $\GL_{\infty}(\Cp)$ containing a fixed splitting maximal torus are in one-to-one correspondence with arbitrary linear orders on a countable set \cite{DimitrovPenkov3}.

Another way of seeing how generalized flags appear is by trying to extend to infinity a finite flag in a finite-dimensional space by the following infinite process: at each stage one adds another flag to a part of the flag obtained at the previous stage so that every time the length of the flag grows up at most by one. The variety of choices for the places where the finite flags spread to become longer leads to the fact that the output of this procedure is a generalized flag; the precise mathematical expression of this ``spreading procedure'' is given by formula (\ref{formula:iota_n_another}) in Subsection~\ref{sst:gen_flags}.

As a result, we arrive to a nice characterization of the ind-varieties of generalized flags as ind-varieties exhausted by usual flag varieties. We need to mention also the issue of commensurability of generalized flags: two points on the same ind-variety of generalized flags correspond to generalized flags which are ``close to each other'', i.e., $E$-commensurable in technical terms, see Subsection~\ref{sst:gen_flags}. The idea of commensurability goes back to A. Pressley and G. Segal.

We believe that the interplay between the ind-varieties of generalized flags and the representation theory of the ind-groups $\GL_{\infty}(\Cp)$, $\SL_{\infty}(\Cp)$, $\SO_{\infty}(\Cp)$, $\Sp_{\infty}(\Cp)$ will be at least as rich as in the finite-dimensional case. However, both theories are very far from having reached the maturity of the finite-dimensional theory, so a great deal of work lies ahead. Our moderate goal is to provide a first survey on results which should serve as an invitation for further research.

A brief description of the contents of the paper is as follows. In Section~\ref{sect:definitions} we give the definitions of and ind-variety and an ind-group and consider our main examples of ind-groups. We define generalized flags (and isotropic generalized flags) and show that the generalized flags which are $E$-commensurable to a given generalized flag $\Fo$ form an ind-variety. In Section~\ref{sect:ind_grassmannians} we show that ind-grassmannians can be defined in a purely algebraic-geometric terms as inductive limits of linear embeddings of finite-dimensional grassmannians. In Section~\ref{sect:homo} we discuss Cartan, Borel and parabolic subgroups of classical ind-groups and show that ind-varieties of generalized flags are homogeneous spaces $G/P$ where $G=\GL_{\infty}$, $\SL_{\infty}(\Cp)$ (or $G=\SO_{\infty}(\Cp)$, $\Sp_{\infty}(\Cp)$ in the case of isotropic generalized flags) and $P$ is a splitting parabolic ind-subgroup. We also briefly discuss general, i.e. non-splitting, Borel and parabolic subgroups.

In Section~\ref{sect:Schubert_calculus} we provide the Bruhat decomposition of a classical ind-group as well as the Schubert decomposition of its ind-varieties of generalized flags. In particular, we present a criterion for the smoothness of a Schubert subvariety. Section~\ref{sect:alg_geom} is devoted to a generalization of the Bott--Borel--Weil theorem and the Barth--Van de Ven--Tyurin--Sato theorem to the case of generalized flags. Finally, in Section~\ref{sect:orbits_real_forms} we extend well-known results of J.A. Wolf on orbits of real forms of semisimple Lie groups to ind-varieties of generalized flags for $G=\SL_{\infty}(\Cp)$.

\textsc{Acknowledgements}. The first author has been supported in part by the Russian Foundation for Basic Research through grants no. 14--01--97017 and 16--01--00154, and by the Ministry of Science and Education of the Russian Federation, project no. 204. A part of this work was done during the stay of the first author at Jacobs University Bremen, and the first author thanks this institution for its hospitality. Both authors have been supported in part by DFG grant no. PE 980/6--1.

\section{Definitions and examples}\label{sect:definitions}\fakesect

In this section we give precise definitions, provide the main examples needed for the sequel, and establish some first properties of generalized flags. The material below is taken from \cite{DimitrovPenkov1} and \cite{FressPenkov1}.

\subsection{Ind-varieties and ind-groups}\fakesst \label{sst:ind_varieties} All varieties and algebraic groups are defined over the field of complex numbers $\Cp$.

\defi{An \emph{ind-variety} is the inductive limit $X=\ilm X_n$ of a chain of morphisms of algebraic varieties
\begin{equation}
\label{formula:ind_variety}
X_1\stackrel{\vfi_1}\to X_2\stackrel{\vfi_2}\to\ldots\stackrel{\vfi_{n-1}}{\to}
X_n\stackrel{\vfi_n}\to X_{n+1}\stackrel{\vfi_{n+1}}\to\ldots.
\end{equation}}

Obviously, the inductive limit of a chain (\ref{formula:ind_variety}) does not change if we replace the sequence $\{X_n\}_{n\geq1}$ by a subsequence $\{X_{i_n}\}_{n\geq 1}$, and the morphisms $\vfi_n$ by the compositions $$\wt\vfi_{i_n}=\vfi_{i_{n+1}-1}\circ\ldots\circ\vfi_{i_n+1}\circ\vfi_{i_n}.$$ Let $Y$ be a second ind-variety obtained as the inductive limit of a chain
\begin{equation*}
Y_1\stackrel{\psi_1}\to Y_2\stackrel{\psi_2}\to\ldots\stackrel{\psi_{n-1}}\to
Y_n\stackrel{\psi_n}\to Y_{n+1}\stackrel{\psi_{n+1}}\to\ldots.
\end{equation*}
A \emph{morphism of ind-varieties} $f\colon X\to Y$ is a map from $\ilm X_n$ to $\ilm Y_n$ induced by a collection of morphisms of algebraic varieties $$\{f_n\colon X_n\to Y_{i_n}\}_{n\geq1}$$ such that $$\wt\psi_{i_n}\circ f_n=f_{n+1}\circ\vfi_n$$ for all $n\geq1$. The \emph{identity morphism} $\id_X$ is a morphism that induces the identity map $X\to X$. As usual, a morphism $f\colon X\to Y$ is called an \emph{isomorphism} if there exists a morphism $g\colon Y\to X$ such that $$f\circ g=\id_Y,~g\circ f=\id_X.$$

Any ind-variety $X$ is endowed with a topology by declaring a subset $U\subset X$ \emph{open} if its inverse image by the natural map $X_m\to\ilm X_n$ is open for all $m$. Clearly, any open (or also closed or locally closed) subset $Z$ of $X$ has a structure of ind-variety induced by the ind-variety structure on $X$. We call $Z$ an \emph{ind-subvariety} of $X$. In what follows we only consider chains (\ref{formula:ind_variety}) where the morphisms $\vfi_n$ are embeddings, so that we can write $X=\bigcup_{n\geq 1} X_n$. Then the sequence $\{X_n\}_{n\geq 1}$ is called an \emph{exhaustion} of the ind-variety $X$.

Next, we introduce the notion of a smooth point of an ind-variety. Let $x\in X$, so that $x\in X_n$ for $n$ large enough. Let $\mt_{n,x}\subset \mathcal{O}_{X_n,x}$ be the maximal ideal of
the localization at $x$ of $\mathcal{O}_{X_n}$. For each $k\geq1$ there is an epimorphism
\begin{equation*}
\alpha_{n,k}\colon S^k(\mt_{n,x}/\mt_{n,x}^2)\to\mt_{n,x}^k/\mt_{n,x}^{k+1}.
\end{equation*}
Note that the point $x$ is smooth in $X_n$ if and only if $\alpha_{n,k}$ is an isomorphism for all $k$. By taking the projective limit, we obtain a map
\begin{equation*}
\wh\alpha_k=\plm\alpha_{n,k}\colon\plm S^k(\mt_{n,x}/\mt_{n,x}^2)\to\plm\mt_{n,x}^k/\mt_{n,x}^{k+1},
\end{equation*}
which is an epimorphism for all $k$. We say that $x$ is a \emph{smooth point} of $X$ if $\wh\alpha_k$ is an isomorphism for all~$k$. We say that $x$ is a \emph{singular point} otherwise. The notion of smoothness of a point is independent of the choice of exhaustion $\{X_n\}_{n\geq1}$ of $X$. We say that the ind-variety $X$ is \emph{smooth} if every point $x\in X$ is smooth.

\exam{i) Suppose that each variety $X_n$ in the chain (\ref{formula:ind_variety}) is an affine space, every image $\vfi_n(X_n)$ is an affine subspace of $X_{n+1}$, and $$\lim_{n\to\infty}\dim X_n=\infty.$$ Then, up to isomorphism, $X=\ilm X_n$ is independent of the choice of $\{X_n,\vfi_n\}_{n\geq1}$ with these properties. We write $X=\Ap^\infty$ and call it the \emph{infinite-dimensional affine space}.
In particular, $\Ap^\infty$ admits the exhaustion $$\Ap^\infty=\bigcup_{n\geq1}\Ap^n,$$ where $\mathbb{A}^n$ stands for the $n$-dimensional affine space. Clearly, $\Ap^\infty$ is a smooth ind-variety.

ii) Now, suppose that each $X_n$ is a projective space, every image $\vfi_n(X_n)$ is a projective subspace of $X_{n+1}$, and $$\lim_{n\to\infty}\dim X_n=\infty.$$ Then $X=\ilm X_n$ is independent of the choice of $\{X_n,\vfi_n\}_{n\geq1}$ with these properties. We write $$X=\Pp^\infty=\bigcup_{n\geq1}\Pp^n$$ and call it the \emph{infinite-dimensional projective space}. It is also a smooth ind-variety.}

In general, an ind-group is a group object in the category of ind-varieties. However, in this paper we consider only locally linear algebraic ind-groups.

\defi{A \emph{locally linear algebraic ind-group} is an ind-variety $G=\bigcup G_n$ such that all $G_n$ are linear algebraic groups and the embeddings are group homomorphisms. In what follows we write \emph{ind-group} for brevity. Clearly, $G$ is a group. By an \emph{ind-subgroup} of $G$ we understand a closed subgroup of $G$. A \emph{morphism of ind-groups} $f\colon G\to H$ is a group homomorphism which is also a morphism of ind-varieties.}

We now introduce some ind-groups which play a central role in this review. The notation of the following example will be used throughout the paper.

\exam{i) Denote \label{exam:ind_groups}by $V$ a countable-dimensional complex vector space with fixed basis $E$. We fix an order on ${E}$ via the ordered set $\Zp_{>0}$, i.e., $$E=\{{e}_1,~{e}_2,~\ldots\}.$$ Let $V_*$ denote the span of the dual system $${E}^*=\{{e}_1^*,~{e}_2^*,~\ldots\}.$$ By definition, the \emph{finitary general linear group} $\GL(V,E)$ is the group of invertible $\Cp$-linear trans\-for\-ma\-tions on $V$ that keep fixed all but finitely many elements of ${E}$. It is not difficult, but important, to verify that $\GL(V,{E})$ depends only on the pair $(V,V_*)$ and not on $E$. We say that a basis $E'$ of $V$ is $G$-\emph{eligible} if $$G(V, E)=G(V,E').$$ Next, it is clear that any operator from $\GL(V,{E})$ has a well-defined determinant. By $\SL(V,E)$ we denote the subgroup of $\GL(V,E)$ of all operators with determinant $1$. We call this subgroup the \emph{finitary special linear group}.

Next, let us express the basis ${E}$ as a union $${E}=\bigcup {E}_n$$ of nested finite
subsets. Then $V$ is exhausted by the finite-dimensional subspaces
$V_n=\langle {E}_n\rangle_{\Cp}$, where $\langle\cdot\rangle_{\Cp}$ stands for linear span over $\Cp$. Equivalently, we can write $V=\varinjlim V_n$. To each linear operator $\vfi\colon V_n\to V_n$ one can assign an operator $$\wh\vfi\colon V_{n+1}\to V_{n+1}$$ such that $\wh\vfi(x)=\vfi(x)\text{ for }x\in V_n,~\wh\vfi(e)=e\text{ for }e\in E\setminus E_n$. This gives embeddings $$\GL(V_n)\hookrightarrow\GL(V_{n+1}),~\SL(V_n)\hookrightarrow\SL(V_{n+1}),$$ so that $$\GL(V,{E})=\varinjlim\GL(V_n),~\SL(V,{E})=\varinjlim\SL(V_n)$$ are ind-groups. Sometimes we write $\GL_{\infty}(\Cp)$ and $\SL_{\infty}(\Cp)$ instead of $\GL(V,{E})$ and $\SL(V,{E})$ res\-pec\-ti\-vely, then we keep in mind the above exhaustion of $V$ by the finite-dimensional subspaces $V_n$.

ii) Suppose now that $V$ is endowed with a nondegenerate symmetric or skew-symmetric bilinear
form $\beta$. We assume that the restriction~$\beta_n$ of the form $\beta$ to $V_n$ is nondegenerate for all $n$ and that ${\beta(e,V_n)=0}$ for $e\in E\setminus E_n$. Let the \emph{finitary orthogonal group} $\Or(V,{E},\beta)$ and the \emph{finitary symplectic group} $\Sp(V,{E},\beta)$ be the respective subgroups of $\GL(V,{E})$ consisting of all invertible operators pre\-ser\-ving the form $\beta$ in the cases when $\beta$ is symmetric or skew-symmetric. Put also $$\SO(V,{E},\beta)=\Or(V,{E},\beta)\cap\SL(V,{E}).$$ If a linear operator $\vfi$ on $V_n$ preserves $\beta_n$, then the linear operator $\wh\vfi$ on $V_{n+1}$ preserves $\beta_{n+1}$. Hence we can express our groups as inductive limits:
\begin{equation*}\predisplaypenalty=0
\Or(V,{E},\beta) = \varinjlim\Or(V_n,\beta_n),~\SO(V,{E},\beta) = \varinjlim\SO(V_n,\beta_n),~\Sp(V,{E},\beta) = \varinjlim\Sp(V_n,\beta).
\end{equation*}
Thus, they are ind-groups. Again, when we write $\Or_{\infty}(\Cp)$, $\SO_{\infty}(\Cp)$ or $\Sp_{\infty}(\Cp)$,
we have in mind the fixed exhaustion of $V$ by finite-dimensional $\beta$-nonsingular subspaces.}

Throughout Sections~\ref{sect:definitions}--\ref{sect:Schubert_calculus} we denote by $G$ one of the ind-groups $\GL_{\infty}(\Cp)$, $\SL_{\infty}(\Cp)$, $\SO_{\infty}(\Cp)$,\break $\Sp_{\infty}(\Cp)$. For example, assume that $G=\SL_{\infty}(\Cp)$. Let $H$ be the subgroup of elements $g\in G$ which are diagonal in the basis ${E}$. Then $H$ is an ind-subgroup of $G$ called a \emph{splitting Cartan subgroup}. An ind-subgroup $B\subset G$ which contains $H$ is called a \emph{splitting Borel subgroup} if it is locally solvable (i.e., every finite-dimensional ind-subgroup of $B$ is solvable) and is maximal with this property. An ind-subgroup which contains such a splitting Borel subgroup $B$ is called a \emph{splitting parabolic subgroup}. Equivalently, an ind-subgroup $P$ of $G$ containing $H$ is a splitting parabolic subgroup of $G$ if and only if $P\cap G_n$ is a parabolic subgroup of $G_n$ for all $n\geq1$, where $G=\bigcup_{n\geq 1}G_n$ is the natural exhaustion of $G$ mentioned above. The quotient $$G/P=\bigcup\nolimits_{n\geq 1}G_n/(P\cap G_n)$$ is a \emph{locally projective} ind-variety (i.e., an ind-variety exhausted by projective varieties); note however that $G/P$ is in general not a \emph{projective} ind-variety, i.e. $G/P$ is not isomorphic to a closed ind-subvariety of $\Pp^\infty$: see Theorem~\ref{theo:projectivity_via_BBW} below. Cartan, Borel and parabolic subgroups of an arbitrary classical ind-group $G$ are discussed in more detail in Section~\ref{sect:homo} below, where we also characterize $G/P$ as an ind-variety of generalized flags.

\subsection{Generalized flags}\fakesst In this subsection \label{sst:gen_flags}we introduce a key notion, namely that of a generalized flag. The obvious notion of an infinite flag (possibly, one-sided or two-sided, see its definition below) is not sufficient for describing the locally projective homogeneous ind-spaces of the classical ind-groups. This notion must be replaced by the somewhat intricate notion of a generalized flag which we now introduce.

Recall that $V$ is a countable-dimensional complex vector space. A~\emph{chain} of subspaces in $V$ is a set $\Co$ of distinct subspaces of $V$ such that if $F$, $F'\in\Co$, then either $F\subset F'$ or $F'\subset F$. Given a chain $\Co$ of subspaces in $V$, we write $\Co'$ (respectively,~$\Co''$) for the subchain of $\Co$ of all $F\in\Co$ with an immediate successor (respectively, an immediate predecessor). Also, we write $\Co^{\dag}$ for the set of all pairs $(F',~F'')$ such that $F''\in\Co''$ is the immediate successor of $F'\in\Co'$.

\defi{A \emph{generalized flag} is a chain $\Fo$ of subspaces in $V$ with the property $\Fo=\Fo'\cup\Fo''$ and such that $$V\setminus\{0\}=\bigcup_{(F',~F'')\in\Fo^{\dag}}F''\setminus F'.$$}

Note that each nonzero vector $v\in V$ determines a unique pair $(F_v',~F_v'')\in\Fo^{\dag}$ for which $v\in F_v''\setminus F_v'$. If $\Fo$ is a generalized flag, then each of the chains $\Fo'$ and $\Fo''$ determines $\Fo$. Indeed, if $(F',F'')\in\Fo^{\dag}$ then
$$F'=\bigcup_{G''\in\Fo'',~G''\subsetneq F''}G'',~F''=\bigcap_{G'\in\Fo',~G'\supsetneq F'}G'$$
(see \cite[Proposition 3.2]{DimitrovPenkov1}). A generalized flag $\Fo$ is called \emph{maximal} if $\Fo$ is not properly contained in another generalized flag. This is equivalent to the condition that $\dim F_v''=\dim F_v'+1$ for all nonzero vectors $v\in V$. A maximal generalized flag is not necessarily a maximal chain of subspaces in $V$, see Example~\ref{exam:gen_flags} v) below. Every generalized flag is contained in some maximal generalized flag.

By definition, a generalized flag $\Fo$ is a \emph{flag} if $\Fo$ is isomorphic as a linearly ordered set to a subset of $\Zp$ (with the natural order).

Given a generalized flag $\Fo$, we fix a linearly ordered set $(\Au,\preceq)$ and an isomorphism of ordered sets $$\Au\to\Fo^{\dag},~\alpha\mapsto(F_{\alpha}',F_{\alpha}'')$$ so that $\Fo$ can be written as $$\Fo=\{F_{\alpha}',F_{\alpha}'',~\alpha\in\Au\}.$$ We will write $\alpha\prec\beta$ if $\alpha\preceq\beta$ and $\alpha\neq\beta$ for $\alpha,\beta\in\Au$.

As above, fix a basis ${E}$ of $V$. We say that a generalized flag $\Fo$ is \emph{compatible} with ${E}$ if there exists a (necessarily surjective) map $\sigma\colon{E}\to\Au$ such that every pair $(F_{\alpha}',F_{\alpha}'')\in\Fo^{\dag}$ has the form
\begin{equation}
F_{\alpha}'=\langle e\in{E}\mid \sigma(e)\prec\alpha\rangle_{\Cp},~F_{\alpha}''=\langle e\in{E}\mid \sigma(e)\preceq\alpha\rangle_{\Cp}.\label{formula:F'_F''}
\end{equation}
By \cite[Proposition 4.1]{DimitrovPenkov1}, every generalized flag admits a com\-pa\-tible basis. Below we recall the proof of this fact.

\propp{Every generalized flag $\Fo$ in $V$ admits\label{prop:comp_basis} a compatible basis.}{First, consider the case when $\Fo$ is maximal. Let $$D=\{d_1,~d_2,~\ldots\}$$ be any basis of $V$. We use induction to construct a basis $L=\{l_1,~l_2,~\ldots\}$ such that $$\langle l_1,~\ldots,~l_n\rangle_{\Cp}=\langle d_1,~\ldots,~d_n\rangle_{\Cp}$$ for all~$n$ and the subspaces $F_{l_n}'$ are pairwise distinct.

Let $l_1=d_1$. Assume that the vectors $l_1$, $\ldots$, $l_n$ have already been constructed. Denote by $W$ the affine subspace of $V$ of the form $$d_{n+1}+\langle l_1,~\ldots,~l_n\rangle_{\Cp}.$$ We claim that there exists $l\in W$ such that $F'_l$ does not coincide with any of the spaces $F'_{l_1},~\ldots,~F'_{l_n}$. Indeed, assume the contrary. Then $W$ is contained in the union $\bigcup_{i=1}^n F_{l_i}''$, hence there exists $k$ for which $W\subset F_{l_k}''$. If $i_1,~\ldots,~i_n$ is the permutation of $1,~\ldots,~n$ such that $F_{l_{i_1}}'\subsetneq\ldots\subsetneq F_{l_{i_n}}'$, then $F_{l_{i_1}}''\subsetneq\ldots\subsetneq F_{l_{i_n}}''$ and $W\subset F_{l_{i_n}}''$. Since $\Fo$ is maximal, $$\dim F_{l_{i_n}}'\cap\langle d_{n+1},l_{i_n}\rangle_{\Cp}\geq1,$$ but $l_{i_n}\notin F_{l_{i_n}}'$, so $W\cap F_{l_{i_n}}'\neq\varnothing$. Furthermore, $\dim W\cap F_{l_{i_n}}'\geq n-1$ because $l_{i_1},~\ldots,~l_{i_{n-1}}\in F_{l_{i_n}}'$. On the other hand, $l_{i_n}\notin W\cap F_{l_{i_n}}'$, hence $\dim W\cap F_{l_{i_n}}'=n-1$. More precisely, $$W\cap F_{l_{i_n}}'=d_{n+1}+\langle l_{i_1},~\ldots,~l_{i_{n-1}}\rangle_{\Cp}.$$ According to our assumption, for any $l\in W\cap F_{l_{i_n}}'$, the subspace $F_l'$ coincides with one of the subspaces $F'_{l_{i_1}},~\ldots,~F'_{l_{i_{n-1}}}$ because $F_l'\neq F_{l_{i_n}}'$.

Arguing in a similar way, we see that for any $2\leq k\leq n$, $$W\cap F_{l_{i_k}}'=d_{n+1}+\langle l_{i_1},~\ldots,~l_{i_{k-1}}\rangle_{\Cp}$$ and, for any $l\in W\cap F_{l_{i_k}}'$, $F_l'$ coincides with one of the subspaces $F_{l_{i_1}}',~\ldots,~F_{l_{i_{k-1}}}'$. In particular, $$W\cap F_{l_{i_2}}'=d_{n+1}+\langle l_{i_1}\rangle_{\Cp}$$ and $F_{d_{n+1}+cl_{i_1}}'=F_{l_{i_1}}'$ for all $c\in\Cp$. This means that $\langle d_{n+1},~e_{l_1}\rangle_{\Cp}\subset F_{l_{i_1}}''$. Since the generalized flag $\Fo$ is maximal, we have $$\dim F_{l_{i_1}}'\cap\langle d_{n+1},~l_{i_1}\rangle_{\Cp}\geq1.$$ However, $l_{i_1}\notin F_{l_{i_1}}'$, thus $W\cap F_{l_{i_1}}'\neq\varnothing$. Taking any $l\in W\cap F_{l_{i_1}}'$, we see that the subspaces $F'_{l_1},~\ldots,~F'_{l_n},~F_l'$ are pairwise distinct, which is a contradiction with our assumption that $F_l'$ coincides with some of the subspaces $F'_{l_1},~\ldots,~F'_{l_n}$.

We now set $l_{n+1}=l$, where $l$ is a vector from $W$ such that the subspaces $F'_{l_1},~\ldots,~F'_{l_n},~F_l'$ are pairwise distinct. This allows us to conclude that a basis $L$ is constructed as required. It is easy to see that, given $F'\in\Fo'$, the set $F''\setminus F'$ contains exactly one vector from $L$. Moreover, $\Fo$ is compatible with $L$ because putting $\sigma(l)=\alpha$ for $l\in L$, where $l\in F_{\alpha}''\setminus F_{\alpha}'$, we obtain a surjection $\sigma\colon L\to\Au$ with the property (\ref{formula:F'_F''}).

If $\Fo$ is a non-necessarily maximal generalized flag, it is enough to consider a basis compatible with a maximal generalized flag $\Go$ containing $\Fo$. Such a basis is automatically compatible with $\Fo$.}

Next, we define a generalized flag $\Fo$ to be \emph{weakly compatible} with the basis $E$ if $\Fo$ is compatible with a basis~$L$ of~$V$ such that the set $E\setminus(E\cap L)$ is finite. Two generalized flags $\Fo$, $\Go$ are $E$-\emph{commensurable} if both of them are weakly compatible with $E$ and there exist an isomorphism of ordered sets $\phi\colon\Fo\to\Go$ and a finite-dimensional subset $U\subset V$ such that
\begin{equation*}
\begin{split}
&\text{i) $\phi(F)+U=F+U$ for all $F\in\Fo$;}\\
&\text{ii) $\dim\phi(F)\cap U=\dim F\cap U$ for all $F\in\Fo$.}
\end{split}
\end{equation*}
Given a generalized flag $\Fo$ compatible with $E$, denote by $\Fl(\Fo,E)$ the set of all generalized flags in $V$, which are $E$-commensurable with $\Fo$.

To endow $\Fl(\Fo,E)$ with an ind-variety structure, denote $$\Fo_n=\{F\cap V_n,~F\in\Fo\}.$$ Given $\alpha\in \Au$, set
\begin{equation*}
\begin{split}
d_{\alpha,n}'&=\dim F_{\alpha}'\cap V_n=|\{e\in E_n\mid\sigma(e)\prec\alpha\}|,\\
d_{\alpha,n}''&=\dim F_{\alpha}''\cap V_n=|\{e\in E_n\mid\sigma(e)\preceq\alpha\}|.
\end{split}
\end{equation*}
Let $\Fl_n$ be the projective variety of flags in $V_n$ of the form $\{U_{\alpha}',U_{\alpha}'',~\alpha\in \Au\}$, where $U_{\alpha}'$, $U_{\alpha}''$ are subspaces of $V_n$ of dimensions $d_{\alpha,n}'$, $d_{\alpha,n}''$ respectively, $U_{\alpha}'\subset U_{\alpha}''$ for all $\alpha\in \Au$, and $U_{\alpha}''\subset U_{\beta}'$ for all $\alpha\prec\beta$. (Of course, if $\Au$ is infinite, then there exists infinitely many $\alpha$, $\beta\in \Au$ such that $U_{\alpha}''=U_{\beta}'$.) We define an embedding $$\iota_n\colon \Fl_n\to \Fl_{n+1},~\{U_{\alpha}'',U_{\alpha}'',~\alpha\in \Au\}\mapsto\{W_{\alpha}'',W_{\alpha}'',~\alpha\in \Au\}$$ by setting
\begin{equation}
\begin{split}\label{formula:iota_n}
W_{\alpha}'&=U_{\alpha}'\oplus\langle e\in E_{n+1}\setminus E_n\mid\sigma(e)\prec\alpha\rangle_{\Cp},\\
W_{\alpha}''&=U_{\alpha}''\oplus\langle e\in E_{n+1}\setminus E_n\mid\sigma(e)\preceq\alpha\rangle_{\Cp}.
\end{split}
\end{equation}
Then $\iota_n$ is an embedding of smooth algebraic varieties, there exists a bijection between $\Fl(\Fo,E)$ and the inductive limit of this chain of embeddings, see \cite[Proposition 5.2]{DimitrovPenkov1} or \cite[Section 3.3]{FressPenkov1}.\break Furthermore, this ind-variety structure on $\Fl(\Fo,E)$ is independent of the exhaustion $\{E_n\}$\break of the basis $E$.

There is another way, sometimes more convenient, to describe the embedding $\iota_n$. First, for all $n\geq1$, set $\Vo_n=\langle e_1,~\ldots,~e_n\rangle_{\Cp}$. For any $\Go\in\Fl(\Fo,E)$ pick a nonnegative integer $n_{\Go}$ such that $\Go$ is compatible with a basis containing $\{e_n,~n>n_{\Go}\}$ and $\Vo_{n_{\Go}}$ contains a subspace which makes these generalized flags $E$-commensurable; clearly, we can put $n_{\Fo}=0$. For $n\geq n_{\Go}$, set also $$\Go(n)=\{W\cap\Vo_n,~W\in\Go\}.$$

Given $n\geq1$, the dimensions of the spaces of the flag $\Fo(n)$ form a sequence of integers $$0=d_{n,0}<d_{n,1}<\ldots<d_{n,d_{n,s_n-1}}<d_{n,s_n}=n=\dim\Vo_n.$$ Let $\Fl(d_n,\Vo_n)$ be the flag variety of type $d_n=(d_{n,1},\ldots,d_{n,s_n-1})$ in $\Vo_n$. Since either $s_{n+1}=s_n$ or $s_{n+1}=s_n+1$, there is a unique $j_n$ such that $d_{n+1,i}=d_{n,i}+1$ for $0\leq i<j_n$ and $d_{n+1,j_n}>d_{n,j_n}$. Then, for $j_n\leq i<s_n$, $d_{n+1,i}=d_{n,i}+1$ in case $s_{n+1}=s_n$, and $d_{n+1,i}=d_{n,i-1}+1$ in case $s_{n+1}=s_n+1$. In other words, $j_n\leq s_n$ is the minimal nonnegative integer for which there is $\alpha\in\Au$ with $$\dim F_{\alpha}''\cap\Vo_{n+1}=\dim F_{\alpha}''\cap\Vo_n+1.$$

For each $n$ we define an embedding $$\xi_n\colon\Fl(d_n,\Vo_n)\hookrightarrow\Fl(d_{n+1},\Vo_{n+1})$$ as follows: given a flag $$\Go_n=\{\{0\}=G_0^n\subset G_1^n\subset\ldots\subset G_{s_n}^n=\Vo_n\}\in\Fl(d_n,\Vo_n),$$ we set $$\xi_n(\Go_n)=\Go_{n+1}=\{\{0\}=G_0^{n+1}\subset G_1^{n+1}\subset\ldots\subset G_{s_{n+1}}^{n+1}=\Vo_{n+1}\}\in\Fl(d_{n+1},\Vo_{n+1}),$$ where
\begin{equation}
G_i^{n+1}=\begin{cases}
G_i^n,&\text{if\label{formula:iota_n_another} $0\leq i<j_n$},\\
G_i^n\oplus\Cp e_{n+1},&\text{if $j_n\leq i\leq s_{n+1}$ and $s_{n+1}=s_n$},\\
G_{i-1}^n\oplus\Cp e_{n+1},&\text{if $j_n\leq i\leq s_{n+1}$ and $s_{n+1}=s_n+1$}.
\end{cases}
\end{equation}
Note that $\xi_n(\Go(n))=\Go(n+1)$ for all $\Go\in\Fl(\Fo,E)$ and all $n\geq n_{\Go}$.

Now recall that we have the exhaustion of $E$ by its finite subsets $\{E_n\}$. Denote $m_n=|E_n|=\dim V_n$. Then, according to (\ref{formula:iota_n}), $$\iota_n=\xi_{m_{n+1}-1}\circ\xi_{m_{n+1}-2}\circ\ldots\circ\xi_{m_n}.$$ The bijection $$\Fl(\Fo,E)\to\ilm\Fl(d_{m_n},V_n)$$ mentioned above now has the form $\Go\mapsto\ilm\Go(n)$. By a slight abuse of notation, in the sequel we will denote the canonical embedding $$\Fl(d_{m_n},V_n)\hookrightarrow\Fl(\Fo,E)$$ by the same letter $\iota_n$.

We complete this subsection by giving some main examples of generalized flag ind-varieties which we will refer to throughout the paper.

\exam{i) A first example \label{exam:gen_flags}of generalized flags is provided by the flag $$\Fo=\{\{0\}\subset F\subset V\},$$ where $F$ is a proper nonzero subspace of $V$. Here $\Fo'=\{\{0\}\subset F\}$, $\Fo''=\{F\subset V\}$. If $\Fo$ is compatible with $E$, then $E\cap F$ is a basis of $F$, i.e., $F=\langle\sigma\rangle_{\Cp}$ for some subset $\sigma$ of $E$. The ind-variety $\Fl(\Fo,E)$ is called an \emph{ind-grassmannian}, and is denoted by~$\Grr{F}{E}$. If $k=\dim F$ is finite, then a flag\break $\{\{0\}\subset F'\subset V\}$ is $E$-commensurable with $\Fo$ if and only if $\dim F=k$, hence $\Grr{F}{E}$ depends only on~$k$, and we denote it by $\Gra{k}=\Grr{k}{V}$. Similarly, if $k=\codim_VF$ is finite, then $\Grr{F}{E}$ depends only on $E$ and $k$ (but not on~$F$) and is isomorphic to $\Grr{k}{V_*}$: an isomorphism $$\Grr{F}{E}\to\{F\subset V_*\mid\dim F=k\}=\Gra{k,V_*}$$ is induced by the map $$\Grr{F}{E}\ni U\mapsto U^{\#}=\{\phi\in V_*\mid\phi(x)=0\text{ for all }x\in U\}.$$ Finally, if $F$ is both infinite dimensional and infinite codimensional, then $\Grr{F}{E}$ depends on $F$ and~$E$, but all such ind-varieties are isomorphic and can be denoted by $\Gra{\infty}$, see \cite{PenkovTikhomirov1} or \cite[Section 4.5]{FressPenkov1} for the details.

ii) Our second example is the generalized flag $$\Fo=\{\{0\}=F_0\subset F_1\subset\ldots\},$$ where $F_i=\langle e_1,\ldots,e_i\rangle_{\Cp}$ for all $i\geq1$. This obviously is a flag. A flag $$\wt\Fo=\{\{0\}=\wt F_0\subset\wt F_1\subset \ldots\}$$ is $E$-commensurable with $\Fo$ if and only if $\dim F_i=\dim\wt F_i$ for all~$i$, and $F_i=\wt F_i$ for large enough $i$. The flag $\Fo$ is maximal, and $\Fo'=\Fo$, $\Fo''=\Fo\setminus\{0\}$.

iii) Next, put $\Fo=\{\{0\}=F_0\subset F_1\subset F_2\subset\ldots\subset F_{-2}\subset F_{-1}\subset V\}$, where $$F_i=\langle e_1,~e_3,~\ldots, e_{2i-1}\rangle_{\Cp},~F_{-i}=\langle\{e_j,~j\text{ odd}\}\cup\{e_{2j},~j>i\}\rangle_{\Cp}$$ for $i\geq1$. This generalized flag is maximal and is clearly not a flag. Here $\Fo'=\Fo\setminus V$, $\Fo''=\Fo\setminus\{0\}$. Note also that $\wt\Fo\in X=\Fl(\Fo,E)$ does not imply that $\wt F_i=F_i$ for $i$ large enough. For example, let $\wt F_1=\Cp e_2$, $\wt F_i=\langle e_2,~e_3,~e_5,~e_7,~\ldots,~e_{2i-1}\rangle_{\Cp}$ for $i>1$, and $$\wt F_{-i}=\langle \{e_j,~j\text{ odd},~j\geq3\}\cup\{e_2\}\cup\{e_{2j},j>i\}\rangle_{\Cp},~i\geq1,$$ then $\wt\Fo\in\Fl(\Fo,E)$, but $\wt F_i\neq F_i$ for all $i$.

iv) Now let $\Fo=\{\{0\}=F_0\subset F_1\subset F_2\subset\ldots\subset F\subset F_{-1}\subset F_{-2}\subset\ldots\}$, where $$F_i=\langle e_1,~e_3,~\ldots, e_{2i-1}\rangle_{\Cp},~F=\langle e_j,~j\text{ odd}\rangle_{\Cp},~F_{-i}=\langle\{e_J,~j\text{ odd}\}\cup\{e_{2j},~1\leq j\leq i\}\rangle_{\Cp}.$$
The chain $\Fo$ is a maximal generalized flag but not a flag. Note that the space $F$ does not have an immediate predecessor. We have $\Fo'=\Fo$, $\Fo''=\Fo\setminus(\{0\}\cup F)$.

v) Finally, let $\Qp$ be the field of rational numbers. Fix a bijection $\sigma\colon E\to\Qp$ and for each $\alpha\in\Qp$ set
\begin{equation*}
F_{\alpha}'=\langle e\in E,~\sigma(e)<\alpha\rangle_{\Cp},~F_{\alpha}''=\langle e\in E,~\sigma(e)\leq\alpha\rangle_{\Cp}.
\end{equation*}
The subspaces $\{F_{\alpha}',~F_{\alpha}''\}_{\alpha\in\Qp}$ form a maximal generalized flag $\Fo$ with $\Au=\Qp$. Of course, $\Fo$ is not a flag. None of the sub\-spaces~$F_{\alpha}'$ has an immediate predecessor, and none of the subspaces $F_{\alpha}''$ has an immediate successor, so $\Fo'=\{F_{\alpha}'\}_{\alpha\in\Qp}$ and $\Fo''=\{F_{\alpha}''\}_{\alpha\in\Qp}$.

Furthermore, note that $\Fo$ is not maximal among all chains of subspaces of $V$. Indeed, given $\gamma\in\Rp\setminus\Qp$, set
$C_{\gamma}=\langle e\in E,~\sigma(e)<\gamma\rangle_{\Cp}$. The subspaces $\{F_{\alpha}',~F_{\alpha}''\}_{\alpha\in\Qp}\sqcup\{C_{\gamma}\}_{\gamma\in\Rp\setminus\Qp}$ form a maximal chain $\Co$ of subspaces in $V$, and it is clear that $\Fo$ is a proper subchain of $\Co$. The chain $\Co$ is parameterized by the ``symmetric Dedekind cuts of $\Qp$'': a point of $\Rp\setminus\Qp$ corresponds to one cut, and a point of $\Qp$ corresponds to two cuts featuring respectively a maximum and a minimum.}

\subsection{Isotropic generalized flags}\fakesst In\label{sst:iso_gen_flags} this subsection we assume that $V$ is endowed with a nondegenerate symmetric or skew-symmetric bilinear form $\beta$ such that $\beta(e,V_n)=0$ for $e\in E\setminus E_n$ as in Example~\ref{exam:ind_groups} ii). Given $U\subset V$, we set $U^{\perp}=\{x\in V\mid\beta(x,y)=0\text{  for all }y\in U\}$, i.e., $U^{\perp}$ is the maximal subspace of $V$ orthogonal to $U$ with respect to $\beta$. We suppose that $E$ is $\beta$-\emph{isotropic} (or, equivalently, \emph{isotropic}), i.e., that there exists an involution $i_E\colon E\to E$ with at most one fixed point and satisfying $\beta(e,e')=0$ for all $e,~e'\in E$ such that $e'\neq i_E(e)$ (here $e$ may equal $e'$).

\defi{A generalized flag $\Fo$ is called $\beta$-\emph{isotropic} (or simply \emph{isotropic}) if $F^\perp\in\Fo$ for all $F\in\Fo$, and if the map $i_{\Fo}\colon\Fo\to\Fo,~F\mapsto F^\perp$ is an involutive anti-automorphism of the ordered set~$\Fo$, i.e., a bijection which reverses the inclusion relation. Note that the involution on $\Fo$ induces the involution $$i_{\Au}\colon\Au\to\Au,~(F_{\alpha}',F_{\alpha}'')\mapsto((F_{\alpha}'')^{\perp},(F_{\alpha}')^{\perp}).$$}

Let a generalized flag $\Fo$ be $\beta$-isotropic. It follows from the definition that each $\beta$-isotropic flag has the form $\Go\cup\Go^{\perp}$, where $\Go$ consists of isotropic subspaces of $V$, and $\Go^{\perp}$ consists of the respective orthogonal spaces. Denote by $U$ the union of all subspaces of $\Fo$ which belong to $\Go$. Then $U$ is an isotropic subspace of $V$, $\Go$ is a generalized flag in $U$ (possibly, containing $U$ as an element), and $\Fo$ is uniquely determined by its intersection $\Go=\Fo\cap U$ with the subspace $U$. Furthermore, $U$ is a maximal isotropic subspace of $V$ if and only if $\Fo$ is a maximal generalized flag in $V$.

Arguing as in the proof of Proposition~\ref{prop:comp_basis}, one can prove that each $\beta$-isotropic generalized flag has a compatible $\beta$-isotropic basis. In the rest of this subsection we assume that $\Fo$ is $\beta$-isotropic and compatible with a fixed $\beta$-isotropic basis~$E$. Now, suppose that a generalized flag~$\wt\Fo$ is $\beta$-isotropic and $E$-commensurable with~$\Fo$. In particular, the set $\wt\Fo^{\dag}$ of consecutive subspaces from $\wt\Fo$ is isomorphic to $\Au$, and the involution $i_{\wt\Fo}$ induces the same involution $i_{\Au}$ on the linearly ordered set $\Au$. The following lemma is proved in \cite[Subsection~3.1, Lemma1]{FressPenkov1}.

\lemmp{\textup{i)} There exists\label{lemm:iso_flags} a $\beta$-isotropic basis $L$ such that $E\setminus(E\cap L)$ is finite and $\wt\Fo$ is compatible with $L$. \textup{ii)} If $\wt\Fo$ is compatible with a $\beta$-isotropic basis $L$ via the surjective map $\sigma\colon L\to\Au$, then $$\sigma\circ i_L=i_{\Au}\circ\sigma.$$}{i) Denote by $L'$ a basis of $V$ such that $E\setminus(E\cap L')$ is finite and $\wt\Fo$ is compatible with $L'$.
Pick a subset $E'\subset E$ stable under the involution $i_E$ such that $i_E$ has no fixed point in $E'$, $E\setminus E'$ is finite, and $E'\subset E\cap L'$. Then $V''=\langle E\setminus E'\rangle_{\Cp}$ is a finite-dimensional space and the restriction of $\beta$ to $V''$ is nondegenerate. The intersections $\{\wt F\cap V'',~\wt F\in\wt\Fo\}$ form an isotropic flag of $V''$. Since $V''$ is finite dimensional, there exists a $\beta$-isotropic basis $E''$ of $V''$ such that this isotropic flag is compatible with $E''$. Then $L=E'\cup E''$ is a required basis.

ii) By definition of compatibility, $e\in\wt F_{\sigma(e)}''\setminus\wt F_{\sigma(e)}'$ for each $e\in L$, hence $$i_L(e)\in (\wt F_{\sigma(e)}')^{\perp}\setminus(\wt F_{\sigma(e)}'')^{\perp}.$$ The result follows.}

Now, denote by $\Fl(\Fo,\beta,E)$ the set of all $\beta$-isotropic generalized flags in $V$ which are $E$-com\-men\-surable to $\Fo$. To endow $\Fl(\Fo,\beta,E)$ with a structure of ind-variety, assume that in the exhaustion $\{E_n\}$ of $E$ all the subsets $E_n$ are $i_E$-stable. Recall the definition of $\Fl_n$ from the previous subsection. Let $U^{\perp,V_n}\subset V_n$ be the orthogonal subspace in $V_n$ of a subspace $U\subset V_n$ with respect to $\beta_n$. Denote by $\Fl_n^{\beta}$ the closed subvariety of $\Fl_n$ consisting of all flags $\{U_{\alpha}',U_{\alpha}'',~\alpha\in\Au\}$ from the flag variety $\Fl_n$ such that $$((U_{\alpha}'')^{\perp,V_n},(U_{\alpha}')^{\perp,V_n})=(U_{i_{\Au}(\alpha)}',U_{i_{\Au}(\alpha)}'')\text{ for all }\alpha\in\Au.$$

Then the embedding $\iota_n\colon\Fl_n\hookrightarrow\Fl_{n+1}$ defined by formula (\ref{formula:iota_n}) restricts to an embedding $$\iota_n^{\beta}\colon\Fl_n^{\beta}\hookrightarrow\Fl_{n+1}^{\beta}.$$
Hence we obtain a chain of embeddings of projective varieties $$\Fl_1^{\beta}\stackrel{\iota_1^{\beta}}{\hookrightarrow}\Fl_2^{\beta}\stackrel{\iota_2^{\beta}}
{\hookrightarrow}\ldots\stackrel{\iota_{n-1}^{\beta}}{\hookrightarrow}\Fl_n^{\beta}
\stackrel{\iota_n^{\beta}}{\hookrightarrow}\Fl_{n+1}^{\beta}\stackrel{\iota_{n+1}^{\beta}}{\hookrightarrow}\ldots.$$
There exists a bijection between $\Fl(\Fo,\beta,E)$ and the inductive limit $\ilm\Fl_n^{\beta}$ of this chain of embeddings. Thus, $\Fl(\Fo,\beta,E)$ is endowed with an ind-variety structure independent of the exhaustion of $E$. Further, $\Fl(\Fo,\beta,E)$ is a closed ind-subvariety of $\Fl(\Fo,E)$, see \cite{DimitrovPenkov1} or \cite{FressPenkov1}.

\exam{i) Let $F$ be an isotropic subspace\label{exam:gen_flags_iso} of $V$, and $U$ be a maximal isotropic subspace of~$V$ containing $F$. Note that $U$ is always infinite dimensional. Set $\Fo=\{\{0\}\subset F\subset F^{\perp}\subset V\}$. Then $\Fo$ is a $\beta$-isotropic flag in $V$. Let $E$ be a $\beta$-isotropic basis of $V$ compatible with the flag $\Fo$.\break If $\dim F=k<\infty$, then we denote the corresponding ind-variety $\Fl(\Fo,\beta,E)$ by $\Grab{k}=\Grb{k}{V}$. In general, we denote the corresponding ind-variety by $\Grb{F}{E}$ and call it an \emph{isotropic ind-grass\-mannian}. Note that if, for example, $F=U$ then there exist maximal isotropic subspaces of $V$ which are not contained in $\Grb{U}{E}$.

ii) Let $\Fo=\{\ldots\subset F_{-2}\subset F_{-1}\subset F_1\subset F_2\subset\ldots\subset F\subset F_1'\subset F_2'\subset\ldots\}$ where, for $i>0$,
\begin{equation*}
\begin{split}
F_{-i}&=\langle e_{3j},~j\geq i\rangle_{\Cp},\\
F_i&=\langle\{e_{3j},~ j\in\Zp_{>0}\}\cup\{e_{3j-1},~1\leq j\leq i\}\rangle_{\Cp},\\
F&=\langle e_{3j},~e_{3j-1},~j\in\Zp_{>0}\rangle_{\Cp},\\
F_i'&=\langle\{e_{3j},~e_{3j-1},~j\in\Zp_{>0}\}\cup\{e_{3j-2},~1\leq j\leq i\}\rangle_{\Cp}.\\
\end{split}
\end{equation*}
Here $\Fo'=\Fo$, $\Fo''=\Fo\setminus\{F\}$, and there is no nondegenerate symmetric or skew-symmetric bilinear form $\beta$ on $V$ with respect to which $\Fo$ is isotropic. Indeed, it is clear that $\Fo$ does not admit an involutive anti-automorphism because there exists a unique subspace of $\Fo$ without an immediate predecessor, while all subspaces of $\Fo$ have immediate successors.}

\section{Linear ind-grassmannians}\label{sect:ind_grassmannians}\fakesect

In this section, we consider a more general approach to the definition of an ind-grassmannian based on the notion of a linear embedding of finite-dimensional grassmannians. Theorem~\ref{theo:class_grass} below claims that any ind-variety obtained as an inductive limit of linear embeddings of grassmannians is isomorphic to an ind-grassmannian as defined in Exam\-ples~\ref{exam:gen_flags}~i) and \ref{exam:gen_flags_iso} i). This material is taken from \cite{PenkovTikhomirov1}.

\subsection{Definition of linear ind-grassmannians}\fakesst \label{sst:linear_ind_grassmannians}Given an algebraic variety $X$, we denote by $\Pic{X}$ its Picard group, i.e. the group of isomorphism classes of line bundles. The group operation here is tensor product. If $X$ is a projective variety with Picard group isomorphic to $\Zp$, then $\Ou_X(1)$ stands for the ample generator of the Picard group, and $\Ou_X(n)=\Ou_X(1)^{\otimes n}$ for $n\in\Zp$. Note that each morphism $\vfi\colon X\to Y$ of algebraic varieties induces a group homomorphism $\vfi^*\colon\Pic{Y}\to\Pic{X}$.

If $X$ is an ind-variety obtained as the inductive limit of a chain of morphisms
\begin{equation*}X_1\stackrel{\vfi_1}\to X_2\stackrel{\vfi_2}\to\ldots\stackrel{\vfi_{n-1}}{\to}
X_n\stackrel{\vfi_n}\to X_{n+1}\stackrel{\vfi_{n+1}}\to\ldots,
\end{equation*} of algebraic varieties, then, by definition, the \emph{Picard group} of $X$ is the projective limit $$\Pic{X}=\plm\Pic{X_n}$$ of the chain of group homomorphisms
\begin{equation*}\Pic{X_1}\stackrel{\vfi_1^*}\gets\Pic{X_2}\stackrel{\vfi_2^*}\gets\ldots\stackrel{\vfi_{n-1}^*}\gets\Pic{X_n}
\stackrel{\vfi_n^*}\gets\Pic{X_{n+1}}\stackrel{\vfi_{n+1}^*}\gets\ldots.
\end{equation*}
Clearly, each morphism $\vfi\colon X\to Y$ of ind-varieties induces a group homomorphism $\vfi^*\colon\Pic{Y}\to\Pic{X}$ of their Picard groups.

\defi{We call a morphism $\vfi\colon X\to Y$ of\label{defi:linear_morphism} algebraic varieties or ind-varieties \emph{linear} if $\vfi^*$ is an epimorphism of Picard groups.}

For example, let $X=\Grr{k}{W}$ be the grassmannian of $k$-dimensional subspaces of an $n$-dimensional vector space $W$. Then $\Pic{X}\cong\Zp$ and $\Ou_X(1)\cong\bigwedge^kS_X^*$, where $S_X$ is the tautological bundle on $X$. If $Y$ is also a grassmannian, then a morphism $\vfi\colon X\to Y$ is linear if and only if $\vfi^*\Ou_Y(1)=\Ou_X(1)$.

In the sequel we will also consider orthogonal and symplectic grassmannians. Assume that a finite-dimensional vector space $W$ is endowed with a nondegenerate symmetric or skew-symmetric bilinear form $\beta$. Given $k\leq[\dim W/2]$, the \emph{isotropic grassmannian} $\Grb{k}{W}$ is the subvariety of $\Grr{k}{W}$ consisting of all $k$-dimensional isotropic subspaces of~$W$. If $\beta$ is symmetric (respectively, skew-symmetric) then $\Grb{k}{W}$ is called \emph{orthogonal} (respectively, \emph{symplectic}).

We will assume throughout this section that if $\beta$ is symmetric, then $\dim W\geq7$ and $k\neq(\dim W)/2$, $k\neq(\dim W)/2-1$. It is well known that $\Grb{k}{W}$ is smooth and
\begin{equation*}
\dim\Grb{k}{W}=\begin{cases}
k\dim W-k(3k+1)/2,&\text{if $\beta$ is symmetric},\\
k\dim W-k(3k-1)/2,&\text{if $\beta$ is skew-symmetric}.\\
\end{cases}
\end{equation*}
Furthermore, $\Pic{\Grb{k}{W}}=\Zp\Ou_{\Grb{k}{W}}(1)$, where $\Ou_{\Grb{k}{W}}(1)$ satisfies the following condition: if $\tau\colon\Grb{k}{W}\hookrightarrow\Grr{k}{W}$ is the tautological embedding, then
\begin{equation*}
\tau^*\Ou_{\Grr{k}{W}}(1)=\begin{cases}
\Ou_{\Grb{k}{W}}(2),&\text{if $\beta$ is symmetric and $k=[\dim W/2]$},\\
\Ou_{\Grb{k}{W}}(1)&\text{otherwise}.\\
\end{cases}
\end{equation*}
Note that a morphism $\vfi\colon X\to Y$ of orthogonal or symplectic grassmannians is linear if and only if $\vfi^*\Ou_Y(1)=\Ou_X(1)$.

\defi{A \emph{linear ind-grassmannian} is an\label{defi:linear_ind_grassmannian} ind-variety $X$ obtained as the inductive limit $\ilm X_n$ of a chain of embeddings
\begin{equation*}X_1\stackrel{\vfi_1}\hookrightarrow X_2\stackrel{\vfi_2}\hookrightarrow\ldots\stackrel{\vfi_{n-1}}{\hookrightarrow}
X_n\stackrel{\vfi_n}\hookrightarrow X_{n+1}\stackrel{\vfi_{n+1}}\hookrightarrow\ldots,
\end{equation*}
where each $X_n$ is a grassmannian or isotropic grassmannian with $\Pic{X}\cong\Zp$, $\lim_{n\to\infty}\dim X_n=\infty$, and all em\-bed\-dings $\vfi_n$ are linear morphisms. Note that we allow a mixture of all three types of grassmannians (usual, orthogonal and symplectic).}

\exam{Let $V$ and\label{exam:grass_like_ind_grass} $E$ be as above, and $F$ be a subspace of the space $V$ such that the flag $\Fo=\{\{0\}\subset F\subset V\}$ is compatible with $E$. Set $\Grr{F}{E}=\Fl(\Fo,E)$ as in Example~\ref{exam:gen_flags} i). Then $\Grr{F}{E}\cap V_n=\Grr{d_n}{V_n}$ for $d_n=\dim F\cap V_n$, and the embedding $\iota_n\colon\Grr{d_n}{V_n}\hookrightarrow\Grr{d_{n+1}}{V_{n+1}}$ has the following simple form: $$\iota_n(A)=A\oplus U_{n+1}\text{ for }A\in\Grr{d_n}{V_n},$$ where $U_{n+1}$ is a fixed subspace of $V_{n+1}$ spanned by certain basis vectors from $E_{n+1}\setminus E_n$, see (\ref{formula:iota_n_another}). Such embeddings $\iota_n$ are clearly linear, so $\Grr{F}{E}$ is a linear ind-grassmannian.}


\subsection{Standard extensions}\fakesst The key\label{sst:standard extensions} idea in the description of linear ind-grassmannians is that each linear ind-grassmannian is isomorphic to the inductive limit of a chain of certain standard embeddings.

\defi{Let $X$, $X'$ be usual grassmannians. An embedding $\wt\vfi\colon X\to X'$ is called a \emph{standard extension} if there exist isomorphisms $$j_X\colon X\to\Grr{k}{W},~j_{X'}\colon X'\to\Grr{k'}{W'}$$ and an embedding $\vfi\colon\Grr{k}{W}\to\Grr{k'}{W'}$ for $\dim W'\geq\dim W$, $k'\geq k$, such that the diagram
\begin{equation*}
\xymatrix{
X\ar[d]^{j_X}\ar@{^{(}->}[r]^{\wt\vfi}&X'\ar[d]^{j_{X'}}\\
\Grr{k}{W}\ar@{^{(}->}[r]^{\vfi}&\Grr{k'}{W'}
}
\end{equation*}
is commutative, and $\vfi$ is given by the formula
\begin{equation}\label{formula:standard_extension}
\vfi(A)=A\oplus U,~A\in\Grr{k}{W}
\end{equation}
for some fixed isomorphism $W'\cong W\oplus U'$ and a fixed subspace $U\subset U'$ of dimension $k'-k$. Note that $$\dim W'-\dim W=\dim U'\geq\dim U=k'-k,$$ hence $\dim W'-k'\geq\dim W-k$.}

For instance, all embeddings considered in Example~\ref{exam:grass_like_ind_grass} are standard extensions. It is clear that the composition of two standard extensions is a standard extension.

In what follows we will say that a standard extension $\vfi\colon\Grr{k}{W}\to\Grr{k'}{W'}$ is \emph{strict} if $j_{\Grr{k}{W}}$ and $j_{\Grr{k'}{W'}}$ are automorphisms. Given a strict standard extension $\vfi\colon\Grr{k}{W}\to\Grr{k'}{W'}$, the isomorphism $W'\cong W\oplus U'$ can always be changed so that $\vfi$ is given simply by formula (\ref{formula:standard_extension}). Obviously, the composition of two strict standard extensions is a strict standard extension. Note that, given a strict standard extension $\vfi\colon\Grr{k}{W}\to\Grr{k'}{W'}$, $U$~can be recovered by the formula $$U=\bigcap\nolimits_{A\in\Grr{k}{V}}\vfi(A).$$ Put also $$U^{\sharp}=\left\langle\vfi(A),~A\in\Grr{k}{W}\right\rangle_{\Cp},$$ then $\vfi$ determines a surjective linear map $\vfi^{\sharp}\colon U^{\sharp}\to W$ with kernel $U$ such that $\left(\vfi^{\sharp}\right)^{-1}(A)=\vfi(A)$ for all $A\in\Grr{k}{V}$. One can easily check that fixing the standard extension $\vfi$ is equivalent to fixing the triple $(U,U^{\sharp},\vfi^{\sharp})$.

Suppose now that $W$ and $W'$ are finite-dimensional spaces endowed with respective nondegenerate bilinear forms $\beta$ and $\beta'$ such that $\beta$ and $\beta'$ are both symmetric or skew-symmetric. An embedding $\vfi\colon\Grb{k}{W}\hookrightarrow\Grb{k'}{W'}$ is called a \emph{standard extension} if $\vfi$ is given by formula (\ref{formula:standard_extension}) where\break $W'\cong W\oplus U'$ is an isometry and $U$ is an isotropic subspace of $U'$.

As for usual grassmannians, a standard extension of isotropic grassmannians can be defined by the following linear-algebraic datum. Pick a flag $\{U\subset U^{\sharp}\}$ in $W'$ for which $U$ is isotropic and there exists a surjective linear map $\vfi^{\sharp}\colon U^{\sharp}\to W$ with kernel $U$ so that the bilinear form $\left(\vfi^{\sharp}\right)^*\beta$ coincides with the restriction of $\beta'$ to $U^{\sharp}$. This datum defines an embedding $\vfi\colon\Grb{k}{W}\hookrightarrow\Grb{k'}{W'}$ by the formula
\begin{equation*}
\vfi(A)=\left(\vfi^{\sharp}\right)^{-1}(A)\subset U^{\sharp}\subset W',~A\in\Grb{k}{W}.
\end{equation*}
Moreover,
\begin{equation*}
\begin{split}
U&=\bigcap\nolimits_{A\in\Grb{k}{W}}\vfi(A),\\
U^{\sharp}&=\langle\vfi(A),~A\in\Grb{k}{W}\rangle_{\Cp}.
\end{split}
\end{equation*}

To formulate the first main result of this section, we need some more definitions. Let $\wt W$ be an isotropic subspace of $W$. For $k\leq\dim\wt W$, we call the natural embeddings $$\Grr{k}{\wt W}\hookrightarrow\Grb{k}{W},~\Grr{\dim\wt W-k}{\wt W^*}\hookrightarrow\Grb{k}{W}$$ \emph{isotropic extensions}. By definition, a \emph{combination of isotropic and standard extensions} is an embedding of the form
\begin{equation*}
\Grb{k}{W}\stackrel{\tau}{\hookrightarrow}\Grr{k}{W}\stackrel{\vfi}{\hookrightarrow}\Grr{k''}{\wt W''}\stackrel{b}{\hookrightarrow}\Grb{k''}{W''}\stackrel{\psi}{\hookrightarrow}\Grb{k'}{W'},
\end{equation*}
where $\wt W''$ is an isotropic subspace of $W''$, $\tau$ is the tautological embedding, $\vfi$ and $\psi$ are standard extensions and $b$ is an isotropic extension. It is easy to check that a composition of combinations of isotropic and standard extensions is a combination of isotropic and standard extensions.

By \emph{projective space in} (or \emph{on}) a variety (or an ind-variety) $X$ we understand a linearly embedded subvariety of $X$ isomorphic to a projective space. Similarly, by a \emph{quadric on} $X$ of dimension $m\geq3$ we understand a linearly embedded subvariety of $X$ isomorphic to a smooth $m$-dimensional quadric. (For quadrics on $X$ of dimensions 1 and 2 the definition is slightly different, see \cite[Sub\-section~2.2]{PenkovTikhomirov1} for the details.)

Now, suppose that $X$ and $Y$ are usual (respectively, orthogonal or symplectic) grassmannians. If these are orthogonal grassmannians of the form $\Grb{k}{W}$ and $\Grb{k'}{W'}$ respectively, then assume in addition that either $$k\leq[(\dim W)/2]-3,~k'\leq[(\dim W')/2]-3,$$ or that both $\dim W$, $\dim W'$ are odd and $$[(\dim W')/2]-k'\leq[(\dim W)/2]-k\leq2.$$

It follows immediately from the definitions that standard extensions and combinations of isotropic and standard extensions are linear morphisms.

The following theorem is the central result which leads to the classification of linear\break ind-grassmannians, see Subsection~\ref{sst:classification_ind_grassmannians}. We invite the reader to read the proof of Theorem~\ref{theo:stand_exts} in the original paper \cite[Theorem 1]{PenkovTikhomirov1}.

\mtheo{Let $\vfi\colon X\to Y$ be\label{theo:stand_exts} a linear morphism. Then some of the following statements hold\textup:
\begin{equation*}
\begin{split}
&\text{\textup{a)} $\vfi$ is a standard extension\textup;}\\
&\text{\textup{b)} $X$ and $Y$ are isotropic grassmannians}\\
&\text{$\hphantom{\text{\textup{b)} }}$and $\vfi$ is a combination of isotropic and standard extensions\textup;}\\
&\text{\textup{c)} $\vfi$ factors through a projective space in $Y$}\\
&\text{$\hphantom{\text{\textup{c)} }}$or\textup, in case of orthogonal grassmannians\textup, through a maximal quadric on $Y$\textup.}
\end{split}
\end{equation*}}

In particular, $\vfi$ is an embedding unless it factors through a projective space in $Y$ or, in case of orthogonal grassmannians, through a maximal quadric on $Y$.

\subsection{Classification of linear ind-grassmannians}\fakesst\label{sst:classification_ind_grassmannians} In this subsection we explain that each linear (possibly, isotropic) ind-grassmannian is isomorphic to one of the standard (possibly, isotropic) ind-grassmannians defined below.

Let $V$ be a countable-dimensional vector space, $E$ be a basis of $V$, $E=\bigcup E_n$ be its exhaustion of finite subsets and $V=\bigcup V_n$ be the corresponding exhaustion of $V$ by its finite-dimensional subspaces $V_n=\langle E_n \rangle_{\Cp}$. Denote by $\Gra{k}$ the inductive limit of a sequence
\begin{equation*}
\Grr{k}{V_1}\stackrel{\vfi_1}\hookrightarrow\Grr{k}{V_2}\stackrel{\vfi_2}\hookrightarrow\ldots\stackrel{\vfi_{n-1}}{\hookrightarrow}
\Grr{k}{V_n}\stackrel{\vfi_n}\hookrightarrow\Grr{k}{V_{n+1}}\stackrel{\vfi_{n+1}}\hookrightarrow\ldots,
\end{equation*}
where $k\geq1$ is an integer, and all $\vfi_n$ are canonical inclusions of grassmannians.

Denote also by $\Gra{\infty}$ the inductive limit of a sequence
\begin{equation*}
\Grr{k_1}{V_1}\stackrel{\vfi_1}\hookrightarrow\Grr{k_2}{V_2}\stackrel{\vfi_2}\hookrightarrow\ldots\stackrel{\vfi_{n-1}}{\hookrightarrow}
\Grr{k_n}{V_n}\stackrel{\vfi_n}\hookrightarrow\Grr{k_{n+1}}{V_{n+1}}\stackrel{\vfi_{n+1}}\hookrightarrow\ldots,
\end{equation*}
where $1\leq k_1<k_2<\ldots$ are integers satisfying $$\lim_{n\to\infty}(\dim V_n-k_n)=\infty,$$ and all $\vfi_n$ are standard extensions.

In the orthogonal and symplectic cases we assume that $V$ is endowed with a respective non\-de\-ge\-ne\-rate symmetric or skew-symmetric bilinear form $\beta$ such that the restriction of $\beta$ to $V_n$ is nondegenerate for all $n$. Here we don't assume that $e$ is orthogonal to $V_n$ for $e\in E\setminus E_n$. For an integer $$1\leq k\leq[(\dim V_1)/2],$$ let $\Grb{k}{\infty}$ denote the inductive limit of the chain
\begin{equation*}
\Grbn{k}{V_1}{\beta_1}\stackrel{\vfi_1}\hookrightarrow\Grbn{k}{V_2}{\beta_2}\stackrel{\vfi_2}\hookrightarrow
\ldots\stackrel{\vfi_{n-1}}{\hookrightarrow}
\Grbn{k}{V_n}{\beta_n}\stackrel{\vfi_n}\hookrightarrow\Grbn{k}{V_{n+1}}{\beta_{n+1}}\stackrel{\vfi_{n+1}}\hookrightarrow\ldots,
\end{equation*}
where $\beta_n$ denotes the restriction of the form $\beta$ to $V_n$, and all morphisms $\vfi_n$ are canonical inclusions of isotropic grassmannians.

Given a sequence of integers $$1\leq k_1<k_2<\ldots$$ such that $k_n<[(\dim V)_n/2]$ for all $n$ (and, consequently, $\lim_{n\to\infty}(\dim V_n-k_n)=\infty$), we denote by $\Grb{\infty}{\infty}$ the inductive limit of a chain
\begin{equation}\label{formula:iso_stand_grass}
\Grbn{k_1}{V_1}{\beta_1}\stackrel{\vfi_1}\hookrightarrow\Grbn{k_2}{V_2}{\beta_2}\stackrel{\vfi_2}\hookrightarrow
\ldots\stackrel{\vfi_{n-1}}{\hookrightarrow}
\Grbn{k_n}{V_n}{\beta_n}\stackrel{\vfi_n}\hookrightarrow\Grbn{k_{n+1}}{V_{n+1}}{\beta_{n+1}}\stackrel{\vfi_{n+1}}\hookrightarrow\ldots
\end{equation}
of standard extensions of isotropic grassmannians.

Next, in the symplectic case, for a sequence of integers $$1\leq k_1<k_2<\ldots$$ satisfying $k_n\leq(\dim V_n)/2$ and $$\lim_{n\to\infty}((\dim V_n)/2-k_n)=k\geq0,$$ we denote by $\Grb{\infty}{k}$ the inductive limit of a chain (\ref{formula:iso_stand_grass}) of standard extensions. In the orthogonal case, suppose first that $\dim V_n$ is even for all $n$; then let $\Grbi{\infty}{k}{0}$ be the inductive limit of a chain (\ref{formula:iso_stand_grass}) where $k_n<(\dim V_n)/2$ and $$\lim_{n\to\infty}((\dim V_n)/2-k_n)=k\geq2.$$ Finally, if $\dim V_n$ is odd for all $n$ in the orthogonal case, we denote by $\Grbi{\infty}{k}{1}$ the inductive limit of a chain (\ref{formula:iso_stand_grass}) with $k_n\leq[\dim V_n/2]$ and $$\lim_{n\to\infty}([(\dim V_n)/2]-k_n)=k\geq0.$$

\defi{The above ind-varieties are called \emph{standard ind-grassmannians}.}
\vspace{-0.2cm}
\lemmp{All standard ind-grassmannians are well-defined\textup, i.e. a standard ind-grassmannian does not depend \textup(up to an isomorphism of ind-varieties\textup) on the specific chain used in its definition.}{Consider $\Gra{\infty}$ (all other cases are similar). Suppose that we have two chains of strict standard extensions
\begin{equation*}
\begin{split}
&\Grr{k_1}{V_1}\stackrel{\vfi_1}\hookrightarrow\Grr{k_2}{V_2}\stackrel{\vfi_2}\hookrightarrow\ldots\stackrel{\vfi_{n-1}}{\hookrightarrow}
\Grr{k_n}{V_n}\stackrel{\vfi_n}\hookrightarrow\Grr{k_{n+1}}{V_{n+1}}\stackrel{\vfi_{n+1}}\hookrightarrow\ldots,\\
&\Grr{k_1'}{V_1'}\stackrel{\vfi_1'}\hookrightarrow\Grr{k_2}{V_2'}\stackrel{\vfi_2'}\hookrightarrow\ldots\stackrel{\vfi_{n-1}}{\hookrightarrow}
\Grr{k_n}{V_n'}\stackrel{\vfi_n'}\hookrightarrow\Grr{k_{n+1}}{V_{n+1}'}\stackrel{\vfi_{n+1}'}\hookrightarrow\ldots,\\
\end{split}
\end{equation*}
where $E=\bigcup E_n'$ is an exhaustion of the basis $E$ by its finite subsets, $V_n'=\langle E_n'\rangle_{\Cp}$, and $$\lim_{n\to\infty}k_n=\lim_{n\to\infty}k_n'=\lim_{n\to\infty}(\dim V_n-k_n)=\lim_{n\to\infty}(\dim V_n'-k_n')=\infty.$$ We must show that the inductive limits $\Gra{\infty}$ and $\Grad{\infty}$ of these two chains are isomorphic.

To do this, we find $n$ such that $$\dim V_n'\geq\dim V_1,~k_n'\geq k_n,~\dim V_n'-k_n\geq\dim V_1-k_1$$ and consider an arbitrary strict standard extension $$f\colon\Grr{k_1}{V_1}\hookrightarrow\Grr{k_n'}{V_n'}.$$ Let $m$ be such that $$k_m\geq k_n',~\dim V_m\geq\dim V_n',~\dim V_m-k_m\geq\dim V_n'-k_n.$$ Denote $$\vfi=\vfi_{m-1}\circ\ldots\vfi_1\colon\Grr{k_1}{V_1}\hookrightarrow\Grr{k_m}{V_m}.$$ It is enough to construct a strict standard extension $$g\colon\Grr{k_n'}{V_n'}\hookrightarrow\Grr{k_m}{V_n}$$ such that $g\circ f=\vfi$. (Then, repeating this procedure, we will construct two mutually inverse morphisms of ind-varieties $\Gra{\infty}$ and $\Grad{\infty}$.)

As we mentioned above, the strict standard extensions $f$ and $\vfi$ are given by triples $(U_f,U_f^{\sharp},f^{\sharp})$ and $(U_{\vfi},U_{\vfi}^{\sharp},\vfi^{\sharp})$ respectively, where $\{U_f\subset U_f^{\sharp}\}$ and $\{U_{\vfi}\subset U_{\vfi}^{\sharp}\}$ are flags in $V_n'$ and $V_m$ respectively, $f^{\sharp}\colon U_f^{\sharp}\to V_1$ and $\vfi^{\sharp}\colon U_{\vfi}^{\sharp}\to V_1$ are linear surjections, and the triples
\begin{equation*}
\begin{split}
&0\to U_f\hookrightarrow U_f^{\sharp}\stackrel{f^{\sharp}}{\twoheadrightarrow}V_1\to0,\\
&0\to U_{\vfi}\hookrightarrow U_{\vfi}^{\sharp}\stackrel{\vfi^{\sharp}}{\twoheadrightarrow}V_1\to0
\end{split}
\end{equation*}
are exact.

Since $k_m>k_n'$, we obtain
\begin{equation*}
\begin{split}
\dim U_{\vfi}^{\sharp}&=\dim V_1+\dim U_{\vfi}=\dim V_1+(k_m-k_1)>\dim V_1+(k_n'-k_1)\\
&=\dim V_1+\dim U_f=\dim U_f^{\sharp}.
\end{split}
\end{equation*}
Since $f^{\sharp}$ and $\vfi^{\sharp}$ are surjective, there exists a linear surjection $$\epsi\colon U_{\vfi}^{\sharp}\twoheadrightarrow U_f^{\sharp}$$ satisfying $\vfi^{\sharp}=f^{\sharp}\circ\epsi$. Then the restriction of $\epsi$ to $U_{\vfi}$ is a well-defined linear surjection $U_{\vfi}\twoheadrightarrow U_f$. Set $U_g=\Ker\epsi$, then the triple
$$0\to U_g\hookrightarrow U_{\vfi}^{\sharp}\stackrel{\epsi}{\twoheadrightarrow}U_f^{\sharp}\to0$$ is exact. Now, set $\wt U_g^{\sharp}=U_g\oplus V_n'$ and set $$\pi\colon\wt U_g^{\sharp}\to V_n'$$ to be the projection on $V_n'$ along $U_g$.

Fix embeddings $j\colon U_{\vfi}^{\sharp}\hookrightarrow\wt U_g^{\sharp}$ and $i\colon\wt U_g^{\sharp}\hookrightarrow V_m$ such that $$\restr{(i\circ j)}{U_{\vfi}^{\sharp}}=\id_{U_{\vfi}^{\sharp}},~\pi\circ j=\restr{\epsi}{U_{\vfi}^{\sharp}}.$$ Such embeddings exist. Indeed,
\begin{equation*}
\begin{split}
\dim\wt U_g^{\sharp}&=\dim U_g+\dim V_n'=(\dim U_{\vfi}^{\sharp}-\dim U_f^{\sharp})+\dim V_n'\\
&=(k_m-k_n')+\dim V_n'=k_m+(\dim V_n'-k_n')\\
&\leq k_m+\dim V_m-k_m=\dim V_m.
\end{split}
\end{equation*}
Now, let $Z$ be a subspace of $U_{\vfi}^{\sharp}$ such that $\wt U_{\vfi}^{\sharp}=U_g\oplus Z$, then $\restr{\epsi}{Z}$ is a linear isomorphism between $Z$ and $U_f^{\sharp}$. Note that $$\epsi(Z)=\epsi(U_{\vfi}^{\sharp})=U_f^{\sharp}.$$ Given $u\in U_g$, $z\in Z$, one can set $j(u+z)=u+\epsi(z)$, then $\pi\circ j=\restr{\epsi}{U_{\vfi}^{\sharp}}$. Next, if $T$ is a subspace of $V_n'$ satisfying $V_n'=U_f^{\sharp}\oplus T$, then, given $u\in U_g$, $\epsi(z)\in U_f^{\sharp}$, $t\in T$, one can set $i(u+\epsi(z)+t)=u+z+\alpha(t)$, where $\alpha\colon T\hookrightarrow V_m$ is an arbitrary embedding such that $U_{\vfi}^{\sharp}\cap\alpha(T)=0$. Clearly, $$\restr{(i\circ j)}{U_{\vfi}^{\sharp}}=\id_{U_{\vfi}^{\sharp}},$$ as required.

Thus, if we denote $U_g^{\sharp}=i(\wt U_g^{\sharp})\subset V_m$ then $\{U_g\subset U_g^{\sharp}\}$ is a flag in $V_m$ equipped with an isomorphism $U_g^{\sharp}/U_g\cong V_n'$. This isomorphism induces a surjection $g^{\sharp}\colon U_g^{\sharp}\twoheadrightarrow V_n'$ with kernel $U_g$. The strict standard extension $g\colon\Grr{k_n'}{V_n'}\hookrightarrow\Grr{k_m}{V_m}$ corresponding to the triple $(U_g,U_g^{\sharp},g^{\sharp})$ satisfies the property $g\circ f=\vfi$, as required.}

Note that the ind-varieties $\Gra{k}$ and $\Gra{\infty}$ considered in Example~\ref{exam:gen_flags} i) are exactly the standard nonisotropic ind-grassmannians defined above, so there is no abuse of notation. Similarly, all standard isotropic ind-grassmannians represent isomorphism classes of ind-varieties $\Grb{\Fo}{E}$ introduced in Example~\ref{exam:gen_flags_iso} i).

We are now ready to classify linear ind-grassmannians. The second main result of this section is as follows (see \cite[Theorem2]{PenkovTikhomirov1}).

\theop{Every\label{theo:class_grass} linear ind-grassmannian is isomorphic as an ind-variety to one of the standard ind-grassmannians $\Gra{k}$, $\Gra{\infty}$\textup, $\Grb{k}{\infty}$\textup, $\Grb{\infty}{\infty}$\textup, $\Grb{\infty}{k}$\textup, $\Grbi{\infty}{k}{0}$\textup, $\Grbi{\infty}{k}{1}$\textup, and the latter are pairwise non-isomorphic.}{Let $X$ be a linear ind-grassmannian given by the inductive limit of a chain of embeddings
$$X_1\stackrel{\vfi_1}\hookrightarrow X_2\stackrel{\vfi_2}\hookrightarrow\ldots\stackrel{\vfi_{n-1}}{\hookrightarrow}
X_n\stackrel{\vfi_n}\hookrightarrow X_{n+1}\stackrel{\vfi_{n+1}}\hookrightarrow\ldots,$$ where all $X_n$ are (possibly, orthogonal or symplectic) grassmannians and $\lim_{n\to\infty}\dim X_n=\infty$. Then for infinitely many $n$, $X_n$ is a grassmannian, or an orthogonal grassmannian, or a symplectic grassmannian. Hence we may assume without loss of generality that all $X_n$ are of one of these three types. Consider the case when all $X_n$ are grassmannians. (Two other cases can be considered similarly with some special features in the orthogonal type.)

There are two different options: either, for infinitely many $n$, the embedding $\vfi_n\colon X_n\hookrightarrow X_{n+1}$ factors through a projective space in $X_{n+1}$, or this is not the case. In the first case, clearly, $X\cong\Gra{1}\cong\Pp^{\infty}$. In the second case, by deleting some first embeddings we can assume that none of the embeddings $\vfi_n$ factors through a projective space in $X_{n+1}$. Thus, it follows from Theorem~\ref{theo:stand_exts} that all $\vfi_m$ are standard extensions, and, consequently, $X$ is isomorphic to $\Gra{k}$ or $\Gra{\infty}$.

For the proof of the fact that the standard ind-grassmannians are pairwise non-isomorphic, see \cite[Lemmas 5.1, 5.2, 5.4]{PenkovTikhomirov1}.}

The isotropic ind-grassmannians which do not appear in Theorem~\ref{theo:class_grass} are of the form $\Grb{F}{E}$ where $\beta$ is symmetric and $F$ is a maximal isotropic subspace such that $F=F^{\perp}$ or $F$ is an isotropic subspace which has codimension $2$ or $4$ in $F^{\perp}$. An ind-grassmannian $\Grb{F}{E}$ with $\codim_{F^{\perp}}F=4$ is a linear ind-grassmannian according to Definition~\ref{defi:linear_ind_grassmannian}. What \cite[Theorem2]{PenkovTikhomirov1} does not verify is that any linear ind-grassmannian $X=\ilm X_n$ such that $X_n$ is isomorphic to a grassmannian $\Grb{l_n-2}{V_n}$ for $m_n=\dim V_n=2l_n$, is isomorphic to $\Grb{F}{E}$ for $\codim_{F^{\perp}}F=4$. We nevertheless expect this to hold. The cases $F=F^{\perp}$ and $\codim_{F^{\perp}}F=2$ do not really fit Definition~\ref{defi:linear_ind_grassmannian} as the requirement $\Pic{X_n}\cong\Zp$ is violated. These cases deserve a special consideration.

We conclude this subsection by remarking that the idea of characterizing the ind-grassmannians\break $\Grr{F}{E}$ (or their isotropic counterparts) purely geometrically, in other words without reference to an action of $\GL_{\infty}(\Cp)$ on $\Grr{F}{E}$, could carry over to arbitrary ind-varieties of generalized flags. For this, one would need to define the notion of a strongly linear embedding of arbitrary usual flag varieties\break $X\hookrightarrow Y$ in purely geometric terms (i.e., strengthen Definition~\ref{defi:linear_morphism} in an appropriate way), and then prove that strongly linear embeddings are nothing but standard extensions (the notion of standard extension admits an obvious generalization to arbitrary flag varieties). This would then imply that any strongly linear flag ind-variety is an ind-variety of generalized flags.

\section{Ind-varieties of generalized flags as homogeneous ind-spaces}\label{sect:homo}\fakesect

In this section, we prove that all ind-varieties of generalized flags (possibly, isotropic) are homo\-ge\-neous ind-spaces of the classical ind-groups $$\GL_{\infty}(\Cp),~\SL_{\infty}(\Cp),~\SO_{\infty}(\Cp),~\Sp_{\infty}(\Cp)$$ defined in Section~\ref{sect:definitions}. Our exposition is based on the papers \cite{DimitrovPenkov1}, \cite{DimitrovPenkov2}, \cite{DanCohen1}, \cite{DanCohenPenkovSnyder1}, \cite{NeebPenkov1} and \cite{DanCohenPenkov1}.

\subsection{Classical ind-groups and their Cartan subgroups}\fakesst Recall\label{sst:Cartan_subgroups} the definitions of the ind-groups $\SL_{\infty}(\Cp)=\SL(V,E)$, $\SO_{\infty}(\Cp)=\SO(V,E,\beta)$ and $\Sp_{\infty}(\Cp)=\Sp(V,E,\beta)$ from Example~\ref{exam:ind_groups}. Let $G$ be one of these ind-groups. Note that in each case we have an exhaustion of $G$ by its finite-dimensional subgroups of the corresponding type described in Example~\ref{exam:ind_groups}. We denote this exhaustion by $G=\bigcup G_n$. For example, if $G=\SL_{\infty}(\Cp)$, then $G_n=\SL(V_n)\cong\SL_n(\Cp)$, etc. The goal of this subsection is to describe the structure of Cartan subgroups of $G$ in some detail.

Sometimes it is more convenient to work with Lie algebras instead of groups. Let $\gt_n\subset\glt(V_n)$ be the Lie algebra of $G_n$. To each linear operator $\vfi$ on $V_n$ one can assign a linear operator $\vfi'$ on $V_{n+1}$ by letting $$\vfi'(x)=x\text{ for $x\in V_n$, and $\vfi'(e)=0$ for $e\in E\setminus E_n$}$$ (cf. the definition of $\wh\vfi$ from Example~\ref{exam:ind_groups}). This gives an embedding $$\gt_n\hookrightarrow\gt_{n+1},~\vfi\mapsto\vfi'$$ for all $n\geq1$. We denote the inductive limit $\ilm\gt_n$ by $\gt$ and call it the \emph{Lie algebra} of $G$. For $G=\GL_{\infty}(\Cp)$, $\SL_{\infty}(\Cp)$, $\SO_{\infty}(\Cp)$ and $\Sp_{\infty}(\Cp)$ we write $\gt=\glt_{\infty}(\Cp)$, $\slt_{\infty}(\Cp)$, $\sot_{\infty}(\Cp)$ and $\spt_{\infty}(\Cp)$ respectively. If $H=\bigcup H_n$ is an ind-subgroup of $G$, where $H_n=H\cap G_n$, and $\htt_n$ is the Lie algebra of $H_n$, then the map $\vfi\mapsto\vfi'$ induces an embedding $$\htt_n\hookrightarrow\htt_{n+1},~\vfi\mapsto\vfi',$$ and the inductive limit $\htt=\ilm\htt_n$ is called the \emph{Lie algebra} of $H$.

In the finite-dimensional setting, there are several equivalent definitions of a Cartan subalgebra of a semisimple Lie algebra. For example, a Cartan subalgebra is a maximal toral subalgebra, or a nilpotent self-normalizing subalgebra. Moreover, all Cartan subalgebras of a finite-dimensional Lie algebra are conjugate. For $\gt$, the situation is somewhat more delicate. One new effect is that there exist maximal toral subalgebras of $\gt$ which do not yield a root decomposition of $\gt$. Another new effect is that there are non-conjugate maximal splitting toral subalgebras, and, as a consequence, different root systems of the same Lie algebra.

Recall that for finite-dimensional Lie algebras we have a notion of a Jordan decomposition, and
in particular of semisimple elements and nilpotent elements. Since these notions are preserved under embeddings $\gt_n\hookrightarrow\gt_{n+1}$, we can talk about a Jordan decomposition of an element of $\gt$. A subalgebra $\ttt$ of $\gt$ is called \emph{toral} if every element of $\ttt$ is semisimple. Every toral subalgebra is abelian. Given a subalgebra $\htt$ of $\gt$, denote by $\htt_{ss}$ the set of semisimple Jordan components of the elements of $\htt$. For example, if $\ttt$ is a toral subalgebra, then $\ttt_{ss}=\ttt$. A subalgebra $\htt$ of $\gt$ is called \emph{locally nilpotent} if it is locally nilpotent as $\htt$-module, i.e., if for all $x$, $y\in\htt$ there exists $m\in\Zp_{>0}$ such that $\ad{x}^m(y)=0$. This is equivalent to requiring that $\htt$ be a nested union of finite-dimensional nilpotent Lie algebras. A toral subalgebra $\ttt$ is called \emph{splitting} if $\gt$ is a weight $\ttt$-module, i.e., $$\gt=\bigoplus_{\lambda\in\htt^*}\gt^{\lambda},$$ where $$\gt^{\lambda}=\{x\in\gt\mid\ad{y}(x)=\lambda(y)x\text{ for all }y\in\htt\}.$$

\exam{i) Recall the\label{exam:gl_V_V_*} definition of $V_*$ from Section~\ref{sect:definitions}. The Lie algebra $\glt_{\infty}(\Cp)$ is isomorphic to the Lie algebra $\glt(V,V_*)=V\otimes V_*$ with the bracket induced by the product $$(v_1\otimes\alpha_1)(v_2\otimes\alpha_2)=\alpha_2(v_1)v_2\otimes\alpha_1.$$
An isomorphism $\eta\colon\glt(V,V_*)\to\glt_{\infty}(\Cp)=\glt(V)$ is given by the usual formula
$$\eta(v\otimes\alpha)(w)=\alpha(w)v,~v,~w\in V,~\alpha\in V_*.$$
Now, $$\ttt=\bigoplus_{n\geq0}\Cp e_n\otimes\Cp e_n^*$$ is a splitting maximal toral subalgebra of $\glt_{\infty}(\Cp)$. In fact, $\ttt$ consists of linear operators from $\glt_{\infty}(\Cp)$ which are diagonal in the basis $E$.

ii) Put $$\ttt=\bigoplus_{n\geq3}\Cp(e_n-e_1)\otimes\Cp(e_n^*-e_2^*).$$ Then $\ttt$ is a maximal toral subalgebra of $\glt_{\infty}(\Cp)$, and $\ttt$ is not splitting.}\newpage

For any subset $A\subset\gt$ and any subalgebra $\htt\subset\gt$ we define the \emph{centralizer} of $A$ in $\htt$ as the subspace $$\zt_{\htt}(A)=\{x\in\htt\mid[x,y]=0\text{ for all }y\in A\}\subset\htt.$$

\defi{A subalgebra $\htt$ of $\gt$ is a \emph{Cartan subalgebra} if $\htt$ is locally nilpotent and $\htt=\zt_{\gt}(\htt_{ss})$. An ind-subgroup $H$ of $G$ is a \emph{Cartan subgroup} if the Lie algebra of~$H$ is a Cartan subalgebra of $\gt$.}

\lemmp{Let $\htt$ be a locally nilpotent subalgebra of $\gt$. Then the following assertions hold\textup:
\begin{equation*}
\begin{split}
&\text{\textup{i)} $\htt\subset\zt_{\gt}(\htt_{ss})$};\\
&\text{\textup{ii)} $\htt_{ss}$ is a toral subalgebra of $\gt$};\\
&\text{\textup{iii)} $\zt_{\gt}(\htt_{ss})$ is a self-normalizing subalgebra of $\gt$}.
\end{split}
\end{equation*}}{i) Pick two elements $x$, $y\in\htt$. Since $\htt$ is locally nilpotent, $\ad{x}^m(y)=0$ for some $m$. Denote by $z_{ss}$ the semisimple part of an element $z\in\gt$. It is well known that $\ad{x_{ss}}$ is a polynomial in $\ad{x}$ with no constant term, hence $$\ad{x_{ss}}(\ad{x}^{m-1}(y))=0.$$ Each element commutes with its semisimple component, so $$\ad{x}^{m-1}(\ad{x_{ss}}(y))=0,$$ and it follows by induction that $\ad{x_{ss}}^m(y)=0$. Thus, $\ad{x_{ss}}(y)=0$, and, consequently, $\htt\subset\zt_{\gt}(\htt_{ss})$.

ii) Similarly, $\ad{y}(x_{ss})=0$ implies that $\ad{y_{ss}}(x_{ss})=0$. Therefore any two elements of $\htt_{ss}$ commute. Since the sum of any two commuting semisimple elements is semisimple, $\htt_{ss}$ is a subalgebra.

iii) Suppose that $z\in\gt$ belongs to the normalizer of $\zt_{\gt}(\htt_{ss})$, then $[x,y]\in\zt_{\gt}(\htt_{ss})$ for all $y\in\zt_{\gt}(\htt_{ss})$, in particular, for all $y\in\htt_{ss}$. Hence if $y\in\htt_{ss}$, then $[[x,y],y]=0$, and as $y$ is semisimple it follows that $[x,y]=0$. Thus, $x\in\zt_{\gt}(\htt_{ss})$, and so $\zt_{\gt}(\htt_{ss})$ is self-normalizing.}

The following theorem is the main general result of \cite{DanCohenPenkovSnyder1} characterizing Cartan subalgebras of $\gt$ (in \cite{DanCohenPenkovSnyder1} it is proved in more general context of locally reductive Lie algebras).

\mtheo{Let $\htt$ be a subalgebra of $\gt$. The following conditions on $\htt$ are equivalent\textup:
\begin{equation*}
\begin{split}
&\text{\textup{i)} $\htt$ is a Cartan subalgebra};\\
&\text{\textup{ii)} $\htt=\zt_{\gt}(\htt_{ss})$ and $\htt_{ss}$ is a subalgebra};\\
&\text{\textup{iii)} $\htt=\zt_{\gt}(\ttt)$ for some maximal toral subalgebra $\ttt$ of $\gt$\textup, and $\ttt=\htt_{ss}$}.
\end{split}
\end{equation*}
In addition\textup, if $\htt$ is a Cartan subalgebra\textup, then $\htt$ is self-normalizing\textup, and both the semisimple and nilpotent Jordan components of an element of $\htt$ belong to $\htt$.}

There is a unified description of Cartan subalgebras of $\gt$ in term of so-called self-dual systems, see \cite[Corollary 4.11]{DanCohenPenkovSnyder1} for the details. The Lie algebras $\glt_{\infty}(\Cp)$, $\slt_{\infty}(\Cp_{\infty})$ and $\spt_{\infty}(\Cp)$ admit only abelian Cartan subalgebras, while $\sot_{\infty}(\Cp)$ has non-abelian ones, see \cite[Subsection 4.2]{DanCohenPenkovSnyder1} for the example.

A Cartan subalgebra $\htt\subset\gt$ is called \emph{splitting} if $\htt_{ss}$ is a splitting toral subalgebra of $\gt$. In this case, $\htt=\htt_{ss}$. A Cartan subgroup $H$ of $G$ is called \emph{splitting} if its Lie algebra is a splitting Cartan subalgebra. In general, a Cartan subgroup (respectively, a Cartan subalgebra) of $G$ (respectively, of $\gt$) is splitting if and only if it a the inductive limit of Cartan subgroups of $G_n'$ (respectively, of Cartan subalgebras of $\gt_n'$), where $G=\bigcup G_n'$ is an exhaustion of $G$ by its finite-dimensional classical subgroups of corresponding type, and $\gt=\bigcup\gt_n'$ is the corresponding exhaustion of $\gt$ by its finite-dimensional classical Lie subalgebras.

\exam{\textbf{(Splitting maximal toral subalgebras)} i) Let $G=\GL_{\infty}(\Cp)$ or $\SL_{\infty}(\Cp)$ (res\-pec\-ti\-vely, $\gt=\glt_{\infty}(\Cp)$ or $\slt_{\infty}(\Cp)$). The\label{exam:Cartans} set $H$ (respectively,~$\htt$) of all linear operators from $G$ (respectively, from $\gt$) which are diagonal with respect to the basis $E$ is a splitting Cartan subgroup (respectively, a splitting Cartan subalgebra), and $\htt=\mathrm{Lie}\,H$.

ii) Put $m_n=\dim V_n$ for $n\geq1$. If $\beta$ is skew-symmetric then assume that each $m_n=2l_n$ is even (and so $\gt_n\cong\spt_{2l_n}(\Cp)$). If $\beta$ is symmetric then assume that either each $m_n=2l_n$ is even (and so $\gt_n\cong\sot_{2l_n}(\Cp)$), or that each $m_n=2l_n+1$ is odd (and so $\gt_n\cong\sot_{2l_n+1}(\Cp)$). If each $m_n$ is even, we renumerate the basis $E$ by letting $$E=\{e_i,~e_{-i}\}_{i\in\Zp_{>0}},~E_n=\{e_i,~e_{-i}\}_{i=1}^{l_n},$$ if each $m_n$ is odd, we renumerate $E$ by letting $$E=\{e_0\}\cup\{e_i,~e_{-i}\}_{i\in\Zp_{>0}},~E_n=\{e_0\}\cup\{e_i,~e_{-i}\}_{i=1}^{l_n}.$$ We may assume without loss of generality that
\begin{equation*}
\beta(u,v)=\begin{cases}
\sum\nolimits_{i=1}^{l_n}(u_iv_{-i}+u_{-i}v_i)&\text{for }\sot_{2l_n}(\Cp),\\
u_0v_0+\sum\nolimits_{i=1}^{l_n}(u_iv_{-i}+u_{-i}v_i)&\text{for }\sot_{2l_n+1}(\Cp),\\
\sum\nolimits_{i=1}^{l_n}(u_iv_{-i}-u_{-i}v_i)&\text{for }\spt_{2l_n}(\Cp).
\end{cases}
\end{equation*}
Here $u$, $v\in V_n$ and $x_i$ denotes the coordinate of a vector $x$ corresponding to $e_i$. Then the set $H$ (respectively, $\htt$) of all linear operators from $G$ (respectively, from $\gt$) which are diagonal with respect to the basis $E$ is a splitting Cartan subgroup (respectively, a splitting Cartan subalgebra), and $\htt=\mathrm{Lie}\,H$.}

The main result about splitting Cartan subalgebras is as follows \cite{DanCohenPenkovSnyder1}.
\mprop{\textup{i) }If $\ttt$ is a maximal splitting toral subalgebra of $\gt$\textup, then $\ttt=\zt_{\gt}(\ttt)$ is a splitting Cartan subalgebra of $\gt$. \textup{ii) }If $\htt$ is a splitting Cartan subalgebra of $\gt$\textup, then $\htt=\htt_{ss}$. \textup{iii) }For $\gt=\glt_{\infty}(\Cp)$\textup, $\slt_{\infty}(\Cp)$ and $\spt_{\infty}(\Cp)$\textup, all splitting Cartan subalgebras are conjugate by the group $\Aut\gt$. For $\gt=\sot_{\infty}(\Cp)$\textup, there exists exactly two $\Aut\gt$-conjugacy classes of splitting Cartan subalgebras corresponding to the exhaustions $\gt=\ilm\sot_{2l_n}(\Cp)$ and $\gt=\ilm\sot_{2l_n+1}(\Cp)$ described above.}

Note also that $G\subsetneq\Aut\gt$. Indeed, $G=\SL(V,E)$ consists of automorphisms of $V$ with determinant~$1$ which keep all but finitely many elements of $E$ fixed, while $\Aut\slt_{\infty}(\Cp)$ in this case contains the group of all automorphisms of $V$ which induce automorphisms on $V_*$.

\subsection{Splitting Borel and parabolic subgroups of classical ind-groups}\fakesst In\label{sst:Borel_parabolic_subgroups} this subsection, we discuss Borel and parabolic subgroups and subalgebras of classical ind-groups and their Lie algebras res\-pec\-tively. In particular, we classify all splitting Borel subgroups in terms of their roots.

In the finite-dimensional setting, a Borel subalgebra of a semisimple Lie algebra is a maximal solvable subalgebra, and a parabolic subalgebra is a subalgebra containing some Borel subalgebra. We say that a subalgebra $\bt$ of $\gt$ is \emph{locally solvable} if every finite subset of $\bt$ is contained in a finite-dimensional solvable subalgebra, i.e. if $\bt$ is a union of its finite-dimensional solvable subalgebras.

\defi{A \emph{Borel subalgebra} $\bt$ of $\gt$ is a maximal locally solvable subalgebra of $\gt$. A Borel subalgebra $\bt$ of $\gt$ is called \emph{splitting} if $\bt$ contains a splitting Cartan subalgebra of $\gt$. A \emph{Borel subgroup} of $G$ is an ind-subgroup $B$ such that the Lie algebra $\bt=\Lie B$ is a Borel subalgebra of $\gt$. A Borel subgroup $B$ is called \emph{splitting} if its Lie algebra is a splitting Borel subalgebra, or, equivalently, if $B$ contains a splitting Cartan subgroup of $G$.}

In general, Borel subalgebras may have properties which are very unusual from the finite-di\-men\-si\-onal point of view. For instance, there exists a Borel subalgebra of $\glt_{\infty}(\Cp)$ (see Example~\ref{exam:Borel_without_toral} below) which contains no nonzero semisimple elements, and hence no non-trivial toral subalgebras!

In contrast, splitting Borel subalgebras have a very nice description presented below. Note that a Borel subalgebra of $\gt$ (respectively, a Borel subgroup of $G$) is splitting if and only if it a the inductive limit of Borel subgroups of $G_n'$ (respectively, of Borel subalgebras of~$\gt_n'$), where $G=\bigcup G_n'$ is an exhaustion of the ind-group $G$ by its finite-dimensional classical subgroups of corresponding type, and $\gt=\bigcup\gt_n'$ is the corresponding exhaustion of $\gt$ by its finite-dimensional classical Lie sub\-al\-gebras.

Let $\htt$ be a splitting Cartan subalgebra of $\gt$ described in Examples~\ref{exam:Cartans} i) and ii). Any splitting Borel subalgebra is conjugate via $\mathrm{Aut}\,\gt$ to a splitting Borel subalgebra containing $\htt$. Therefore, in what follows we restrict ourselves to considering only Borel subalgebras $\bt$ which contains $\htt$. We have a root decomposition $$\gt=\htt\oplus\bigoplus_{\alpha\in\Phi}\gt^{\alpha}$$ where $\Phi$ is \emph{the root system of} $\gt$ with respect to $\htt$, and $\gt^{\alpha}$ are the \emph{root spaces}. The root system $\Phi$ is simply the union of the root systems of $\gt_n$, and equals one of the following infinite root systems:
\begin{equation*}
\begin{split}
A_{\infty}&=\pm\{\epsi_i-\epsi_j,~i,j\in\Zp_{>0},~i<j\},\\
B_{\infty}&=\pm\{\epsi_i-\epsi_j,~i,j\in\Zp_{>0},~i<j\}\\
&\cup\pm\{\epsi_i+\epsi_j,~i,j\in\Zp_{>0},~i<j\}\cup\pm\{\epsi_i,~i\in\Zp_{>0}\},\\
C_{\infty}&=\pm\{\epsi_i-\epsi_j,~i,j\in\Zp_{>0},~i<j\}\\
&\cup\pm\{\epsi_i+\epsi_j,~i,j\in\Zp_{>0},~i<j\}\cup\pm\{2\epsi_i,~i\in\Zp_{>0}\},\\
D_{\infty}&=\pm\{\epsi_i-\epsi_j,~i,j\in\Zp_{>0},~i<j\}\\
&\cup\pm\{\epsi_i+\epsi_j,~i,j\in\Zp_{>0},~i<j\}.
\end{split}
\end{equation*}
The linear functions $\epsi_i-\epsi_j$, $\epsi_i+\epsi_j$, $\epsi_i$, $2\epsi_i$ on $\htt$ are defined as follows: given $h\in\htt$,
\begin{equation*}
\begin{split}
(\epsi_i-\epsi_j)(h)&=\begin{cases}h_{i,i}-h_{j,j}&\text{in case }A_{\infty},\\
h_{i,i}-h_{-i,-i}-h_{j,j}+h_{-j,-j}&\text{otherwise},\\
\end{cases}\\
(\epsi_i+\epsi_j)(h)&=h_{i,i}-h_{-i,-i}+h_{j,j}-h_{-j,-j},\\
\epsi_i(h)&=2(h_{i,i}-h_{-i,-i}),\\
2\epsi_i(h)&=h_{i,i}-h_{-i,-i}.
\end{split}
\end{equation*}
Here, in the orthogonal and symplectic cases, we enumerate the basis vectors from $E$ as in Example~\ref{exam:Cartans} and denote by $x_{i,j}$ the $(i,j)$th element of a matrix $x$.

Recall \cite{DimitrovPenkov3} that a linear order on $\{0\}\cup\{\pm\epsi_i\}$ is $\Zp_2$-\emph{linear} if multiplication by $-1$ reverses the order. By \cite[Proposition 3]{DimitrovPenkov3}, there exists a bijection between splitting Borel subalgebras of $\gt$ containing $\htt$ and certain linearly ordered sets as follows:
\begin{equation*}\predisplaypenalty=0
\begin{split}
&\text{for }A_{\infty}\text{: linear orders on }\{\epsi_i\};\\
&\text{for }B_{\infty}\text{ and }C_{\infty}\text{: }\Zp_2\text{-linear orders on }\{0\}\cup\{\pm\epsi_i\};\\
&\text{for }D_{\infty}\text{: }\Zp_2\text{-linear orders on }\{0\}\cup\{\pm\epsi_i\}\text{ with the property that}\\
&\text{a minimal positive element (if it exists) belongs to }\{\epsi_i\}.
\end{split}
\end{equation*}
In the sequel we denote these linear orders by $\prec$. To write down the above bijection, denote $\teta_i=\epsi_i$, if $\epsi_i\succ0$, and $\teta_i=-\epsi_i$, if $\epsi_i\prec0$ (for $A_{\infty}$, $\teta_i=\epsi_i$ for all $i$). Then put $\bt=\htt\oplus\nt$, where $$\nt=\bigoplus\limits_{\alpha\in\Phi^+}\gt^{\alpha}$$ and, by definition,
\begin{equation*}\predisplaypenalty=0
\begin{split}
A_{\infty}^+&=\{\teta_i-\teta_j,~i,j\in\Zp_{>0},~\teta_i\succ\teta_j\},\\
B_{\infty}^+&=\{\teta_i-\teta_j,~i,j\in\Zp_{>0},~\teta_i\succ\teta_j\}\\
&\cup\{\teta_i+\teta_j,~i,j\in\Zp_{>0},
~\teta_i\succ\teta_j\}\cup\{\teta_i,~i\in\Zp_{>0}\},\\
C_{\infty}^+&=\{\teta_i-\teta_j,~i,j\in\Zp_{>0},~\teta_i\succ\teta_j\}\\
&\cup\{\teta_i+\teta_j,~i,j\in\Zp_{>0},~\teta_i\succ\teta_j\}\cup\{2\teta_i,~i\in\Zp_{>0}\},\\
D_{\infty}^+&=\{\teta_i-\teta_j,~i,j\in\Zp_{>0},~\teta_i\succ\teta_j\}\\
&\cup\{\teta_i+\teta_j,~i,j\in\Zp_{>0},~\teta_i\succ\teta_j\}.\\
\end{split}
\end{equation*}

A subalgebra $\pt\subset\gt$ is called \emph{parabolic} (respectively, \emph{splitting parabolic}) if it contains a Borel (respectively, a splitting Borel) subalgebra of $\gt$. An ind-subgroup $P\subset G$ is called a \emph{parabolic} (re\-spec\-ti\-ve\-ly, a \emph{splitting parabolic}) \emph{subgroup} if it contains a Borel (respectively, a splitting Borel) subgroup of $G$, or, equivalently, if the Lie algebra $\pt=\Lie P$ is a parabolic (respectively, a splitting parabolic) subalgebra of~$\gt$.

Note that a parabolic subalgebra of $\gt$ (respectively, a parabolic subgroup of $G$) is splitting if and only if it a the inductive limit of parabolic subgroups of $G_n'$ (respectively, of parabolic subalgebras of~$\gt_n'$), where $G=\bigcup G_n'$ is an exhaustion of the ind-group $G$ by its finite-dimensional classical subgroups of corresponding type, and $\gt=\bigcup\gt_n'$ is the corresponding exhaustion of $\gt$ by its finite-dimensional classical Lie sub\-al\-gebras.


\subsection{Homogeneous ind-spaces}\fakesst Here\label{sst:homo_ind_spaces} we establish the main result of this section which claims that each splitting parabolic subgroup ia s stabilizer of a generalized flag, and, vice versa, each ind-variety of (isotropic) generalized flags is a homogeneous ind-space of the group $G$.

Let $H$ and $\htt$ be as in Example~\ref{exam:Cartans}. In the isotropic case, we will use all the notation from Example~\ref{exam:Cartans}~ii). Note that in this case $E$ is a $\beta$-isotropic basis with respect to the involution $$i_E\colon e_i\mapsto e_{-i},~e_i\in E.$$
Each splitting parabolic subalgebra $\pt$ of $\gt$ is conjugate to a splitting parabolic subalgebra of $\gt$ containing $\htt$ via $\Aut\gt$, so we may assume without loss of generality that $\pt$ contains $\htt$ (and, consequently, the corresponding parabolic subgroup $P$ of $G$ contains $H$).

Let $\Fo$ be a generalized flag in $V$ compatible with $E$ (and $\beta$-isotropic whenever $E$ is $\beta$-isotropic). The ind-group $G$ naturally acts on the ind-variety $\Fl(\Fo,E)$ (and on $\Fl(\Fo,\beta,E)$ in the isotropic case). Denote by $P_{\Fo}$ the stabilizer of $\Fo$ in $G$.

\theop{\textup{i)} The subgroup\label{theo:splitting_parabs_flags} $P_{\Fo}$ is a splitting parabolic subgroup of $G$ containing $H$. \textup{ii)} The map $\Fo\mapsto P_{\Fo}$ is a bijection between the set of generalized flags in $V$ compatible with the basis $E$ and the set of splitting parabolic subgroups of $G$ containing $H$.}{i) The inclusion $H\subset P_{\Fo}$ follows immediately from the definition of $H$ and the compatibility of $\Fo$ and $E$. Since each $P_n=P\cap G_n$ is the stabilizer of the flag $\Fo\cap V_n$, $P_n$ is a parabolic subgroup of $G_n$. Hence $P=\ilm P_n$ is a splitting parabolic subgroup of the ind-group $G$.

ii) Conversely, let $P=\ilm P_n$ be a parabolic subgroup of $G$ containing $H$, where $P_n$ is a parabolic subgroup of $G_n$ for $n\geq1$. Denote by $\Fo(n)$ the flag in $V_n$ whose stabilizer coincides with $P_n$. Then $\iota_n(\Fo(n))=\Fo(n+1)$ (and $\iota_n^{\beta}(\Fo(n))=\Fo(n+1)$ in the isotropic case). The inductive limit $\ilm\Fo(n)$ defines a generalized flag $\Fo$ in $V$, and it is straightforward to check that $P=P_{\Fo}$.}

Note that maximal generalized flags in $V$ compatible with $E$ correspond to splitting Borel subgroups of $G$ containing $H$ under the above bijection.

Since $G/P_{\Fo}=\bigcup G_n/P_n$, where $P_n=P_{\Fo}\cap V_n$, we conclude that $G/P_{\Fo}$ is a locally projective ind-variety, as we mentioned in Section~\ref{sect:definitions}. We are now ready to endow $\Fl(\Fo,E)$ and $\Fl(\Fo,\beta,E)$ with respective structures of homogeneous ind-spaces.

\theop{There is an isomorphism of ind-varieties $\Fl(\Fo,E)\cong G/P_{\Fo}$ \textup(and of ind-varieties $\Fl(\Fo,\beta,E)\cong G/P_{\Fo}$ in the isotropic case\textup).}{Given $\Go\in\Fl(\Fo,E)$ or $\Go\in\Fl(\Fo,\beta,E)$, let $U$ be a finite-dimensional subspace of $V$ whose existence is provided by the $E$-commensurability of $\Fo$ and $\Go$. We may assume without loss of generality that $U=V_n$ for some $n\geq1$. Then $\Fo(n)=\Fo\cap V_n$ and $\Go(n)=\Go\cap V_n$ are flags of the same type in the finite-dimensional vector space $V_n$, hence there exists $g_n\in G_n$ such that $g_n(\Fo(n))=\Go(n)$. We can extend $g_n$ to an element $g_{n+1}=\wh g_n$ by letting $\wh g_n(e)=e$ for $e\in E\setminus E_n$, etc. Let $g$ be the corresponding element of $G$. Then $$\eta\colon\Fl(\Fo,E)\to G/P_{\Fo}\text{ (or $\eta\colon\Fl(\Fo,\beta,E)\to G/P_{\Fo}$)},~\Go\mapsto gP,$$ is a well-defined map. One can easily check that $\eta$ is an isomorphism of ind-varieties.}

\subsection{Borel and parabolic subalgebras: general case}\fakesst In\label{sst:Borel_parabolic_general} this subsection we briefly discuss a description of (possibly, non-splitting) Borel and parabolic subalgebras of $\gt$ (or, equivalently, Borel and parabolic subgroups of $G$) in terms of so-called closed generalized flags and taut pairs of semiclosed generalized flags. This material is taken from the papers \cite{DanCohen1} and \cite{DanCohenPenkov1}.

If a parabolic subgroup of $G$ (or, equivalently, a parabolic subalgebra of $\gt$) is not splitting, it can not be the stabilizer of a generalized flag compatible with the basis $E$ of $V$ or with any other $G$-eligible basis $E'$ of $V$. Thus, in order to relate general non-splitting parabolic subgroups and subalgebras to generalized flags, we should consider generalized flags which do not admit a $G$-eligible compatible basis.

Recall the identification $$\glt(V,V_*)\cong\glt_{\infty}(\Cp)=\glt(V)$$ from Example~\ref{exam:gl_V_V_*}. Under this isomorphism, $\slt_{\infty}(\Cp)$ is identified with the commutator subalgebra $\slt(V,V_*)$ of $\glt(V,V_*)$. Note that if $U$ is a countable-dimensional complex vector space and $$\langle\cdot,\cdot\rangle\colon V\times U\to\Cp$$ is a non\-de\-ge\-ne\-rate pairing, then we can set $\glt(V,U)\cong\glt_{\infty}(\Cp)$ to be the Lie algebra $V\otimes U$ with the bracket induced by the product $$(v_1\otimes u_1)(v_2\otimes u_2)=\langle v_1,u_2\rangle v_2\otimes u_1.$$ Then $\glt(V,V_*)$ is a particular case of this construction, where nondegenerate paring $V\times V_*\to\Cp$ is given by $\langle v,\alpha\rangle=\alpha(v)$.

Now, if $\beta$ is symmetric (respectively, skew-symmetric) nondegenerate bilinear form on $V$, then $\beta$ defines a non\-de\-ge\-ne\-rate pairing $V\times V\to\Cp$, and we can identify $\sot_{\infty}(\Cp)$ (respectively, $\spt_{\infty}(\Cp)$) with the Lie subalgebra $\sot(V,V)=\bigwedge^2V$ (respectively, $\spt(V,V)=\mathrm{Sym}^2\,V$) of the Lie algebra $\glt(V,V)$.

Given a nondegenerate pairing $\langle\cdot,\cdot\rangle\colon V\times U\to\Cp$ and a subspace $F$ of $V$ or $U$, we denote by $F^{\perp}$ the $\langle\cdot,\cdot\rangle$-orthogonal complement of $F$ in $U$ or $V$. (If the pairing is given by a nondegenerate symmetric or skew-symmetric bilinear form $\beta$, this definition coincides with the definition of $F^{\perp}$ given is Section~\ref{sect:definitions}.) To each chain $\Co$ of subspaces of $V$ one can assign the chain $\Co^{\perp}=\{F^{\perp},~F\in\Co\}$ of subspaces of $U$, and vice versa.

A subspace $F\subset V$ is called \emph{closed} (in the Mackey topology on $V$) if $F=\overline{F}$, where $\overline{F}=F^{\perp\perp}$ is\break called the \emph{closure} of $F$. A generalized flag $\Fo=\{F_{\alpha}',~F_{\alpha}''\}_{\alpha\in\Au}$ in $V$ is called \emph{semiclosed} if $$\overline{F_{\alpha}'}\in\{F_{\alpha}',~F_{\alpha}''\}$$ for all $\alpha\in\Au$. A semiclosed generalized flag $\Fo$ is \emph{closed} if, in addition, $\overline{F_{\alpha}''}=F_{\alpha}''$ for all $\alpha\in\Au$. Note that if a generalized flag $\Fo$ in $V$ is weakly compatible with the basis $E$ defining $V_*$, then $\Fo$ is automatically closed.

Note that each of spaces $V$ and $U$ is a $\glt(V,U)$-module. Hence $\glt(V,U)$ naturally acts on chains\break in $V$ and $U$, and the stabilizer $\Stab{\,}{\Fo}$ of a generalized flag $\Fo$ of $V$ in $\glt(V,U)$ is given by formula $$\Stab{\,}{\Fo}=\sum_{\alpha\in\Au}F_{\alpha}''\otimes(F_{\alpha}')^{\perp}.$$
If $\gt=\sot(V,V)$ or $\spt(V,V)$, then we write $\St_{\Fo}^{\gt}=\Stab{\,}{\Fo}\cap\gt$.

The following result describes Borel subalgebras of classical infinite-dimensional simple Lie algebras (or, equivalently, Borel subgroups of classical ind-groups), see \cite[Theorems 4.3, 4.10, 4.16]{DanCohen1}.

\mtheo{\textup{i) }Suppose $\gt=\glt(V,V_*)$ \textup(respectively\textup, $\spt(V,V)$ or $\sot(V,V)$\textup). A subalgebra $\bt$ of $\gt$ is a Borel subalgebra if and only if $\bt$ is a stabilizer of a maximal closed \textup(respectively\textup, of a maximal closed isotropic\textup) generalized flag in $V$. \textup{ii) }For $\gt=\glt(V,V_*)$ \textup(respectively\textup, $\gt=\spt(V,V)$\textup{),} the map $\Fo\mapsto\Stab{\,}{\Fo}$ \textup(respectively\textup, $\Fl\mapsto\St_{\Fo}^{\gt}$\textup) from the set of maximal closed \textup(respectively\textup, the set of maximal closed isotropic\textup) generalized flags in $V$ to the set of Borel subalgebras of $\gt$ is bijective. For $\gt=\sot(V,V)$\textup, a fiber of the map $\Fo\mapsto\St_{\Fo}^{\gt}$ from the set of maximal closed isotropic flags in $V$ to the set of Borel subalgebras of $\gt$ contains at most two elements.}

(An explicit description of fibers of the latter map is given in \cite[Subsection 4.2]{DanCohen1}.)

\defi{We say that two semiclosed generalized flags $\Fo$ and $\Go$ in $V$ and $U$ respectively form a \emph{taut pair} if the chain $\Fo^{\perp}$ (respectively, $\Go^{\perp}$) is stable under $\Stab{\,}{\Go}$ (respectively, under $\Stab{\,}{\Fo}$). Given a nondegenerate symmetric or skew-symmetric form $\beta$ on $V$, we say that a semiclosed generalized flag $\Fo$ is \emph{self-taut} if $\Fo^{\perp}$ is stable under the stabilizer of $\Fo$ in $\glt(V,V)$ (i.e., if the generalized flag $\Fo$ form a taut pair with itself).}

To each taut pair $\Fo$, $\Go$ one can assign a subalgebra $$\St_{\Fo,\Go}=\Stab{\,}{\Fo}\cap\Stab{\,}{\Go}$$ and a certain subalgebra $(\St_{\Fo,\Go})_-$ defined in \cite[p. 23]{DanCohenPenkov1}. If $\Fo$ is a self-taut generalized flag in $V$ and $\gt=\sot(V,V)$ or $\spt(V,V)$, then we write $$(\St_{\Fo}^{\gt})_-=(\St_{\Fo,\Fo})_-\cap\gt.$$ The following result describes parabolic subalgebras of classical infinite-dimensional simple Lie algebras (or, equivalently, parabolic subgroups of classical ind-groups), see \cite[Theorems 5.6, 6.6]{DanCohenPenkov1}.

\mtheo{\textup{i) }Let $\gt=\glt(V,V_*)$ or $\slt(V,V_*)$\textup, and $\pt$ be a subspace of $\gt$. Then $\pt$ is a parabolic subalgebra of $\gt$ if and only if there exists a \textup(unique\textup) taut couple $\Fo$\textup, $\Go$ such that $$(\St_{\Fo,\Go})_-\subset\pt\subset\St_{\Fo,\Go}.$$ \textup{ii) }Let $\gt=\sot(V,V)$ or $\spt(V,V)$\textup, and $\pt$ be a subspace of $\gt$. Then $\pt$ is a parabolic subalgebra of $\gt$ if and only if there exists a \textup(unique\textup) self-taut generalized flag $\Fo$ in $V$ such that $$(\St_{\Fo}^{\gt})_-\subset\pt\subset\St_{\Fo}^{\gt}.$$}

\exam{Let\label{exam:Borel_without_toral} $V=\langle e_{\alpha},~\alpha\in\Qp\rangle_{\Cp}$, $U=\langle f_{\beta},~\beta\in\Qp\rangle_{\Cp}$, and
\begin{equation*}
\langle e_{\alpha},f_{\beta}\rangle=\begin{cases}
1,&\text{if }\alpha>\beta,\\
0,&\text{if }\alpha\leq\beta.
\end{cases}
\end{equation*}
Then $\langle\cdot,\cdot\rangle$ is a nondegenerate pairing, so $$\glt(V,U)\cong\glt_{\infty}(\Cp).$$ For $\alpha\in\Qp$, put $$F_{\alpha}'=\langle e_{\gamma},~\gamma<\alpha\rangle_{\Cp},~F_{\alpha}''=\langle e_{\gamma},~\gamma\leq\alpha\rangle_{\Cp}.$$ Then $\Fo=\{F_{\alpha}',~F_{\alpha}''\}_{\alpha\in\Qp}$ is a maximal closed generalized flag in $V$. Its stabilizer $\bt$ in $\glt(V,U)$ is a Borel subalgebra of $\glt(V,U)$ which contains no nonzero semisimple elements and, consequently, no nonzero toral subalgebras.}

Note that, for a non-splitting parabolic subgroup $P$ of $G$, $G/P$ can be endowed with an ind-variety structure. These ind-varieties have not been yet explored. The main difference with the splitting case is that $P\cap G_n$ is not necessarily a parabolic subgroup in $G_n$, hence $G/P$ is not exhausted by compact varieties. The role of arbitrary nonsplitting Borel subgroups and subalgebras in representation theory also remains unclear.

\section{Schubert decomposition}\label{sect:Schubert_calculus}\fakesect

In the finite-dimensional setting, Schubert decomposition plays a central role in the study of the geometry of flag varieties.

Recall the definition of $G_n$ from Section~\ref{sect:definitions}. Fix a maximal torus $H_n$ of the group $G_n$, a Borel subgroup $B_n$ of $G_n$ containing $H_n$, and a parabolic subgroup $P_n$ of $G_n$ containing $B_n$. Let $G_n/P_n$ be the corresponding flag variety. Denote by $N_{G_n}(H_n)$ the normalizer of $H_n$ in $G_n$. Then $$W_n=N_{G_n}(H_n)/H_n$$ is the \emph{Weyl group} of $G_n$. Since $H_n$ and $B_n$ are fixed, we have a set of \emph{simple generators} of $W_n$, and, consequently, a \emph{length function} $\ell_n$ and a \emph{Bruhat order} $\leq_n$ on $W_n$. The details see for instance in \cite{Bourbaki} or \cite{BjornerBrenti1}.

Given $w\in W_n$, we denote by $\dot{w}$ an arbitrary representative of $w$ in $N_{G_n}(H_n)$. Let $\Fo_n\in G_n/P_n$ be the flag in $V$ whose stabilizer $\Stab{G_n}{\Fo_n}$ in $G_n$ is $P_n$. Denote by $W_{P_n}$ the parabolic subgroup of $W_n$ corresponding to $P_n$, and by $W^{P_n}$ the set of minimal length representatives in the right cosets of $W_{P_n}$ in $W_n$ ($W^{P_n}$ is in bijection with $W_{P_n}\backslash W_n$). Then for the $W_n$-action on $G_n/P_n$ we have $$w\Go=\dot{w}(\Go),~\Go\in G_n/P_n.$$ In what follows we write $$\Fo_w=w\Fo_n,~w\in W_n.$$ The description of $B_n$-orbits on $G_n/P_n$ is given by the following \emph{Schubert decomposition}: $$G_n/P_n=\bigsqcup_{w\in W^{P_n}}B_n\Fo_w.$$ Furthermore, each \emph{Schubert cell} $X_w^{\circ}=B_n\Fo_w$ is isomorphic to the affine space $\Ap^{\ell_n(w)}$, and, given $\sigma$, $\tau\in W_n$, the Schubert cell $X_{\sigma}^{\circ}$ is contained in the \emph{Schubert subvariety} $X_{\tau}$ (by definition, $X_{\tau}$ is the closure of $X_{\tau}^{\circ}$ in $G_n/P_n$) if and only if $\sigma\leq_n\tau$. The \emph{Bruhat decomposition} of $G_n$ claims that $$G_n=\bigsqcup_{w\in W^{P_n}}B_n\dot{w}P_n.$$

In this section we show how these classical results can be extended to the case of ind-varieties of generalized flags. This material is based on the paper \cite{FressPenkov1}.

\subsection{Analogues of the Weyl group}\fakesst In this subsection we present some combinatorial results which are analogues to the Weyl group combinatorics used in usual Schubert calculus. Denote\label{sst:analogues_Weyl_group} by $W=W(E)$ the group of permutations of $E$ which fix all but finitely many elements of $E$. In the isotropic case we assume in addition that each $w\in W$ commutes with the involution $i_E$. Note that $$W=\ilm W_n,$$ where $H$ is the splitting Cartan subgroup of $G$ consisting of all diagonal operators from $G$ with respect to $E$, $H_n=H\cap G_n$, and the embedding $W_n\hookrightarrow W_{n+1}$ is induced by the embedding $H_n\hookrightarrow H_{n+1}$. It is also not difficult to check that $$W\cong N_G(H)/H,$$ where $N_G(H)$ denotes the normalizer of $H$ in the ind-group $G$.

Next, let $B$ and $P$ be a splitting Borel and a splitting parabolic subgroup of $G$ containing $H$ respectively. (We do not assume that $B$ is conjugate to a subgroup of $P$!) For brevity, denote by $\Fl$ the ind-variety of generalized flags $\Fl(\Fo,E)$ (or $\Fl(\Fo,\beta,E)$ in the isotropic case). As we know from Theorem~\ref{theo:splitting_parabs_flags}, there exists a unique generalized flag $\Go\in\Fl$ such that $$P=P_{\Go}=\Stab{G}{\Go}.$$ We may assume without loss of generality that $\Go=\Fo$, i.e., that $P=P_{\Fo}$.

Recall the linearly ordered set $\Au$ from Section~\ref{sect:definitions}. Set $\So$ to be the set $\So(E,\Au)$ (respectively, the set $\So(E,\beta,\Au)$) to be of all surjective maps from $E$ to~$\Au$ (respectively, all surjective maps $\sigma$ from $E$ to $\Au$ satisfying $\sigma\circ i_E=i_{\Au}\circ\sigma$). To each $\sigma\in\So$ we assign the generalized flag $$\Fo_{\sigma}=\{\Fo_{\sigma.\alpha}',~\Fo_{\sigma,\alpha}'',~\alpha\in\Au\},$$ where $$\Fo_{\sigma,\alpha}'=\langle e\in E,~\sigma(e)<\alpha\rangle_{\Cp},~\Fo_{\sigma,\alpha}''=\langle e\in E,~\sigma(e)\leq\alpha\rangle_{\Cp}.$$ In this way, $\{\Fo_{\sigma},~\sigma\in\So\}$ are all generalized flags from $\Fl$ compatible with the basis $E$. (One should apply Lemma~\ref{lemm:iso_flags} to check this in the isotropic case.)

Let $\sigma_0\in\So$ be the unique surjection such that $\Fo=\Fo_{\sigma_0}$. Then $\sigma_0$ defines a partial order $\leq_P$ on the basis $E$ by letting $$e\leq_Pe'\text{ if }\sigma_0(e)\leq\sigma_0(e').$$ This partial order has the property that the relation ``$e=e'$ or $e$ is not comparable with $e'$'' is an equivalence relation on $E$. In fact, fixing a splitting parabolic subgroup $P$ of $G$ containing $H$ is equivalent to fixing a partial order $\leq_P$ on $E$ with this property. Moreover, $P$ is a splitting Borel subgroup if and only if the order $\leq_P$ is linear.

We say that a pair $(e,e')\in E\times E$ is an \emph{inversion} for a surjection $\sigma\in\So$ if $e<_Be'$ and $\sigma(e)>\sigma(e')$. (In the isotropic case we assume in addition that $e<_Bi_E(e')$ and $e'\neq i_E(e')$.) Note that the group $W$ acts on $\So$ by $$w\cdot\sigma=\sigma\circ w^{-1},$$ and that if $\sigma$ belongs to the $W$-orbit of $\sigma_0$, then the condition $\sigma(e)>\sigma(e')$ is equivalent to the condition $w(e)>_Pw(e')$ where $\sigma=w^{-1}\cdot\sigma_0$. The \emph{inversion number} of $\sigma\in\So$ is $$n(\sigma)=n_B^P(\sigma)=\#\{(e,e')\in E\times E\mid (e,e')\text{ is an inversion of }\sigma\}.$$ Of course, $n(\sigma)$ can be infinite.

The inversion number cannot be directly interpreted as Bruhat length because we do not assume~$B$ to be conjugate to a subgroup of $P$. Nevertheless, we put $$\wh E=\{(e,e')\in E'\times E'\mid e\neq e'\},$$ where, by definition,
\begin{equation*}
E'=\begin{cases}E&\text{for $\GL_{\infty}(\Cp)$ and $\SL_{\infty}(\Cp)$},\\
\{e\in E\mid e\neq i_E(e)\}&\text{for $\SO_{\infty}(\Cp)$ and $\Sp_{\infty}(\Cp)$.}
\end{cases}
\end{equation*}
Let $t_{e,e'}$ be the permutation of $E$ such that $$t_{e,e'}(e)=e',~t_{e,e'}(e')=e$$ and $t_{e,e'}(e'')=e''$ for all other $e''\in E$. Set \begin{equation*}
s_{e,e'}=\begin{cases}t_{e,e'}\circ t_{i_E(e),i_E(e')},&\text{if $e'\neq i_E(e)$ in the isotropic case},\\
t_{e,e'}&\text{otherwise}.
\end{cases}
\end{equation*}
Clearly, $\{s_{e,e'},~(e,e')\in\wh E\}$ is a set of generators of $W$.

Now we denote $$S_B=\{s_{e,e'}\mid\text{ $e$, $e'$ are consecutive elements of the partial ordered set $(E',\leq_B)$}\}.$$ In general, $S_B$ does not generate $W$. For an element $w\in W$, we define
\begin{equation*}
\ell(w)=\ell_B(w)=\begin{cases}
\min\{l\geq0\mid w=s_1\ldots s_l\text{ for some }s_1,~\ldots,~s_l\in S_B\},&\text{ if such an $l$ exists},\\
\infty&\text{otherwise}.
\end{cases}
\end{equation*}
Note that $W_n$ can be considered as the subgroup of $W$ generated by the set $$S_B^n=\{s_{e,e'}\mid\text{ $e$, $e'$ are consecutive elements of $(E_n\cap E',\leq_B)$}\},$$ which is a set of simple generators of $W_n$. Let $\ell_n$ be the corresponding length function on $W_n$.

\mprop{Let $w\in W$.\label{prop:l_B_w} Then
\begin{equation*}
\begin{split}
&\text{\textup{i)} $\ell(w)=\lim_{n\to\infty}\ell_n(w)$};\\
&\text{\textup{ii)} $\ell(w)=n_B^B(w^{-1}\cdot\sigma_0)$};\\
&\text{\textup{iii)} $\ell(w)=\infty$ if and only if there exists $e\in E$}\\
&\hphantom{\text{\textup{iii)} }}\text{such that the set $\{e'\in E\mid e<_Be'<_Bw(e)\}$ is infinite}.\\
\end{split}
\end{equation*}}
See \cite[Proposition 8]{FressPenkov1} for the proof of this proposition.
\corop{The following conditions\label{coro:l_B_w} are equivalent\textup{:}
\begin{equation*}
\begin{split}
&\text{\textup{i)} $S_B$ generates $W$};\\
&\text{\textup{ii)} $\ell(w)<\infty$ for all $w\in W$};\\
&\text{\textup{iii)} $(E,\leq_B)$ is isomorphic as an ordered set to a subset of $\Zp$}.\\
\end{split}
\end{equation*}}{The equivalence {\rm i)}$\Leftrightarrow${\rm ii)} is immediate. Note that condition {\rm iii)} is equivalent to requiring that, for all $e,~e'\in E$, the interval $\{e''\in E\mid e\leq_Be''\leq_Be'\}$ is finite. Thus the implication {\rm iii)}$\Rightarrow${\rm ii)} is guaranteed by Proposition~\ref{prop:l_B_w} iii). Conversely, if {\rm ii)} holds true, then we get $\ell(s_{e,e'})<\infty$ for all $(e,e')\in\wh E$, whence, again by Proposition~\ref{prop:l_B_w} iii), the set $\{e''\in E\mid e\leq_Be''\leq_Be'\}$ is finite. This implies {\rm iii)}.}

Let $\dot w$ denote a representative of $w$ in $N_G(H)$.
\propp{Let $w\in W$. Then\label{prop:conj_B_P}
\begin{equation*}
\begin{split}
&\text{\textup{i)} $B\subset\dot wP\dot w^{-1}$ if and only if $n_B^P(w\cdot\sigma_0)=0$};\\
&\text{\textup{ii)} there exists $w\in W$ such that $B\subset\dot wP\dot w^{-1}$ if and only if}\\
&\hphantom{\text{\textup{ii)} }}\text{there exists $w\in W$ such that $n_B^P(w^{-1}\cdot\sigma_0)<\infty$}.\\
\end{split}
\end{equation*}}{i) By the definition of the generalized flag $\Fo_{\sigma_0}$, the condition $B\subset P$ is equivalent to the condition that the linear order $\leq_B$ on $E$ refines the partial order $\leq_P$ on $E$, i.e., that $e\leq_Pe'$ implies $e\leq_Be'$ for all $e$, $e'\in E$. The latter condition is equivalent to the condition that the map $\sigma_0$ is nondecreasing, i.e., that $e\leq_Be'$ implies $\sigma_0(e)\leq\sigma_0(e')$ for all $e$, $e'\in E$. Since $$\dot wP\dot w^{-1}=\Stab{G}\Fo_{w\cdot\sigma_0},$$ part i) follows.

ii) Part i) implies that if $B\subset\dot wP\dot w^{-1}$, then $n_B^P(w^{-1}\cdot\sigma_0)<\infty$. For the proof of the remaining implication, see \cite[Proposition 9]{FressPenkov1}.}

\subsection{Schubert decomposition}\fakesst Denote\label{sst:Schubert_decomposition} by $W_P$ the subgroup of $W$ consisting of all $\sigma\in W$ satisfying $w\cdot\sigma_0=\sigma_0$. It is straightforward to check that the map $w\mapsto\Fo_{w\cdot\sigma_0}$ induces a bijection between the left coset space $W/W_P$ and the set of $E$-compatible generalized flags from $\Fl$.

We now define a partial order on $\So$ analogous to the Bruhat order. Let $\sigma$, $\tau\in\So$. We write $\sigma\to\tau$ if there exists $(e,e')\in\wh E$ such that $$e<_Be',~\sigma(e)<\sigma(e')$$ and $\tau=\sigma\circ s_{e,e'}$. We set $\sigma<\tau$ if there exist $k\geq1$ and elements $\tau_1$, $\ldots$, $\tau_k\in\So$ such that $$\sigma\to\tau_1\to\ldots\to\tau_k=\tau.$$

Given a generalized flag $$\Go=\{G_{\alpha}',~G_{\alpha}'',~\alpha\in\Au\}\in\Fl,$$ we define a map $\sigma_{\Go}\colon E\to\Au$ which measures the relative position of $\Go$ to the maximal generalized flag $\Fo_0$, where $B=\Stab{G}{\Fo_0}$. Namely, for an arbitrary $e\in E$, set $$\sigma_{\Go}(e)=\min\{\alpha\in\Au\mid G_{\alpha}''\cap F_{0,e}''\neq G_{\alpha}''\cap F_{0,e}'\}.$$
Here $\Fo_0=\{F_{0,e}',~F_{0,e}'',~e\in E\}$ and
$$F_{0,e}'=\langle e'\in E\mid e'<_Be\rangle_{\Cp},~F_{0,e}''=\langle e'\in E\mid e'\leq_Be\rangle_{\Cp}.$$
It can be checked directly that $\sigma_{\Go}$ belongs to $\So$; furthermore, $\sigma\in W\sigma_0$ where $W\sigma_0=\{w\cdot\sigma_0,~w\in W\}$ denotes the $W$-orbit of $\sigma_0$.

We are now ready to formulate the main result of this section. Given a generalized flag $\Go\in\Fl$, we denote by $B\Go$ its $B$-orbit under the natural action of $B$ on $\Fl$. Let $W^P$ be a set of representatives of the left coset space $W/W_P$.
\theop{Let $P=P_{\Fo}$ and\label{theo:Schubert} let $B$ be any splitting Borel subgroup of $G$ containing $H$. Then
\begin{equation*}
\begin{split}
&\text{\textup{i)} }G/P=\Fl=\bigsqcup_{\sigma\in W\sigma_0}B\Fo_{\sigma}=\bigsqcup_{w\in W^P}B\Fo_{w\cdot\sigma_0};\\
&\text{\textup{ii)} given $\sigma\in W\sigma_0$\textup, a generalized flag $\Go\in\Fl$ belongs to $B\Fo_{\sigma}$}\\
&\text{$\hphantom{\text{\textup{ii) }}}$if and only if $\sigma_{\Go}=\sigma$};\\
&\text{\textup{iii)} for $\sigma\in W\sigma_0$\textup, the orbit $B\Fo_{\sigma}$ is a locally closed ind-subvariety of $\Fl$}\\
&\text{$\hphantom{\text{\textup{iii) }}}$isomorphic to the affine space $\Ap^{n_P^B(\sigma)}$};\\
&\text{\textup{iv)} for $\sigma,~\tau\in W\sigma_0$\textup, the inclusion $B\Fo_{\sigma}\subset\overline{B\Fo}_{\tau}$ holds if and only if $\sigma\leq\tau$.}
\end{split}
\end{equation*}}{Consider the case of $\Fl(\Fo,E)$; the case of $\Fl(\Fo,\beta,E)$ can be considered similarly.

i) The statement follows from the finite-dimensional Schubert decomposition of $\Fl_n=\Fl(d_n,V_n)$, where $d_n$ is the type of the flag $\Fo(n)=\Fo\cap V_n=\{F\cap V_n,~F\in\Fo\}$. Indeed, if $n$ is large enough so that the flag $\Go(n)=\Go\cap V_n$ belongs to $\Fl_n$, then the $B_n$-orbit of $\Go(n)$ contains a unique element of the form $\Fo_{w\cdot\sigma_0}\cap V_n$ with $w\in W_n$.

ii) Let $\Go\in\Fl$.
According to part i), there is a unique $\sigma\in W\sigma_0$ such that $\Go\in B\Fo_{\sigma}$, say $\Go=b(\Fo_{\sigma})$, where $b\in B$. Hence
$$G_{\alpha}''\cap F_{0,e}'=b(F_{\sigma,\alpha}''\cap F_{0,e}')\text{ and }G_{\alpha}''\cap F_{0,e}''=b(F_{\sigma,\alpha}''\cap F_{0,e}'')\text{ for all }e\in E,~\alpha\in\Au,$$
because $F_{0,e}'$ and $F_{0,e}''$ are $b$-stable. This implies $\sigma_{\Go}=\sigma_{\Fo{_\sigma}}$.
Moreover, from the definition of $\Fo_{\sigma}$ we see that $F_{\sigma,\alpha}''\cap F_{0,e}''\neq F_{\sigma,\alpha}''\cap F_{0,e}'$ if and only if $\sigma(e)\leq\alpha$. Whence $$\sigma(e)=\min\{\alpha\in\Au\mid F_{\sigma,\alpha}''\cap F_{0,e}''\neq F_{\sigma,\alpha}''\cap F_{0,e}'\}=\sigma_{\Fo_{\sigma}}(e)\text{ for all }e\in E.$$ Thus, $\sigma_{\Go}=\sigma$. The equality $\sigma_{\Go}=\sigma$ guarantees in particular that $\sigma_{\Go}\in\So$.

iii) The statement follows again from the finite-dimensional case. Note that, given $w\in W_n$, the image of the Schubert cell $B_n(\Fo_{w\cdot\sigma_0}\cap V_n)$ by the embedding $\iota_n$ is an affine subspace of the Schubert cell $B_{n+1}(\Fo_{w\cdot\sigma_0}\cap V_{n+1})$.

iv) We consider $\sigma,~\tau\in W\sigma_0$, $\sigma\leq\tau$, and let $n\geq1$ be such that $\Fo_{\sigma}(n)$ and $\Fo_{\tau}(n)$ are contained in $\Fl_n$. We may assume without loss of generality that $\sigma\to\tau$, i.e., that $\tau=\sigma\circ s_{e,e'}$ for a pair $(e,e')\in E\times E$ with $$e<_Be',~\sigma(e)<\sigma(e').$$ Up to choosing $n$ larger if necessary, we may assume that $e,~e'\in E_n$. Then, by finite-dimensional results, we get $$B_n\Fo_{\sigma}(n)\subset\overline{B_n\Fo_{\tau}(n)}.$$ Therefore,
$B\Fo_{\sigma}\subset\overline{B\Fo}_{\tau}$. Conversely, assume that $\Fo_{\sigma}\in\overline{B\Fo}_{\tau}$. Then $$\Fo_{\sigma}(n)\in\overline{B_n\Fo_{\tau}(n)}$$ for $n\geq1$ large enough. Once again, by finite-dimensional results, this yields $\sigma\leq\tau$. The proof is complete.}

In what follows, we call $X_{\sigma}^{\circ}=B\Fo_{\sigma}$ and $X_{\sigma}=\overline{X_{\sigma}^{\circ}}$ the \emph{Schubert cell} and the \emph{Schubert subvariety} of~$\Fl$ corresponding to $\sigma\in\So$ respectively.

The following corollary (the Bruhat decomposition of the ind-group $G$) is an immediate consequence of the Schubert decomposition of $\Fl$ established by the theorem above. Note that, in general, $B$ is not conjugate to a subgroup of $P$, so in fact one has many different Bruhat decompositions of $G$ depending on the choice of $B$ and $P$.

\mcoro{\textup{\textbf{(Bruhat decomposition of the ind-group $G$)} }Let $G$ be one of the ind-groups $\GL_{\infty}(\Cp)$\textup, $\SL_{\infty}(\Cp)$\textup, $\SO_{\infty}(\Cp)$ and $\Sp_{\infty}(\Cp)$\textup, and $P$ and $B$ be respectively a splitting parabolic and a splitting Borel subgroup containing $H$. Then we have a decomposition $$G=\bigsqcup_{w\in W^P}B\dot wP.$$}

In general, the $B$-orbits in Theorem~\ref{theo:Schubert} are infinite dimensional. The following two results determine the situations in which finite-dimensional orbits appear.

\theop{Let $G$ be\label{theo:fin_dim_B_orbits_Schubert} one of the ind-groups $\GL_{\infty}(\Cp)$\textup, $\SL_{\infty}(\Cp)$\textup, $\SO_{\infty}(\Cp)$ and $\Sp_{\infty}(\Cp)$\textup, and $P$ and $B$ be respectively a splitting parabolic and a splitting Borel subgroup containing $H$. The following conditions are equivalent\textup:
\begin{equation*}
\begin{split}
&\text{\textup{i)} $B$ is conjugate via $G$ to a subgroup of $P$};\\
&\text{\textup{ii)} at least one $B$-orbit on $G/P$ is finite dimensional};\\
&\text{\textup{iii)} one $B$-orbit on $G/P$ is a point \textup(and this orbit is necessarily unique\textup)}.\\
\end{split}
\end{equation*}}{Condition i) means that there exists $g\in G$ such that $B\subset gPg^{-1}$, or, equivalently, such that the element $gP\in G/P$ is fixed by $B$, i.e., that $G/P$ comprises a $B$-orbit reduced to a single point. We proved the equivalence {\rm i)}$\Leftrightarrow${\rm iii)}. The implication {\rm iii)}$\Rightarrow${\rm ii)} is immediate, while
the implication {\rm ii)}$\Rightarrow${\rm i)} follows from Proposition~\ref{prop:conj_B_P} and Theorem~\ref{theo:Schubert}.}

\corop{Assume $P\neq G$. The following conditions are equivalent\textup:
\begin{equation*}
\begin{split}
&\text{\textup{i)} $B$ is conjugate via $G$ to a subgroup of $P$, and $\Fo_0$ is a flag};\\
&\text{\textup{ii)} every $B$-orbit on $G/P$ is finite dimensional}.\\
\end{split}
\end{equation*}}{The implication {\rm i)}$\Rightarrow${\rm ii)} is a consequence of Theorem~\ref{theo:fin_dim_B_orbits_Schubert}, Corollary~\ref{coro:l_B_w}, Theorem~\ref{theo:Schubert} and the following fact: if there exists $w_0\in W$ such that $n_B^P(w_0\cdot\sigma_0)=0$ then, for all $w\in W$, $$n_B^P(w^{-1}\cdot\sigma_0)=\inf_{w'\in W_P}l(w_0w'w).$$ (See \cite[Proposition 10]{FressPenkov1} for the proof of this fact.)

Suppose that {\rm ii)} holds. By Theorem~\ref{theo:fin_dim_B_orbits_Schubert}, there exists $g\in G$ such that $B\subset gPg^{-1}$, so we may assume without loss of generality that $B\subset P$. Arguing by contradiction, assume that $\Fo_0$ is not a flag, i.e., $(E,\leq_B)$ is not isomorphic to a subset of $\Zp$. Then there exist $e,~e'\in E$ such that the set $$\{e''\in E\mid e<_Be''<_Be'\}$$ is infinite. Since the surjective map $\sigma_0\in\So$ corresponding to $P$ is nondecreasing (see the proof of Proposition\ref{prop:conj_B_P}) and nonconstant (because $P\neq G$), we find $\wh e,\wh e'$ for which $\wh e\leq_Be<_Be'\leq_B\wh e'$ and $\sigma_0(\wh e)<\sigma_0(\wh e')$. Thus, $\dim B\Fo_{w\cdot\sigma_0}=\infty$ for $w=s_{\wh e,\wh e'}^{-1}$ by Theorem~\ref{theo:Schubert}, a~contradiction.}

\exam{\textbf{(Schubert decomposition of ind-grassmannians)} Let $G=\GL_{\infty}(\infty)$\textup, $\SL_{\infty}(\Cp)$. Consider the case of $\Grr{F}{E}$, i.e. let $$\Fo=\{\{0\}\subset F\subset V\}.$$ Here $\Au=\{1,2\}$, and if $\Fo$ is compatible with $E$, then the surjective map $\sigma_0\colon E\to\Au$ for which $F=\langle e\in E\mid\sigma_0(e)=1\rangle_{\Cp}$ can be simply viewed as the subset $\sigma_0\subset E$ such that $F=\langle\sigma_0\rangle_{\Cp}$. Hence we can identify $\So$ with the set $\So(E)$ of subsets of $E$, and this identification is $W$-equivariant. Note that $$W\sigma_0=\{\sigma\in\So(E)\mid|\sigma\setminus\sigma_0|=|\sigma_0\setminus\sigma|<\infty\}.$$ We will write $F_{\sigma}=\langle\sigma\rangle_{\Cp}$ for $\sigma\in\So(E)$.}\newpage

i) First, suppose that $\dim F=k<\infty$. Then $\Grr{F}{E}\cong\Gra{k}$. Denote by $\So_k(E)$ the set of all subsets of $E$ of cardinality $k$. According to Theorem~\ref{theo:Schubert}, $$\Grr{F}{E}=\bigsqcup_{\sigma\in\So_k(E)}B\Fo_{\sigma}.$$ The cell $X_{\sigma}^{\circ}=B\Fo_{\sigma}$ is finite-dimensional if and only if $\sigma$ is contained in a finite ideal of the ordered set $(E,\leq_B)$. It follows that there are finite-dimensional $B$-orbits on $\Grr{F}{E}$ if and only if the maximal flag $\Fo_0$ corresponding to $B$ contains a $k$-dimensional subspace. By Theorem~\ref{theo:fin_dim_B_orbits_Schubert}, in the latter case $B$ is conjugate to a subgroup of $\Stab{G}{\Fo}$. Moreover, all cells $X_{\sigma}^{\circ}$ are finite dimensional if and only if $(E,\leq_B)$ is isomorphic to $\Zp_{>0}$ as an ordered set, i.e., if and only if $\Fo_0$ has the form $$\Fo_0=\{F_{0,0}\subset F_{0,1}\subset\ldots\}$$ with $\dim F_{0,i}=i$.

ii) Next, suppose that $\codim_VF=k<\infty$. Again, $\Grr{F}{E}\cong\Gra{k}$. If, as above, $$\Fo_0=\{F_{0,0}\subset F_{0,1}\subset\ldots\},$$ where $\dim F_{0,i}=i$ for all $i$, and $B=\Stab{G}{\Fo_0}$, then Theorem~\ref{theo:fin_dim_B_orbits_Schubert} implies that all $B$-orbits on $\Grr{F}{E}$ are infinite dimensional. However, if $B'$ is the stabilizer of a maximal generalized flag containing a subspace $F'$ such that the flag $\{\{0\}\subset F'\subset V\}$ is $E$-commensurable with~$\Fo$, then there exist finite-dimensional $B'$-orbits on $\Grr{F}{E}$. Moreover, there is no Borel subgroup which has only finite-dimensional orbits on both ind-grassmannians from i) and ii).

iii) Finally, suppose that both $\dim F$ and $\codim_VF$ are infinite, so $\Grr{F}{E}\cong\Gra{\infty}$. Assume that the basis $E$ is parameterized by $\Zp$. We consider the splitting Borel subgroup $B\subset G$ corresponding to the natural order on $\Zp$. If $F=\langle e_i,~i\leq0\rangle_{\Cp}$ then $B\subset\Stab{G}{\Fo}$, hence every $B$-orbit on $\Grr{F}{E}$ is finite dimensional. On the other hand, if $F=\langle e_{2i},~i\in\Zp\rangle_{\Cp}$, then every $B$-orbit on $\Grr{F}{E}$ is infinite dimensional.

\subsection{Smoothness of Schubert subvarieties}\fakesst In\label{sst:smoothness_Schubert} this subsection we study the smoothness of Schubert subvarieties $X_{\sigma}=\overline{B\Fo}_{\sigma}$ of the ind-variety $\Fl$, where $\Fl=\Fl(\Fo,E)$ or $\Fl=\Fl(\Fo,\beta,E)$. The general principle is straightforward: $X_{\sigma}$ is smooth if and only if its intersections with suitable finite-dimensional flag subvarieties of $\Fl$ are smooth.

For classical finite-dimensional groups there is a remarkable characterization of smooth Schubert subvarieties of flag varieties in terms of pattern avoidance (see, e.g., \cite[Chapter 8]{BilleyLakshmibai1}). For example, if $G_n=\SL(V_n)$ (an hence the Weyl group $W_n$ is isomorphic to the symmetric group $S_n$), then, given $\sigma\in W_n$, the Schubert subvariety $X_{\sigma}$ of the flag variety $G_n/B_n$ is smooth if and only if $\sigma$ \emph{avoids the patterns} 3412 \emph{and} 4231, i.e., there are no $i,~j,~k,~l$, $q\leq i<j<k<l\leq\dim V_n$, such that
\begin{equation*}
\sigma(k)<\sigma(l)<\sigma(i)<\sigma(j)~\text{or }\sigma(l)<\sigma(j)<\sigma(k)<\sigma(i).
\end{equation*}

The notion of smooth point of an ind-variety is given in Subsection~\ref{sst:ind_varieties}. It is known \cite{Kumar1} that if $X$ is and-variety with exhaustion $\bigcup_{n\geq1}X_n$ by its finite-dimensional subvarieties, $x\in X$, and there exists a subsequence $\{X_{n_k}\}_{k\geq1}$ such that $x$ is a smooth point of $X_{n_k}$ for all $k\geq1$, then $x$ is a smooth point of $X$. For example, $\Ap^{\infty}$ and $\Pp^{\infty}$ are smooth ind-varieties.

The converse of the statement above is clearly false. For instance, for every $n$ fix an embedding $\Ap^n\hookrightarrow\Ap^{n+1}$. Pick a point $x\in\Ap^1$, and for each $n\geq1$ let $X_n'\subset\Ap^{n+1}$ be an $n$-dimensional subspace of $\Ap^{n+1}$ containing $x$ and distinct of $\Ap^n$. Now, if we set $X_n=X_n'\cup\Ap^n$, then the subvarieties $X_n$ exhaust the smooth ind-variety $\Ap^{\infty}$ but $x$ is a singular point of every $X_n$. Nevertheless, the following partial converse is true: if each embedding $X_n\hookrightarrow X_{n+1}$ has a left inverse in the category of algebraic varieties, and $x\in X$ is a singular point of $X_n$ for some $n\geq1$, then $x$ is a singular point of $X$ \cite[Lemma 6]{FressPenkov1}.

The singularity criterion given below requires a technical assumption on $B$ and $\Fo_{\sigma}$. We assume that at least one of the following conditions holds: $\Fo_0$ is a flag, or $\Fo_{\sigma}$ is a flag and $\dim F_{\sigma,\alpha}''/F_{\sigma,\alpha}'<\infty$ whenever $$\{0\}\neq F_{\sigma,\alpha}'\subset F_{\sigma_{\alpha}}''\neq V.$$ For example, this holds for ind-grassmannians. Recall the notion of $\Fl_n$ and $\Fl_n^{\beta}$ from Section~\ref{sect:definitions}; the intersection $X_{\sigma,n}=X_{\sigma}\cap\Fl_n$ (respectively, $X_{\sigma,n}=X_{\sigma}\cap\Fl_n^{\beta}$) is a Schubert subvariety of $\Fl_n$ (respectively, of $\Fl_n^{\beta}$) in the usual sense. By $\Sing{X}$ we denote the set of singular point of a variety or an ind-variety $X$.

\mtheo{Let $G$ be\label{theo:smoothness_Schubert} one of the ind-groups $\GL_{\infty}(\Cp)$\textup, $\SL_{\infty}(\Cp)$\textup, $\SO_{\infty}(\Cp)$ and $\Sp_{\infty}(\Cp)$\textup, and $B$\textup, $\Fo_{\sigma}$\textup, $\Fo_0$ be as above. Then exactly one of the two following statements holds\textup:
\begin{equation*}
\begin{split}
&\text{\textup{i)} the variety $X_{\sigma,n}$ is smooth for all $n$\textup, and the ind-variety $X_{\sigma}$ is smooth};\\
&\text{\textup{ii)} there exists $n_0\geq1$ such that the variety $X_{\sigma,n}$ is singular for all $n\geq n_0$},\\
&\text{$\hphantom{\text{\textup{ii) }}}$and the ind-variety $X_{\sigma}$ is singular with $\Sing{X_{\sigma}}=\bigcup_{n\geq n_0}X_{\sigma,n}$}.
\end{split}
\end{equation*}}

\noindent See \cite[Theorem 4]{FressPenkov1} for the proof of this theorem. We do not know whether Theorem~\ref{theo:smoothness_Schubert} is valid in general, i.e. without conditions on $\Fo_0$ and $\Fo_{\sigma}$.

\exam{It turns out that the smoothness criteria for Schubert subvarieties of finite-dimensional flag varieties in terms of pattern avoidance may pass to the limit at infinity. For instance, let $\Fo=\Fo_{\sigma}$ be a maximal generalized flag compatible with $E$. In this case we have two linear orders on $E$: the first one $\leq_B$ corresponds to the splitting Borel subgroup $B$, while the second one $\leq_P$ corresponds to the splitting parabolic (in fact, Borel) subgroup $P=P_{\Fo}$, i.e., $$F_e'=\langle e'\in E\mid e'<_Pe\rangle_{\Cp},~F_e''=\langle e'\in E\mid e'\leq_Pe\rangle_{\Cp}.$$

From Theorem~\ref{theo:Schubert} we know that the Schubert ind-subvarieties $X_{\sigma}$ of $\Fl(\Fo,E)$ are parameterized by the elements of $W\sigma_0$, where $\sigma_0\colon E\to\Au$ is the surjection corresponding to $\Fo=\Fo_{\sigma_0}$. Since $\Fo$ is maximal, $\Au=E$ and $\sigma_0$ is a bijection, so if $\sigma\in W\sigma_0$, then $\sigma\in W$. We now also that, given $\sigma\in W$, $$\dim X_{\sigma}=n_B^P(\sigma)=\#\{(e,e')\in E\mid e<_Be',~\sigma(e)>_P\sigma(e')\}.$$

It follows from Theorem~\ref{theo:smoothness_Schubert} and from the finite-dimensional criterion that if $\Fo_0$ is a flag, or $\Fo$ is a flag and $\dim F_{\sigma,\alpha}''/F_{\sigma,\alpha}'<\infty$ whenever $\{0\}\neq F_{\sigma,\alpha}'\subset F_{\sigma_{\alpha}}''\neq V$, then the Schubert ind-subvariety $X_{\sigma}$ is singular if and only if there exist $e_1$, $e_2$, $e_3$, $e_4\in E$ such that $e_1<_Be_2<_Be_3<_Be_4$ and $$\sigma(e_3)<_P\sigma(e_4)<_P\sigma(e_1)<_P\sigma(e_2)\text{ or }\sigma(e_4)<_P\sigma(e_2)<_P\sigma(e_3)<_P\sigma(e_1).$$ In particular, if the basis $E$ comprises infinitely many pairwise disjoint quadruples $e_1$, $e_2$, $e_3$, $e_4\in E$ such that $e_1<_Be_2<_Be_3<_Be_4$ and, say, $\sigma(e_3)<_P\sigma(e_4)<_P\sigma(e_1)<_P\sigma(e_2)$, then for every $\sigma$ the Schubert ind-subvariety $X_{\sigma}$ is singular. Hence there exist pairs $(B,\Fo)$ for which all Schubert ind-subvarieties of the corresponding ind-variety of generalized flags are singular.}

\subsection{Concluding remarks}\fakesst The\label{sst:concl_remarks_Schubert} main difference with the finite-dimensional case is that an ind-variety of generalized flags $G/P$ admits many non-conjugate Schubert de\-com\-po\-si\-tions. This is a consequence of the fact that the Borel subgroups $B$, whose orbits on $G/P$ form Schubert de\-com\-po\-si\-tions, are non-conjugate under the group of automorphisms of $G$. In particular, the infinite-dimensional projective space admits a Schubert decomposition with all finite-dimensional cells, as well as a Schubert decomposition with all infinite-dimensional cells.

We note that in the finite-dimensional case, in addition to the group-theoretic point of view on Schubert decomposition which we adopt here, there is also the purely geometric approach of fixing a reference maximal flag and studying in what ways a varying flag (in other words, a point on $G_n/P_n$) differs from the reference flag. This approach is also valid in the case we consider: the reference generalized flag (possibly, containing $G/P$ as a subchain) can be chosen to be $B$-stable. The resulting theory is equivalent to the one presented above.

\section{Finite-rank vector bundles}\label{sect:alg_geom}\fakesect

In this section we consider two related topics. First, we discuss an infinite-dimensional analogue of the Bott--Borel--Weil theorem. We also establish projectivity criterion for $G/P$: the ind-variety $G/P$ is projective if and only if $P=P_{\Fo}$, where $\Fo$ is a flag. Next, the Barth--Van de Ven--Tyurin--Sato theorem claims that any finite rank vector bundle on $\Pp^{\infty}$ is isomorphic to a direct sum of line bundles. We present a generalization of this result to a large class of ind-varieties, in particular to the ind-grassmannian $\Gra{\infty}$. The material of this section is taken from the papers \cite{DimitrovPenkovWolf1}, \cite{PenkovTikhomirov2}.

\subsection{Bott--Borel--Weil theorem}\fakesst We first recall\label{sst:Borel_Bott_Weil} the classical Bott--Borel--Weil theorem. Here we assume that $G_n$ is a connected simply-connected simple algebraic group, $B_n$ is a Borel subgroup and $P_n$ is a parabolic subgroup containing $B_n$. More precisely, $G_n=\SL_n(\Cp)$, $\Spin_n(\Cp$) or $\Sp_{2n}(\Cp)$. (The group $\Spin_n(\Cp)$ is the universal cover of $\SO_n(\Cp)$, and in this subsection we replace $\SO_n(\Cp)$ by $\Spin_n(\Cp)$.) It is a well-known theorem that the Picard group $\Pic{G_n/B_n}$ is identified with the lattice of integral weights of $\gt_n$ via the correspondence $\lambda\mapsto\Ou(\lambda)$, where $\Ou(\lambda)$ is the line bundle induced from the character of $B_n$ whose differential restricted to the Lie algebra $\htt_n$ of a maximal torus $H_n$ contained in $B_n$ coincides with $\lambda$. In particular, every line bundle on $G_n/B_n$ admits a canonical $G_n$-linearization.

For an integral weight $\lambda\in\htt_n^*$, consider the weight $\lambda+\rho$, where $\rho$ is the half-sum of the positive roots of $G_n$ with respect to $B_n$. If $\lambda+\rho$ is regular, then there is a unique element $w_{\lambda}$ in the Weyl group $W_n$ of $G_n$ such that $w_{\lambda}(\lambda+\rho)$ is dominant. The following theorem computes the cohomology groups $H^q(G_n/B_n,\Ou(-\lambda))$ for all integral weights $\lambda$ (the consideration of $\Ou(-\lambda)$ instead of $\Ou(\lambda)$ is convenient in the isomorphism (\ref{formula:BBW_fin}) below).

\mtheo{\textup{\textbf{(Bott--Borel--Weil)}} Let $\lambda$ be\label{theo:cohom_fin_B_n} an integral weight. If $\lambda+\rho$ is not regular\textup, then $$H^q(G_n/B_n,\Ou(-\lambda))=0\text{ for $q=0,~\ldots,~\dim G_n/B_n$},$$ i.e.\textup, the sheaf corresponding to the line bundle $\Ou(\lambda)$ is acyclic. If $\lambda+\rho$ is regular\textup, then $$H^q(G_n/B_n,\Ou(-\lambda))=0\text{ for $q\neq l(w_{\lambda})$},$$ where $l(w_{\lambda})$ is the length of $w_{\lambda}$ with respect to the simple roots of $B_n$. For $q=l(w_{\lambda})$\textup, there is a canonical $G_n$-isomorphism
\begin{equation}
H^q(G_n/B_n,\Ou(-\lambda))\cong V(\mu)^*,\label{formula:BBW_fin}
\end{equation}
where $\mu=w(\lambda+\rho)-\rho$ and $V(\mu)$ is the simple $G_n$-module with highest weight $\mu$.}

This theorem allows to compute also the cohomology of any finite-rank vector bundle on $G_n/P_n$ which is simple as a linearized $G_n$-vector bundle, i.e., is induced from a simple finite-dimensional $P_n$-module. Note that any such bundle has the form $p_*\Ou(\lambda)$ for an appropriate weight $\lambda$, where $${p\colon G_n/B_n\to G_n/P_n}$$ is the canonical epimorphism.\newpage

\mtheo{For any $\lambda$ such that $p_*\Ou(-\lambda)\neq0$, there is a canonical\label{theo:cohom_fin_P_n} isomorphism of $G_n$-modules $$H^q(G_n/B_n,\Ou(-\lambda))\cong H^q(G_n/P_n,p_*\Ou(-\lambda))$$ for any $q\geq0$.}

Theorem~\ref{theo:cohom_fin_B_n} was proved by R. Bott in \cite{Bott1} (the case $q=0$ was due to A. Borel and A. Weil). An alternative proof (which we strongly recommend to the reader) was given by M. Demazure in \cite{Demazure1}, and was then shortened in \cite{Demazure2}.

We now turn our attention to the infinite-dimensional situation. Let $G$ be the inductive limit of $G_n$, i.e. $G=\SL_{\infty}(\Cp)$, $\Spin_{\infty}(\Cp)$, $\Sp_{\infty}(\Cp)$. In particular, the ind-group $\Spin_{\infty}(\Cp)$ is well defined as there are natural injective homomorphisms $\Spin_n(\Cp)\hookrightarrow\Spin_{n+1}(\Cp)$ for $n\geq3$.

Note that the Picard group of $G/P$ is naturally isomorphic to the projective limit $\plm\Pic{G_n/P_n}$ of Picard groups. Moreover, the groups $\Pic{\Fl(\Fo,E)}$ and $\Pic{\Fl(\Fo,\beta,E)}$ are naturally isomorphic to the respective groups of integral weights of the Lie algebra of the ind-group $P=P_{\Fo}$.

Indeed, it is well known that $$\Pic{(G_n/P_n)}\cong\Hom{P_n,\Cp^{\times}},$$ where $\Hom{P_n,\Cp^{\times}}$ denotes the group of morphisms from $P_{\Fo}$ to the multiplicative group of $\Cp$. Immediate verification shows that the diagram
\begin{equation*}
\xymatrix{
\Pic{(G_{n+1}/P_{n+1})}\ar[d]\ar[r]^{\cong}&\Hom{P_{n+1},\Cp^{\times}}\ar[d]\\
\Pic{(G_n/P_n)}\ar[r]^{\cong}&\Hom{P_n,\Cp^{\times}}
}
\end{equation*}
is commutative. Thus,
\begin{equation*}
\begin{split}
\Pic{\Fl(\Fo,E)}&=\plm\Pic{G_n/P_n}\\
&=\plm\Hom{P_n,\Cp^{\times}}=\Hom{P_{\Fo},\Cp^{\times}},
\end{split}
\end{equation*}
and the latter group is nothing but the group of integral weights of the Lie algebra of $P=P_{\Fo}$. (See \cite[Proposition 7.2]{DimitrovPenkov1} and \cite[Proposition 15.1]{DimitrovPenkovWolf1} for the details.)

Let us now specialize to the case $P=B$ for an arbitrary splitting Borel subgroup $B$ of $G$ containing a fixed Cartan subgroup $H=\ilm H_n$, where $H_n$ is a Cartan subgroup of $G_n$ for all $n\geq1$. Let $\htt$ be the Lie algebra of $H$. A weight $\lambda\in\htt^*$ is \emph{integral} if, for any $n\geq1$, its restriction $\restr{\lambda}{\htt_n}$ is an integral weight of $\gt_n$, where $\htt_n$ is the Lie algebra of $H_n$. An integral weight $\lambda\in\htt^*$ is \emph{regular} if its restrictions $\restr{\lambda}{\htt_n}$ are regular for all $n$.

Let $\lambda$ be an integral weight of $\gt$. Then, using our description of the Picard group of $G/B$, we see that the line bundles $\Ou(\restr{\lambda}{\htt_n})$ form a well-defined projective system. We denote the projective limit by $\Ou(\lambda)$.

We now turn our attention to the cohomology groups $H^q(G/B,\Ou(-\lambda))$. A first natural question is whether the group $H^q(G/B,\Ou(-\lambda))$ is the projective limit of the groups $H^q(G_n/B_n,\Ou(-\restr{\lambda}{\htt_n}))$. This question is answered affirmatively by use of the Mittag-Leffler condition which we now state.

Let $X=\ilm X_n$ be an ind-variety, and $$\ldots\stackrel{\zeta_{n+1}}{\to}\Fo_n\stackrel{\zeta_n}{\to}\Fo_{n-1}\stackrel{\zeta_{n-1}}{\to}\ldots\stackrel{\zeta_2}{\to}\Fo_1\to0$$ be an projective system of sheaves of $\Ou_X$-modules such that the support of $\Fo_n$ is contained in $X_n$. Assume that for some $q\geq0$ the projective system of vector spaces $$\ldots\stackrel{\zeta_{n+1}^q}{\to}H^q(X,\Fo_n)\stackrel{\zeta_n^q}{\to}H^q(X,\Fo_{n-1})\stackrel{\zeta_{n-1}}{\to}\ldots\to0$$ satisfies the condition that for every $n$ the filtration on the vector space $H^q(X,\Fo_n)$ by the subspaces $\zeta_m\circ\ldots\circ\zeta_{n+1}(H^q(X,\Fo_m))$ is eventually constant (Mittag-Leffler condition). Then there is a canonical isomorphism
\begin{equation}
H^q(X,\plm\Fo_n)\cong\plm H^q(X,\Fo_n).\label{formula:cohom_iso_Mittag_Leffler}
\end{equation}
Furthermore, assume that an ind-group $G'=\ilm G_n'$ acts on $X$ (in the category of ind-varieties) so that $G_n'$ acts on $X_n$. If the sheaves $\Fo_n$ are $G_n'$-sheaves and the morphisms $\zeta_n$ are morphisms of $G_n'$-sheaves, then the isomorphism (\ref{formula:cohom_iso_Mittag_Leffler}) is an isomorphism of $G'$-modules. A standard reference for the Mittag-Leffler condition is \cite{Hartshorne1} (see also \cite{DimitrovPenkovWolf1}).

The Mittag-Leffler condition is obviously satisfied in our case since the varieties $X_n=G_n/B_n$ are compact and, consequently, the vector spaces $H^q(G_n/B_n,\Ou(-\restr{\lambda}{\htt_n}))$ are finite dimensional.

A next natural question is for which $\lambda$ and $q$ is the cohomology groups $H^q(G/B,\Ou(-\lambda))$ nonzero. The answer is given by the following theorem.

Let $W=\ilm W_n$ be the group defined in Subsection~\ref{sst:analogues_Weyl_group}. Recall the definition of $\ell_B(w)$ for $w\in W$ from Subsection~\ref{sst:analogues_Weyl_group}. We say that an integral weight $\nu$ is $B$-\emph{dominant} if its restrictions $\restr{\nu}{\htt_n}$ are $B_n$-dominant for all $n$. If $\nu$ is a dominant integral weight, then the finite-dimensional $G_n$-modules $V(\restr{\nu}{H_n})$ form a inductive system, and the corresponding inductive limit $G$-module $\ilm V(\restr{\nu}{\htt_n})$ is denoted by~$V(\nu)$.

\mtheo{The group $H^q(G/B,\Ou(-\lambda))$ is\label{theo:cohom_infin_B_n} nonzero if and only if there exists $w\in W$ of length $\ell_B(w)=q$ such that $\mu=w(\lambda)-\sum_{\alpha\in\Phi^+,~w(\alpha)\notin\Phi^+}\alpha$ is a dominant integral weight. In this case\textup, $$H^q(G/B,\Ou(-\lambda))\cong V(\mu)^*.$$}

Note that in the finite-dimensional case, $$w(\lambda)-\sum_{\alpha\in\Phi^+,~w(\alpha)\notin\Phi^+}\alpha=w(\lambda+\rho)-\rho,$$ so this condition on $\mu$ is completely analogous to the finite-dimensional condition. It is also easy to see that the pair $(w,q)$ is unique whenever it exists, so $\Ou(-\lambda)$ has at most one nonvanishing cohomology group.

Informally, one may say that $\Ou(-\lambda)$ is acyclic unless the weight $\lambda$ is almost dominant, i.e., there is a dominant weight $\mu$ of the form $$\mu=w(\lambda)-\sum_{\alpha\in\Phi^+,~w(\alpha)\notin\Phi^+}\alpha$$ for $w\in W$. In this way, ``most'' sheaves of the form $\Ou(-\lambda)$ are acyclic.

Theorem~\ref{theo:cohom_infin_B_n} follows from more general results proved in \cite{DimitrovPenkovWolf1}. More precisely, Theorem~\ref{theo:cohom_infin_B_n} is a direct corollary of \cite[Proposition 14.1]{DimitrovPenkovWolf1}. Let us note that in \cite{DimitrovPenkovWolf1} much more general $G$-line bundles on $G/P$ are considered. These bundles are induced from (possibly infinite-dimensional) pro-rational $P$-modules which may or may not be dual to highest weight $P$-modules. Therefore, the theory built up in \cite{DimitrovPenkovWolf1} is analogous to considering the bundles of the form $p_*\Ou(-\lambda)$ as in Theorem~\ref{theo:cohom_fin_P_n}. Not that in recent paper \cite{HristovaPenkov1} cohomologies of equivariant finite rank vector bundles on $G/P$ are studied.

In the finite-dimensional case, any line bundle $\Ou(-\lambda)$ on $G_n/P_n$ with $$H^0(G_n/P_n,\Ou(-\lambda))\neq0$$ induces a morphism from $G_n/P_n$ to the projective space $\Pp(V(\lambda))$. This morphism is an embedding if the weight $\lambda$ is regular. A similar statement holds for $G/P$, namely, if $$H^0(G/P,\Ou(-\lambda))\neq0$$ then $\Ou(\lambda)$ induces a morphism $j_{\lambda}$ from $G/P$ to the projective ind-space $\Pp(V(\lambda))$. However, not for all $P$ there exists a line bundle $\Ou(-\lambda)$ such that $\lambda$ is regular and $B$-dominant for a Borel subgroup $B$ contained in $P$. The following theorem is the precise result in this direction, see \cite[Section 15]{DimitrovPenkovWolf1}.

\mtheo{The morphism $j_{\lambda}$ is\label{theo:projectivity_via_BBW} an embedding if and only if $\lambda$ is regular dominant integral for some Borel subgroup $B\subset P$. Such a weight $\lambda$ exists if and only if $P=P_{\Fo}$ where $\Fo$ is an $E$-compatible flag in $V$. The ind-variety $G/P$ is projective if and only if $P=P_{\Fo}$ for $\Fo$ as above.}

As a consequence, most ind-varieties $\Fl(\Fo,E)$ are not projective as the condition for $\Fo$ to be a flag is very restrictive.

\subsection{Vector bundles}\fakesst A\label{sst:vector_bundles} classical result of G. Birkhoff and A. Grothendieck claims that each any (finite-rank) vector bundle on $\Pp^1$ is isomorphic to a direct sum of line bundles $\Ou(n)$, $n\in\Zp$. For $n\geq2$ the classification of vector bundles of finite rank on $\Pp^n$ remains unfinished. On the other hand, a remarkable theorem of Barth--Van de Ven--Tyurin--Sato states that any finite rank bundle on the ind-variety $\Pp^{\infty}$ is isomorphic to a direct sum of line bundles. For rank-two bundles this was established by W. Barth and A. Van de Ven in \cite{BarthVandeVen1}, and for finite rank bundles it was proved by A.N. Tyurin in \cite{Tyurin1} and E. Sato in \cite{Sato1}.

Here we present the results from \cite{PenkovTikhomirov2} establishing sufficient conditions on a locally complete linear ind-variety
$X$ which ensure that the Barth--Van de Ven--Tyurin--Sato Theorem holds on $X$. We then show a class of ind-varieties of generalized flags which satisfy these sufficient conditions.

Let an ind-variety $X$ be the inductive limit of a chain of embeddings
\begin{equation*}
X_1\stackrel{\vfi_1}\hookrightarrow X_2\stackrel{\vfi_2}\hookrightarrow\ldots\stackrel{\vfi_{n-1}}{\to}
X_n\stackrel{\vfi_n}\hookrightarrow X_{n+1}\stackrel{\vfi_{n+1}}\hookrightarrow\ldots
\end{equation*}
of complete algebraic varieties (we call such ind-varieties \emph{locally complete}). By $\Ou_X=\plm\Ou_{X_n}$ we denote the \emph{structure sheaf} of the ind-variety $X$.

\defi{A \emph{vector bundle} $Q$ \emph{of rank} $r\in\Zp_{>0}$ on $X$ is the projective limit $Q=\plm Q_n$ of an projective system of vector bundles $Q_n$ or rank $r$ on $X_n$, i.e., of a system of vector bundles $Q_n$ with fixed isomorphisms $$\psi_n\colon Q_n\to \vfi^*_nQ_{n+1}.$$ (We consider only vector bundles of finite rank.) Here and below $\vfi^*$ stands for the inverse image of vector bundles under a morphism $\vfi$.}

We call the ind-variety $X$ \emph{linear} (cf. Subsection~\ref{sst:linear_ind_grassmannians}) if, for large enough $n$, the induced ho\-mo\-mor\-phisms of Picard groups $$\vfi_n^*\colon\Pic{X_{n+1}}\to\Pic{X_n}$$ are epimorphisms. Any ind-variety of generalized flags is linear. To formulate the main result of this subsection, we need the following three technical conditions.

First, let $X=\ilm X_n$ be a linear ind-variety such that $\Pic{X_n}$ is a free abelian group for all $n$. Assume that there is a finite or countable set $\Theta_X$ and a collection $$\{L_i=\plm L_{i,n}\}_{i\in\Theta_X}$$ of nontrivial line bundles on $X$ such that, for any $n$, $L_{i,n}\cong\Ou_{X_n}$ for all but finitely many indices $i_1(n),...,i_{j(n)}(n)$, and the images of $L_{i_1(n),n},...,L_{i_{j(n)}(n),n}$ in $\Pic{X_n}$ form a basis of $\Pic{X_n}$. In this case $\Pic{X}$ is isomorphic to a direct product of infinite cyclic
groups which generators are the images of $L_i$. Denote by $\bigotimes_{i\in\Theta_X}L_i^{\otimes a_i}$ the line bundle on $X$ whose restriction to $X_n$ equals $$\bigotimes\nolimits_{i\in\Theta_X}L_{i,n}^{\otimes a_i}=\bigotimes\nolimits_{k=1}^{j(n)}L_{i_k(n),n}^{\otimes a_{i_k(n)}}.$$
We say that $X$ \emph{satisfies the property} L if, in addition to the above condition, $$H^1(X_n,\bigotimes\nolimits_{i\in\Theta_X}L_{i,n}^{\otimes a_i})=0$$ for any $n\ge1$ if some $a_i$ is negative.

Assume that $X$ satisfies the property L. For a given $i\in\Theta_X$, a smooth rational curve $C\cong\Pp^1$ on $X$ is a \emph{projective line of the $i$-th family on} $X$ (or, simply, a \emph{line of the $i$-th family}), if
\begin{equation*}
\restr{L_j}{C}\cong\Ou_{\Pp^1}(\delta_{i,j})\text{ for all }j\in\Theta_X,
\end{equation*}
where $\delta_{i,j}$ is the usual Kronecker delta. By $B_i$ we denote the set of all projective lines of the $i$-th family on $X$. It has a natural structure of an ind-variety: $B_i=\ilm B_{i,n}$, where
$$B_{i,n}=\{C\in B_i\mid C\subset X_n\}.$$ For any point $x\in X$ the subset $$B_i(x)=\{C\in B_i\mid C\ni x\}$$ inherits an induced structure of an ind-variety.

Assume also that for any $i\in\Theta_X$ there exists an ind-variety $\Pi_i$ which is a inductive limit of a chain of embeddings $$Pi_{i,1}\stackrel{\pi_{i,1}}\hookrightarrow\Pi_{i,2}\stackrel{\pi_{i,2}}\hookrightarrow\ldots\stackrel{\pi_{i,n-1}}{\hookrightarrow}
\Pi_{i,n}\stackrel{\pi_{i,n}}\hookrightarrow\Pi_{i,n+1}\stackrel{\pi_{i,n+1}}\hookrightarrow\ldots,$$ where the points of
$\Pi_{i,n}$ are projective subspaces $\Pp^{m_n}$ of $B_{i,n}$, together with linear embeddings $$\Pp^{m_n}\hookrightarrow\Pp^{m_{n+1}}={\pi_{i,n}}(\Pp^{m_n})$$ induced by the embeddings $B_{i,n}\hookrightarrow B_{i,n+1}$, so that each point of $\Pi_i$ is considered as a projective ind-subspace
$\Pp^\infty=\ilm\Pp^{m_n}$ of $B_i$. Given $x\in X$, consider the following conditions:
\begin{equation*}
\begin{split}
\text{(A.i) }&\text{for each $n\geq1$ such that $x\in X_n$, each nontrivial sheaf $L_{i,n}$ defines a morphism}\\
\text{$\hphantom{\text{(A.i) }}$}&\text{$\psi_{i,n}\colon X_n\to\Pp^{r_{i,n}}=\Pp(H^0(X_n,L_{i,n})^*)$ which maps the family of lines $B_{i,n}(x)$}\\
\text{$\hphantom{\text{(A.i) }}$}&\text{isomorphically to a subfamily of lines in $\mathbb{P}^{r_{i,n}}$ passing through the point $\psi_{i,n}(x)$;}\\
\text{(A.ii) }&\text{the variety $\Pi_{i,n}(x)=\{\Pp^{m_n}\in\Pi_{i,n}\mid\Pp^{m_n}\subset B_{i,n}(x)\}$ is connected for any $n\geq1$;}\\
\text{(A.iii) }&\text{the projective ind-subspaces $\Pp^\infty\in\Pi_i(x)=\ilm\Pi_{i,n}(x)$ fill $B_i(x)$;}\\
\text{(A.iv) }&\text{for any $d\in\mathbb{Z}_{>0}$ there exists $n_0(d)\in\mathbb{Z}_{>0}$ such that, for any $d$-dimensional
variety $Y$}\\
\text{$\hphantom{\text{(A.iv) }}$}&\text{and any $n\geq n_0(d)$, any morphism $\Pi_{i,n}(x)\to Y$ is a constant map.}
\end{split}
\end{equation*}
In particular, (A.ii) and (A.iii) imply that the varieties $\Pi_{im},$ $B_{im},$ $B_{im}(x)$ are connected. If all these conditions are satisfied for all $x\in X$, we say that $X$ \emph{satisfies the property} A.

Finally, suppose that $X$ satisfy the properties L and A as above. A vector bundle $Q$ on $X$ is called $B_i$-\emph{uniform}, if for any projective line $\Pp^1\in B_i$ on $X$, the restricted bundle $\restr{Q}{\Pp^1}$ is isomorphic to $\bigoplus_{j=1}^{\rk Q}\Ou_{\Pp^1}(k_j)$
for some integers $k_j$ not depending on the choice of $\Pp^1$. If in addition all $k_j=0$, then $Q$ is called $B_i$-\emph{linearly trivial}. We call $Q$ \emph{uniform} (respectively, \emph{linearly trivial}) if $Q$ is $B_i$-uniform (respectively, $B_i$-linearly trivial) for any
$i\in\Theta_X$. We say that $X$ \emph{satisfies the property} T if any linearly trivial vector bundle on $X$ is trivial.

The main result of this subsection is as follows.

\mtheo{\textup{\cite[Theorem 1]{PenkovTikhomirov2}} Let $Q$ be a vector bundle on a linear ind-variety $X$. \textup{i) }If $X$ satisfies the properties \textup{L} and \textup{A} some fixed line bundles $\{L_i\}_{i\in\Theta_X}$ and corresponding families $\{B_i\}_{i\in\Theta_X}$ of projective lines on $X$\textup, then $Q$ has a filtration by vector subbundles $$0=Q_0\subset Q_1\subset\ldots\subset Q_t=Q$$ with uniform quotient bundles $Q_k/Q_{k-1}$\textup, $1\leq k\leq t$. \textup{ii)} If\textup, in addition\textup, $X$ satisfies the property \textup{T} then the above filtration of $Q$ splits and its quotients are of the form $$Q_k/Q_{k-1}\cong\rk(Q_k/Q_{k-1})\bigotimes\nolimits_{i\in\Theta_X}L_i^{\otimes a_{i,k}}$$ for some $a_{i,k}\in\Zp$\textup, $1\leq k\leq t$. In particular\textup, $Q$ is isomorphic to a direct sum of line bundles.}

If $X$ is one of linear ind-grassmannians $\Gra{k}$ or $\Grb{k}{\infty}$ considered in Section~\ref{sect:ind_grassmannians}, then there exists the tautological bundle $S$ on $X$ with $\rk S=k$. For $k\geq2$, this bundle is not isomorphic to a direct sum of line bundles, hence the Barth--Van de Ven--Tyurin--Sato theorem does not hold  for these ind-grassmannians. On the other hand, the following theorem holds.

\corop{\textup{i)} Suppose that $X=\Gra{\infty}$\textup, $\Grb{\infty}{k}$\textup, $\Grbi{\infty}{k}{0}$ or $\Grbi{\infty}{k}{1}$. Then any vector bundle on $X$ is isomorphic to $\bigoplus_i\Ou_X(k_i)$ for some $k_i\in\Zp$. \textup{ii)} Let $X=\Fl(\Fo,E)$\textup, where $\Fo$ is a flag in~$V$ such that $\codim_{F''}F'=\infty$ for all $(F',F'')\in\Fo^{\dag}$. Then any vector bundle on $X$ is isomorphic to a direct sum of line bundles.}{i) It follows from our description of $\Pic{X}$ that $\Pic{X}\cong\Zp$ and that any line bundle on $X$ is isomorphic to $\Ou_X(m)$, where by definition $\restr{\Ou_X(m)}{X_n}\cong\Ou_{X_n}(m)$. Since $\Ou_X(1)$ is very ample, condition (A.i) holds. It is checked in \cite[Section 4]{PenkovTikhomirov2} that conditions (A.ii)--(A.iv) also hold, and so $X$ satisfies the property A. The property T for $X$ is verified in \cite[Section 4]{PenkovTikhomirov2} as well.

ii) It is shown in \cite[Subsection 6.3]{PenkovTikhomirov2} that $X$ satisfies the properties L, A, T.}

A characterization of bundles on an arbitrary ind-variety of generalized flags remains unknown. Nevertheless, the second author has made the following conjecture: if $\Fo$ is a maximal generalized flag, then, for each finite-rank vector bundle $Q$ on $\Fl(\Fo,E)$, $Q$ admits a filtration by subbundles so that the successive quotients are line bundles.

Note also that in \cite{PenkovTikhomirov3} vector bundles on so-called twisted ind-grassmannians are studied. Twisted ind-grassmannians are not ind-varieties of generalized flags: they are inductive limits $\ilm X_n$ of grass\-mannians for nonlinear embeddings $$X_n\hookrightarrow X_{n+1}.$$ It is proved in \cite{PenkovTikhomirov3} that any finite-rank vector bundle on a twisted ind-grassmannian is trivial.

\section{Orbits of real forms}\label{sect:orbits_real_forms}\fakesect

In this section we fix a real form $G^0$ of the group $G=\SL_{\infty}(\Cp)$ and study the structure of $G^0$-orbits on the ind-variety $\Fl=\Fl(\Fo,E)$. Our exposition is based on \cite{IgnatyevPenkovWolf1}.

In the finite-dimensional setting, the study of orbits of real form of semisimple complex Lie groups on their flag varieties has its roots in linear algebra. Witt's Theorem claims that, given a finite-dimensional vector space $W$ endowed with a nondegenerate bilinear for (symmetric or skew-symmetric) or a nondegenerate Hermitian form, two subspaces $W_1$, $W_2$ of $W$ are isometric within $W$ (i.e., one is obtained from the other via an isometry of $W$) if and only if $W_1$ and $W_2$ are isometric.

When $W$ is a Hermitian space, this is a statement about the orbits of the unitary group $U(W)$ on the complex grassmannian $\Grr{k}{W}$, where $$k=\dim W_1=\dim W_2.$$ More precisely, the orbits of $U(W)$ on $\Grr{k}{W}$ are parameterized by the possible signatures of a, possibly degenerate, Hermitian form on a $k$-dimensional space of $W$.

A general theory of orbits of a real form $G_n^0$ of a (finite-dimensional) semisimple complex Lie group $G_n$ on a flag variety $G_n/P_n$ was developed by J.A. Wolf in \cite{Wolf1} and \cite{Wolf2}. This theory has become a standard tool in semisimple representation
theory and complex algebraic geometry. An important further development has been the theory of cycle spaces initiated by A. Huckleberry and J.A. Wolf, see the monograph \cite{FelsHuckleberryWolf1}.

Our main result in this section is the fact that any $G^0$-orbit in $G/P$, when intersected with a finite-dimensional flag variety $G_n/P_n$ from a given exhaustion of $G/P$, yields a single $G_n^0$-orbit. This means that the mapping
\begin{equation*}
\text{$\{G_n^0$-orbits on $G_n/P_n\}$ $\to$ $\{G_{n+1}^0$-orbits on $G_{n+1}/P_{n+1}\}$}
\end{equation*}
is injective. Using this feature, we are able to answer the following questions.
\begin{enumerate}
\item When are there finitely many $G^0$-orbits on $G/P$?
\item When is a given $G^0$-orbit on $G/P$ closed?
\item When is a given $G^0$-orbit on $G/P$ open?
\end{enumerate}
The answers depend on the type of real form and not only on the parabolic subgroup $P\subset G$. For instance, if $P=B$ is the upper-triangular Borel ind-subgroup of $\SL_{\infty}(\Cp)$ with positive roots $$\{\epsi_i-\epsi_j\},~i,~j\in\Zp_{>0},~i<j,$$ then $G/B$ has no closed $\SU(\infty,\infty)$-orbit and no open $\SL(\infty,\Rp)$-orbit.

\subsection{Finite-dimensional case}\fakesst\label{sst:finite_dim_case} Let $W$ be a finite-dimensional complex vector space. Recall that a~\emph{real structure} on $W$ is an antilinear involution $\tau$ on $W$. The set $$W^0=\{v\in W\mid\tau(v)=v\}$$ is a \emph{real form} of $W$, i.e. $W^0$ is a real vector subspace of $W$ such that $$\dim_{\Rp}W^0=\dim_{\Cp}W$$ and the $\Cp$-linear span $\langle W^0\rangle_{\Cp}$ coincides with $W$. A real form $W^0$ of $W$ defines a unique real structure $\tau$ on $W$ such that $W^0$ is the set of fixed points of~$\tau$. A~\emph{real form} of a complex finite-dimensional Lie algebra $\gt$ is a real Lie subalgebra $\gt^0$ of $\gt$ such that $\gt^0$ is a~real form of $\gt$ as a complex vector space.

Let $G_n$ be a complex semisimple connected algebraic group and $G_n^0$ be a \emph{real form} of $G_n$, i.e. $G_n^0$ is a real closed algebraic subgroup of $G_n$ such that its Lie algebra $\gt_n^0$ is a real form of the Lie algebra $\gt_n$ of $G_n$. Let $P_n$ be a parabolic subgroup of $G_n$ and $G_n/P_n$ be the corresponding flag variety. The group $G_n^0$ naturally acts on $G_n/P_n$. In \cite{Wolf1} the following facts about the $G_n^0$-orbit structure of $G_n/P_n$ are proved, see \cite[Theorems 2.6, 3.3, 3.6, Corollary 3.4]{Wolf1} (here we use the usual differentiable manifold topology on $G_n/P_n$).

\mtheo{\label{theo:finite_dim_case}Let $G_n$, $P_n$, $G_n^0$ be as above.
\begin{equation*}
\begin{split}
&\text{\textup{i)} Each $G_n^0$-orbit is a real submanifold of $G_n/P_n$}.\\
&\text{\textup{ii)} The number of $G_n^0$-orbits on $G_n/P_n$ is finite}.\\
&\text{\textup{iii)} The union of the open $G_n^0$-orbits is dense in $G_n/P_n$}.\\
&\text{\textup{iv)} There is a unique closed orbit $\Omega$ on $G_n/P_n$}.\\
&\text{\textup{v)} The inequality $\dim_{\Rp}\Omega\geq\dim_{\Cp}G_n/P_n$ holds}.
\end{split}
\end{equation*}}

Here is how this theorem relates to Witt's Theorem in the case of a Hermitian form. Let $W$ be an $n$-dimensional complex vector space and $G_n=\SL(W)$. Fix a nondegenerate Hermitian form $\omega$ of signature $(p,n-p)$ on the vector space $W$ and denote by $G_n^0=\SU(W,\omega)$ the group of all linear operators on $W$ of determinant~$1$ which preserve $\omega$. Then $G_n^0$ is a real form of $G_n$. Given $k\leq n$, the group~$G_n$ naturally acts on the grassmannian $\Grr{k}{W}$ of all $k$-dimensional complex subspaces of $W$. To each $U\in\Grr{k}{W}$ one can assign its \emph{signature} $s(U)=(a,b,c)$, where the restricted form $\restr{\omega}{U}$ has rank $a+b$ with $a$ positive squares and $b$ negative ones, and $c$ equals the nullity of $\restr{\omega}{U}$. By Witt's Theorem, two subspaces $U_1$, $U_2\in\Grr{k}{W}$ belong to the same $G_n^0$-orbit if and only if their signatures coincide.

Moreover, one can verify the following formula for the number $|\Grr{k}{W}/G_n^0|$ of $G_n^0$-orbits on $\Grr{k}{W}$. Set $l=\min\{p,~n-p\}$. Then
\begin{equation*}
|\Grr{k}{W}/G_n^0|=\begin{cases}
(-k^2-2l^2-n^2+2kn+2ln+k+n+2)/2,&\text{if }n-l\leq k,\\
(l+1)(2k-l+2)/2,&\text{if }l\leq k\leq n-l,\\
(k+1)(k+2)/2,&\text{if }k\leq l.\\
\end{cases}
\end{equation*}
A $G_n^0$-orbit of a subspace $U\in\Grr{k}{W}$ is open if and only if the restriction of $\omega$ to $U$ is nondegenerate, i.e., if $c=0$. Therefore, the number of open orbits equals $\min\{k+1,l+1\}$. There is a unique closed\break $G_n^0$-orbit $\Omega$ on $\Grr{k}{W}$, and it consists of all $k$-dimensional subspaces of $W$ such that $c=\min\{k,l\}$ (the condition $c=\min\{k,l\}$ maximizes the nullity of the form $\restr{\omega}{U}$ for $k$-dimensional subspaces\break $U\subset W$). In particular, if $k=p\leq n-p$ then $\Omega$ consists of all isotropic $k$-dimensional complex subspaces of $W$. See \cite{Wolf1} for more details in this latter case.

\subsection{Orbits of real forms as ind-manifolds}\fakesst Below\label{sst:orbits_ind_manifolds} we recall the~classification of real forms of $G=\SL_{\infty}(\Cp)$ due to A. Baranov \cite{Baranov1}. By definition, a real ind-subgroup $G^0$ of $G$ is called a \emph{real form} of~$G$ if $G$ can be represented as a nested union $G=\bigcup\Go_n$ of its finite-dimensional Zariski closed subgroups such that $\Go_n$ is a semisimple algebraic group and $G^0\cap\Go_n$ is a real form of~$\Go_n$ for each $n$. To define real forms of $G$, pick a basis $\Eo=\{\epsi_1,~\epsi_2,~\ldots\}$ of $V$ and its exhaustion $$\Eo=\bigcup_{n\geq1}\Eo_n$$ by finite subsets such that $V_n=\langle\Eo_n\rangle_{\Cp}$ and $$\langle\Eo_{n+1}\setminus\Eo_n\rangle_{\Cp}=\langle E_{n+1}\setminus E_n\rangle_{\Cp}$$ for all $n\geq1$. Recall that the embedding $G_n\hookrightarrow G_{n+1}$ is given by the formula $\vfi\mapsto\wh\vfi$, where $\wh\vfi(x)=\vfi(x)$ for $x\in V_n$ and $\wh\vfi(e)=e$ for $e\in E\setminus E_n$.

Fix a real structure $\tau$ on $V$ such that $\tau(\epsi)=\epsi$ for all $\epsi\in\Eo$. Then each finite-dimensional space $V_n$ is $\tau$-invariant. Denote by $\GL(V_n,\Rp)$ (respectively, by~$\SL(V_n,\Rp)$) the group of invertible (respectively, of determinant 1) operators on $V_n$ defined over $\Rp$. Recall that a linear operator on a complex vector space with a real structure is \emph{defined over} $\Rp$ if it commutes with the real structure or, equivalently, if it maps the real form to itself. For each $n$, the map $\vfi\mapsto\wh\vfi$ gives an embedding $$\SL(V_n,\Rp)\hookrightarrow\SL(V_{n+1},\Rp),$$ so the inductive limit $$G^0=\ilm\SL(V_n,\Rp)$$ is well defined. We denote this real form of $G$ by $\SL(\infty,\Rp)$.

Fix a nondegenerate Hermitian form $\omega$ on $V$. Suppose that its restriction $\omega_n=\restr{\omega}{V_n}$ is non\-de\-ge\-ne\-rate for all $n$, and that $\omega(\epsi,V_n)=0$ for $\epsi\in\Eo\setminus\Eo_n$. Denote by $p_n$ the dimension of a maximal $\omega_n$-positive definite subspace of $V_n$, and put $q_n=\dim V_n-p_n$. Let $\SU(p_n,q_n)$ be the subgroup of $G_n$ consisting of all operators preserving the form $\omega_n$. For each $n$, the map $\vfi\mapsto\wh\vfi$ induces an embedding $$\SU(p_n,q_n)\hookrightarrow\SU(p_{n+1},q_{n+1}),$$ so we have a inductive limit $$G^0=\ilm\SU(p_n,q_n).$$ If there exists $p$ such that $p_n=p$ for all sufficiently large $n$ (respectively, if $\lim_{n\to\infty}p_n=\lim_{n\to\infty}q_n=\infty$), then we denote this real form of $G$ by $\SU(p,\infty)$ (respectively, by $\SU(\infty,\infty)$).

Finally, fix a \emph{quaternionic structure} $J$ on $V$, i.e. an antilinear automorphism of $V$ such that $J^2=-\id_V$. Assume that the complex dimension of $V_n$ is even for $n\geq1$, and that the restriction $J_n$ of $J$ to $V_n$ is a quaternionic structure on $V_n$. Furthermore, suppose that $$J(\epsi_{2i-1})=-\epsi_{2i},~J(\epsi_{2i})=\epsi_{2i-1}$$ for $i\geq1$. Let $\SL(V_n,\Hp)$ be the subgroup of $G_n$ consisting of all linear operators commuting with $J_n$. Then, for each $n$, the map $\vfi\mapsto\wh\vfi$ induces an embedding of the groups $$\SL(V_n,\Hp)\hookrightarrow\SL(V_{n+1},\Hp),$$ and we denote the corresponding inductive limit by $$G^0=\SL(\infty,\Hp)=\ilm\SL(V_n,\Hp).$$ The ind-group $\SL(\infty,\Hp)$ is also a real form of $G$.

The following result is a corollary of \cite[Theorem 1.4]{Baranov1} and \cite[Corollary 3.2]{DimitrovPenkov2}.

\mtheo{If $G=\SL_{\infty}(\Cp)$ then\textup, up to isomorphism\textup, $\SL(\infty,\Rp)$\textup, $\SU(p,\infty)$\textup, $0\leq p<\infty$\textup, $\SU(\infty,\infty)$\textup, $\SL(\infty,\Hp)$ are all distinct real\label{theo:list_of_real_forms} forms of $G$.}

Now, let $\Fo$ be a generalized flag compatible with the basis $E$, $\Fl=\Fl(\Fo,E)$, and $\Fl_n=\Fl(d_n,V_n)$, where $d_n$ is the type of the flag $\Fo\cap V_n$. Then $\Fl=\ilm\Fl_n$, where the embedding $$\iota_n\colon\Fl_n\hookrightarrow\Fl_{n+1}$$ is given by formula (\ref{formula:iota_n}) (or, equivalently, $\iota_n$ is the composition of embeddings given by formula (\ref{formula:iota_n_another})), see Subsection~\ref{sst:gen_flags}.

Let $G^0$ be a real form of~$G$ (see Theorem \ref{theo:list_of_real_forms}). The group $G_n=\SL(V_n)$ naturally acts on $\Fl_n$, and the map $\iota_n$ is equivariant: $$g\cdot\iota_n(x)=\iota_n(g\cdot x),~g\in G_n\subset G_{n+1},~x\in\Fl_n.$$ Put also $G_n^0=G^0\cap G_n$. Then $G_n^0$ is a real form of~$G_n$. For the rest of the section we impose some further assumptions on $V_n$, which we now describe case-by-case.

Let $G^0=\SU(p,\infty)$ or $\SU(\infty,\infty)$. Recall that the restriction $\omega_n$ of the fixed nondegenerate Hermitian form $\omega$ to $V_n$ is nondegenerate. From now on we assume that if $e\in E_{n+1}\setminus E_n$ then $e$ is orthogonal to~$V_n$ with respect to $\omega_{n+1}$. Next, let $G^0=\SL(\infty,\Rp)$. Here we assume that $m_n$ is odd for each $n\geq1$, and that $\langle E_n\rangle_{\Rp}$ is a real form of $V_n$. Finally, for $G^0=\SL(\infty,\Hp)$, we assume that $m_n$ is even for all $n\geq1$ and that $$J(e_{2i-1})=-e_{2i},~J(e_{2i})=e_{2i-1}$$ for all $i$. These additional assumptions align the real form $G^0$ with the ind-variety $\Fl$.

Our main result in this section is as follows.

\theop{If \label{theo:Omega_n_cap}
$\iota_n(\Fl_n)$ has nonempty intersection with a $G_{n+1}^0$-orbit
then that intersection is a single $G_n^0$-orbit.}
{The proof goes case by case. Here we consider the case $G^0=\SU(\infty,\infty)$. The proof for $G^0=\SU(p,\infty)$, $0\leq p<\infty$, is completely similar, while the proof for $G^0=\SL(\infty,\Rp)$ and $\SL(\infty,\Hp)$ is based on other ideas, see \cite[Theorem 3.1]{IgnatyevPenkovWolf1}.

Pick two flags
\begin{equation*}
\begin{split}
\Du&=\{\{0\}=D_0\subset D_1\subset\ldots\subset D_s=V_n\},\\
\Bu&=\{\{0\}=B_0\subset B_1\subset\ldots\subset B_s=V_n\}
\end{split}
\end{equation*}
in $\Fl_n$ such that $\wt\Du=\iota_n(\Du)$ and $\wt\Bu=\iota_n(\Bu)$ belong to a given $G_{n+1}^0$-orbit. Let
\begin{equation*}
\begin{split}
\wt\Du&=\{\{0\}=\wt D_0\subset \wt D_1\subset\ldots\subset\wt D_{\wt s}=V_{n+1}\},\\
\wt\Bu&=\{\{0\}=\wt B_0\subset \wt B_1\subset\ldots\subset\wt B_{\wt s}=V_{n+1}\}.
\end{split}
\end{equation*}
There exists $\wt\vfi\in\SU(\omega_{n+1},V_{n+1})$ satisfying $\wt\vfi(\wt\Du)=\wt\Bu$, i.e. $\wt\vfi(\wt D_i)=\wt B_i$ for $i=0,~\ldots,~\wt s$. To prove the result, we must construct an isometry $\vfi\colon V_n\to V_n$ satisfying $\vfi(\Du)=\Bu$. Of course, one can then scale $\vfi$ to obtain an isometry of determinant $1$. By \cite[Theo\-rem~6.2]{Huang1}, an isometry $\vfi\colon V_n\to V_n$ with  $\vfi(\Du)=\Bu$ exists if and only if $D_i$ and $B_i$ are isometric for all $i$ from $1$ to $s$, and
\begin{equation}
\dim(D_i\cap D_j^{\perp,V_n})=\dim(B_i\cap B_j^{\perp,V_n})\label{formula:Witt_for_flags}
\end{equation}
for all $i<j$ from $1$ to $s$. (Here $U^{\perp,V_n}$ denotes the $\omega_n$-orthogonal complement within $V_n$ of a subspace $U\subset V_n$.) Pick $i$ from $1$ to $s$. Since $e_{n+1}$ is orthogonal to $V_n$ and $\wt\vfi$ establishes an isometry between $\wt D_i$ and $\wt B_i$, the first condition is satisfied. So it remains to prove (\ref{formula:Witt_for_flags}).

To do this, denote $$C_n=\langle E_{n+1}\setminus E_n\rangle_{\Cp}.$$ Since $C_n$ is orthogonal to $V_n$, for given subspaces $U\subset V_n$, $W\subset C_n$ one has $$(U\oplus W)^{\perp,V_{n+1}}=U^{\perp,V_n}\oplus W^{\perp,C_n}.$$ Hence, if $$\wt D_k=D_k\oplus W_k,~\wt B_k=B_k\oplus W_k$$ for $k\in\{i,j\}$ and some subspaces $W_i$, $W_j\subset C_n$, then $$\wt D_i\cap\wt D_j^{\perp,V_{n+1}}=(D_i\oplus W_i)\cap(D_j^{\perp,V_n}\oplus W_j^{\perp,C_n})=(D_i\cap D_j^{\perp,V_n})\oplus(W_i\cap W_j^{\perp,C_n}),$$ and the similar equality holds for $\wt B_i\cap\wt B_j^{\perp,V_{n+1}}$. The result follows.}

The following result is an immediate corollary of this theorem. (The definition of a real ind-manifold in the $C^{\infty}$-category is completely similar to the definition of an ind-variety.)

\corop{Let $\Omega$ be a $G^0$-orbit on $\Fl$\textup, and $\Omega_n=\iota_n^{-1}(\Omega)\subset\Fl_n$. Then
\begin{equation*}
\begin{split}
&\text{\textup{i)} $\Omega_n$ is a single $G_n^0$-orbit};\\
&\text{\textup{ii)} $\Omega$ is an infinite-dimensional real ind-manifold}.\\
\end{split}
\end{equation*}}{i) Suppose $\Du$, $\Bu\in\Omega_n$. Then there exists $m\geq n$ such that images of $\Du$ and $\Bu$ under the morphism $\iota_{m-1}\circ\iota_{m-2}\circ\ldots\circ\iota_n$ belong to the same $G_m^0$-orbit. Applying Theorem~\ref{theo:Omega_n_cap} subsequently to $\iota_{m-1}$, $\iota_{m-2}$, $\ldots$, $\iota_n$, we see that $\Du$ and~$\Bu$ belong to the same $G_n^0$-orbit.

ii) By definition, $\Omega=\ilm\Omega_n$. Next, (i) implies that the orbit $\Omega$ is a real ind-manifold. By Theorem~\ref{theo:finite_dim_case}~(v), we have $\dim_{\Rp}\Omega_n\geq\dim_{\Cp}\Fl_n$. Since $$\lim_{n\to\infty}\dim_{\Cp}\Fl_n=\infty,$$ we conclude that the orbit $\Omega$ is infinite dimensional.}

We can now give a criterion for $\Fl$ to have finitely many of $G^0$-orbits. We define a generalized flag $\Go$ to be \emph{finite} if it consists of finitely many (possibly infinite-dimensional) subspaces. We say that a generalized flag $\Go$ has \emph{finite type} if it consists of finitely many subspaces of $V$ each of which has either finite dimension or finite codimension in $V$. A finite type generalized flag is clearly a flag. An ind-variety $\Fl(\Go,E)$ is \emph{of finite type} if $\Go$ (equivalently, any $\wt\Go\in\Fl(\Go,E)$) is of finite type.

\propp{For $G^0=\SU(\infty,\infty)$, $\SL(\infty,\Rp)$ and $\SL(\infty,\Hp)$, there are finitely many $G^0$-orbits\label{prop:finiteness} on $\Fl$ if and only if $\Fl$ is of finite type. For $G^0=\SU(p,\infty)$, $0<p<\infty$, there are finitely many $G^0$-orbits on $\Fl$ if and only if $\Fo$ is finite. For $G^0=\SU(0,\infty)$, the ind-variety $\Fl$ is a single $G^0$-orbit.}{Consider the case $G^0=\SU(p,\infty)$, $0<p<\infty$. First suppose that $\Fo$ is finite, i.e.,that $|\Fo|=N<\infty$. Given $n\geq1$, denote
\begin{equation*}
\begin{split}
S_n&=\{s(A)\mid A\subset V_n\},\\
P_n&=\{\dim A\cap B^{\perp,V_n}\mid A\subset B\subset V_n\}.
\end{split}
\end{equation*}
Let $s(A)=(a,b,c)$ for some subspace $A$ of $V_n$. Then, clearly, $a\leq p$ and $c\leq p$, hence $|S_n|\leq p^2$. On the other hand, if $A\subset B$ are subspaces of $V_n$ then $A^{\perp,V_n}\supset B^{\perp,V_n}$, so $$A\cap B^{\perp,V_n}\subset A\cap A^{\perp,V_n}.$$ But $$\dim A\cap A^{\perp,V_n}=c\leq p.$$ Thus $|P_n|\leq p$. Now \cite[Theo\-rem~6.2]{Huang1} shows that the number of $G_n^0$-orbits on $\Fl_n$ is less or equal to $${N\cdot|S_n|\cdot N^2\cdot|P_n|\leq N^3p^3}.$$ Hence, by Theorem~\ref{theo:Omega_n_cap}, there are finitely many $G^0$-orbits on the ind-variety $\Fl$.

Suppose next that $\Fo$ is infinite. In this case, given $m\geq1$, there exists $n$ such that the length of each flag from $\Fl_n$ is at least~$m$, the positive index of $\restr{\omega}{V_n}$ (i.e. the dimension of a maximal positive definite subspace of $V_n$) equals $p$, and $\codim_{V_n}F_m\geq p$, where $$\Fo_n=\{F_1\subset\ldots\subset F_m\subset\ldots\subset V_n\}.$$ It is easy to check that the number of $G_n^0$-orbits on $\Fl_n$ is at least~$m$. Then, by Theorem~\ref{theo:Omega_n_cap}, there are at least $m$ $G^0$-orbits on $\Fl$ or, in other words, there are infinitely many $G^0$-orbits on $\Fl$. The proof for $\SU(p,\infty)$, $p>0$ is complete. For the proof of other cases see \cite[Proposition 4.1]{IgnatyevPenkovWolf1}.}

\subsection{Open and closed orbits}\fakesst In\label{sst:open_closed_orbits} this subsection we provide necessary and sufficient conditions for a given $G^0$-orbit on $\Fl=\Fl(\Fo,E)$ to be open or closed. It turns out that, for all real forms except $\SU(p,\infty)$, $\Fl$ has both an open and a closed orbit if and only if the number of orbits is finite.

First, consider the case of open orbits. Pick any $n$. Recall \cite{HuckleberryWolf1} that the $G_n^0$-obit of a flag $$\Fo_n=\{A_1\subset A_k\subset\ldots\subset A_k\}\in\Fl_n$$ is open if and only if
\begin{equation*}
\begin{split}
&\text{for\label{formula:nondegen_fin_dim} $G^0=\SU(p,\infty)$ or $\SU(\infty,\infty)$: all $A_i$'s are nondegenerate with respect to $\omega$;}\\
&\text{for $G^0=\SL(\infty,\Rp)$: for all $i$, $j$, $\dim A_i\cap\tau(A_j)$ is minimal},\\
&\hphantom{\text{for $G^0=\SL(\infty,\Rp)$: }}\text{i.e., $\dim A_i\cap\tau(A_j)=\max\{\dim A_i+\dim A_j-\dim V_n,~0\}$;}\\
&\text{for $G^0=\SL(\infty,\Hp)$: for all $i$, $j$, $\dim A_i\cap J(A_j)$ is minimal in the above sense}.\\
\end{split}
\end{equation*}

Note that, for any two generalized flags $\Fo_1$ and $\Fo_2$ in $\Fl$, there is a canonical identification of $\Fo_1$ and $\Fo_2$ as linearly ordered sets. For a space $F\in\Fo_1$, we call the image of $F$ under this identification \emph{the space in $\Fo_2$ corresponding to} $F$.

Fix an antilinear operator $\mu$ on $V$. A point $\Go$ in $\Fl$ is \emph{in general position with respect to} $\mu$ if $F\cap\mu(H)$ does not properly contain $\wt F\cap\mu(\wt H)$ for all $F,~H\in\Go$ and all $\wt\Go\in\Fl$, where $\wt F$, $\wt H$ are the spaces in $\wt\Go$ corresponding to $F$, $H$ respectively. A similar definition can be given for flags in~$\Fl_n$. Note that, for $G^0=\SL(\infty,\Rp)$ or $\SL(\infty,\Hp)$, the $G_n^0$-orbit of $\Fo_n\in\Fl_n$ is open if and only if $\Fo_n$ is in general position with respect to $\tau$ or~$J$ respectively.

With the finite-dimensional case in mind, we give the following
\defi{A generalized\label{defi:gen_position} flag $\Go$ is \emph{nondegenerate} if
\begin{equation*}
\begin{split}
\text{for }G^0&=\SU(p,\infty)\text{ or $\SU(\infty,\infty)$:}\\
&\text{each $F\in\Go$ is nondegenerate with respect to $\omega$;}\\
\text{for }G^0&=\SL(\infty,\Rp)\text{ or $\SL(\infty,\Hp)$}:\\
&\text{$\Go$ is in general position with respect to $\tau$ or $J$ respectively.}
\end{split}
\end{equation*} A nondegenerate generalized flag can be thought of as being ``in general position with respect to $\omega$''. Therefore, all conditions in Definition~\ref{defi:gen_position} are clearly analogous.}

For a generalized flag $\Go\in\Fl$, let $n_{\Go}$ be a fixed positive integer such that $\Go$ is compatible with a basis containing $E\setminus E_{n_{\Go}-1}$ (here we put $E_0=\varnothing$; note also that the definition of $n_{\Go}$ here slightly differs from the one in Subsection~\ref{sst:gen_flags}).

\propp{The $G^0$-orbit $\Omega$ of $\Go\in\Fl$ is open if and only if $\Go$ is nondegenerate.}{By the definition of the topology on $\Fl$, the orbit $\Omega$ is open if and only if $\Omega_n=\iota_n^{-1}(\Omega\cap\iota_n(\Fl_n))$ is open for each $n$. Consider the case $G^0=\SU(p,\infty)$ or $\SU(\infty,\infty)$ (for the proof of other cases see \cite[Proposition 5.3]{IgnatyevPenkovWolf1}). To prove the claim in this case, it suffices to show that $A\in\Go$ is nondegenerate with respect to $\omega$ if and only if $\restr{\omega}{A\cap V_n}$ is nondegenerate for all $n\geq n_{\Go}$. This is straightforward. Indeed, if $A$ is degenerate, then there exists $v\in A$ such that $\omega(v,w)=0$ for all $w\in A$. Let $v\in V_{n_0}$ for some $n_0\geq n_{\Go}$. Then $\restr{\omega}{A\cap V_{n_0}}$ is degenerate. On the other hand, if $v\in A\cap V_n$ is orthogonal to all $w\in A\cap V_n$ for some $n\geq n_{\Go}$, then $v$ is orthogonal to all $w\in A$ because $e$ is orthogonal to $V_n$ for $e\in E\setminus E_n$. The result follows.}

We say that two generalized flags \emph{have the same type} if there is an automorphism of $V$ transforming one into the other. Clearly, two $E$-commensurable generalized flags have the same type. On the other hand, it is easy to see that two generalized flags $\Fo$ and $\wt\Fo$ of the same type do not have to admit a basis $\wt E$ so that both $\Fo$ and $\wt\Fo$ are $\wt E$-commensurable.

It turns out that, for $G^0=\SU(p,\infty)$ and $\SU(\infty,\infty)$, the requirement for the existence of an open orbit on an ind-variety of the form $\Fl(\Fo,E)$ imposes no restriction on the type of the generalized flag~$\Fo$. More precisely, we have

\corop{If $G^0=\SU(p,\infty)$\textup, $0\leq p<\infty$\textup, then $\Fl$ always\label{coro:SU_open} has an open $G^0$-orbit. If\break $G^0=\SU(\infty,\infty)$ then there exist a basis $\wt E$ of $V$ and a generalized flag $\wt\Fo$ such that $\Fo$ and $\wt\Fo$ are of the same type and $\wt\Fl=\Fl(\wt\Fo,\wt E)$ has an open $G^0$-orbit.}{For $\SU(p,\infty)$, let $n$ be a positive integer such that the positive index of $\restr{\omega}{V_n}$ equals $p$. Let $\Go_n\in\Fl_n$ be a flag in $V_n$ consisting of nondegenerate subspaces (i.e. the $G_n^0$-orbit of $\Go_n$ is open in $\Fl_n$). Denote by $g$ a linear operator from $G_n$ such that $g(\Fo_n)=\Go_n$, where, as above, $$\Fo_n=\iota_n^{-1}(\Fo)\in\Fl_n.$$ Then $g(\Fo)$ clearly belongs to $\Fl$ and is nondegenerate. Therefore the $G^0$-orbit of $g(\Fo)$ on $\Fl$ is open.

Now consider the case $G^0=\SU(\infty,\infty)$. Let $\wt E$ be an $\omega$-orthogonal basis of $V$. Fix a bijection $E\to\wt E$. This bijection defines an automorphism of $V$. Denote by $\wt\Fo$ the generalized flag consisting of the images of the subspaces from $\Fo$ under this automorphism. Then $\wt\Fo$ and $\Fo$ are of the same type, and each space in $\wt\Fo$ is nondegenerate as it is spanned by a subset of $\wt E$. Thus the $G^0$-orbit of the generalized flag $\wt\Fo$ on $\wt\Fl$ is open.}

The situation is different for $G^0=\SL(\infty,\Rp)$. While an ind-grassmannian $\Grr{F}{E}$ has an open orbit if and only if $\dim F<\infty$ or $\codim_VF<\infty$, it is easy to check that an ind-variety of the form $\wt\Fl=\Fl(\wt\Fo,\wt E)$, where $\wt\Fo$ has the same type as the flag $\Fo$ from Example~\ref{exam:gen_flags} ii) or iii), cannot have an open orbit as long as the basis $\wt E$ satisfies $\tau(\wt e)=\wt e$ for all $\wt e\in\wt E$. For a discussion of the quaternionic case see \cite[Section 5]{IgnatyevPenkovWolf1}.

We now turn our attention to closed orbits. The description of closed orbits is based on the same idea but (similarly to the case for open orbits) the details differ for the various real forms.

Suppose $G^0=\SU(\infty,\infty)$ or $G^0=\SU(p,\infty)$. We call a generalized
flag $\Go$ in $\Fl$ \emph{pseudo-isotropic} if the space $F\cap H^{\perp,V}$ is not
properly contained in $\wt F\cap\wt H^{\perp,V}$ for all $F,~H\in\Go$ and
all $\wt\Go\in\Fl$, where $\wt F,~\wt H$ are the subspaces in $\wt\Go$
corresponding to $F,~H$ respectively. A similar definition can be given
for flags in $\Fl_n$. An isotropic generalized flag is always pseudo-isotropic,
but the converse does not hold. In the particular case when the generalized flag $\Go$ is of the form
$\{\{0\}\subset F\subset V\}$, $\Go$ is pseudo-isotropic if and only if the kernel $\Ker{\restr{\omega}{F}}$
of the form $\restr{\omega}{F}$ is not properly contained in any kernel $\Ker{\restr{\omega}{\wt F}}$ for an $E$-commensurable flag $\{\{0\}\subset\wt F\subset V\}$. Next, suppose $G^0=\SL(\infty,\Rp)$. A generalized flag $\Go$ in $\Fl$ is \emph{real} if $\tau(F)=F$ for all $F\in\Go$. This condition turns out to be equivalent to the following condition: $F\cap\tau(H)$ is not properly contained in $\wt F\cap\tau(\wt H)$ for all $F,~H\in\Go$ and all $\wt\Go\in\Fl$, where $\wt F,~\wt H$ are the subspaces in $\wt\Go$ corresponding to $F,~H$ respectively. Finally, suppose $G^0=\SL(\infty,\Hp)$. We call a generalized flag $\Go$ in $\Fl$ \emph{pseudo-quaternionic} if $F\cap J(H)$ is not properly contained in $\wt F\cap J(\wt H)$ for all $F,~H\in\Go$ and all $\wt\Go\in\Fl$, where $\wt F,~\wt H$ are the subspaces in $\wt\Go$ corresponding to $F,~H$ respectively. If $\Go$ is \emph{quaternionic}, i.e., if $J(F)=F$ for each $F\in\Go$, then $\Go$ is clearly pseudo-quaternionic, but the converse does not hold. If the generalized flag $\Go$ is of the form $\{\{0\}\subset F\subset V\}$, then $\Go$ is pseudo-quaternionic if and only if $\codim_F(F\cap J(F))\leq1$.

\propp{The $G^0$-orbit\label{prop:closed_orbit} $\Omega$ of $\Go\in\Fl$ is closed if and only if
\begin{equation*}
\begin{split}
&\text{$\Go$ is pseudo-isotropic for $G^0=\SU(\infty,\infty)$ and $\SU(p,\infty)$\textup;}\\
&\text{$\Go$ is real for $G^0=\SL(\infty,\Rp)$\textup;}\\
&\text{$\Go$ is pseudo-quaternionic for $G^0=\SL(\infty,\Hp)$.}\\
\end{split}
\end{equation*}}
{First consider the finite-dimensional case, where there is a unique closed $G_n^0$-orbit on $\Fl_n$ (see Theorem~\ref{theo:finite_dim_case}). For all real forms the conditions of the proposition applied to finite-dimensional flags in $V_n$ are easily checked to be closed conditions on points of $\Fl_n$. Therefore, the $G_n^0$-orbit of a flag in $V_n$ is closed if and only if this flag satisfies the conditions of the proposition at the finite level.

Let $G^0=\SU(\infty,\infty)$ or $\SU(p,\infty)$ (for other cases, see \cite[Proposition 5.6]{IgnatyevPenkovWolf1}). Suppose $\Omega$ is closed, so $\Omega_n$ is closed for each $n\geq n_{\Go}$. Assume $\Go$ is not pseudo-isotropic. Then there exist $\wt\Go\in\Fl$ and $A,~B\in\Go$ such that $\wt A\cap\wt B^{\perp,V}\supsetneq A\cap B^{\perp,V}$, where $\wt A,~\wt B$ are the subspaces in $\wt\Go$ corresponding to $A,~B$ respectively. Let $v\in (\wt A\cap\wt B^{\perp,V})\setminus(A\cap B^{\perp,V})$, and $n\geq n_{\Go}$ be such that $v\in V_n$. Then $$v\in(\wt A_n\cap\wt B_n^{\perp,V_n})\setminus(A_n\cap B_n^{\perp,V_n}),$$ where $A_n=A\cap V_n$, $B_n=B\cap V_n$, $\wt A_n=\wt A\cap V_n$, $\wt B_n=\wt B\cap V_n$, because $B^{\perp,V}\cap V_n=B_n^{\perp,V_n}$. This means that $A_n\cap B_n^{\perp,V_n}$ is properly contained in $\wt A_n\cap\wt B_n^{\perp,V_n}$. Hence $\Go_n$ is not pseudo-isotropic, which contradicts the condition that $\Omega_n$ is closed.

Now, assume that $\Omega_n$ is not closed for some $n\geq n_{\Go}$. Then there exist $A_n,~B_n\in\Go_n=\iota_n^{-1}(\Go)$ and $\wt\Go_n\in\Fl_n$ such that $A_n\cap B_n^{\perp,V_n}$ is properly contained in $\wt A_n\cap\wt B_n^{\perp,V_n}$, where $\wt A_n,~\wt B_n$ are the subspaces in~$\wt\Go_n$ corresponding to $A_n$, $B_n$ respectively. Since each $e\in E_{n+1}\setminus E_n$ is orthogonal to $V_n$, $A_{n+1}\cap B_{n+1}^{\perp,V_{n+1}}$ is properly contained in $\wt A_{n+1}\cap\wt B_{n+1}^{\perp,V_{n+1}}$ where $A_{n+1}$, $B_{n+1}$, $\wt A_{n+1}$, $\wt B_{n+1}$ are the respective images of $A_n$, $B_n$, $\wt A_n$, $\wt B_n$ under the embedding $\Fl_n\hookrightarrow\Fl_{n+1}$. Repeating this procedure we see that $\Go$ is not pseudo-isotropic. The result follows.}

\corop{If $G^0=\SL(\infty,\Rp)$\textup, then $\Fl$ always\label{coro:SU_closed} has a closed orbit.}{The $G^0$-orbit of the generalized flag $\Fo$ is closed because the condition $\tau(e)=e$ is satisfied for all basis vectors $e\in E$.}

For $G^0=\SU(p,\infty)$, $0\leq p<\infty$, $G^0=\SU(\infty,\infty)$ or $\SL(\infty,\Hp)$, $\Fl$ may or may not have a closed orbit.

Combining our results on the existence of open and closed orbits, we now obtain the following corollary for all other real forms \cite[Corollary 5.8]{IgnatyevPenkovWolf1}.

\mcoro{For a given real form $G^0$ of $G=\SL_{\infty}(\Cp)$\textup, $G^0\neq\SU(p,\infty)$\textup, $0<p<\infty$\textup, an ind-variety of generalized flags $\Fl=\Fl(\Fo,E)$ has both an open and a closed $G^0$-orbits if\textup, and only if\textup, there are only finitely many $G^0$-orbits on $\Fl$.}

\subsection{Further results}\fakesst Extending\label{sst:open_problems_real_orbits} the results of Section~\ref{sect:orbits_real_forms} to the real forms of the ind-groups $\SO_{\infty}(\Cp)$ and $\Sp_{\infty}(\Cp)$ is a natural problem. A further natural problem is to study the $K$-orbits on $G/P$, where $G=\SL_{\infty}(\Cp)$, $\SO_{\infty}(\Cp)$, $\Sp_{\infty}(\Cp)$ and $K$ is a symmetric ind-subgroup, i.e. $K$ equals the fixed points of an involution of $G$. In the finite-dimensional case, $K$-orbits and $G^0$-orbits are in Matsuki duality. In the recent paper \cite{FressPenkov2} it was proved that Matsuki duality holds for $G=\SL_{\infty}(\Cp)$, $\SO_{\infty}(\Cp)$, $\Sp_{\infty}(\Cp)$ in the case when $P=B$ is a splitting Borel subgroup. This is a first step in the realization of this program.

The theory of cycle spaces for ind-groups has recently been initiated by J.A. Wolf in \cite{Wolf3} and is another potential source of open problems.

\medskip\textsc{Mikhail Ignatyev: Samara National Research University, Ak. Pavlova 1, 443011 Samara, Russia}

\emph{E-mail address}: \texttt{mihail.ignatev@gmail.com}

\medskip\textsc{Ivan Penkov: Jacobs University Bremen, Campus Ring 1, 28759 Bremen, Germany}

\emph{E-mail address}: \texttt{i.penkov@jacobs-university.de}

\end{document}